\documentclass[11pt,reqno]{amsart}

\usepackage[T1]{fontenc}
\usepackage[utf8]{inputenc}
\usepackage{lmodern}
\usepackage{microtype}
\usepackage{amsmath,amssymb,amsthm,mathtools}
\usepackage{mathrsfs}
\usepackage{xcolor}
\usepackage{cite}
\usepackage[colorlinks=true,citecolor=blue,linkcolor=blue,urlcolor=blue]{hyperref}
\allowdisplaybreaks[1]
\numberwithin{equation}{section}

\newtheorem{theorem}{Theorem}[section]
\newtheorem{proposition}[theorem]{Proposition}
\newtheorem{lemma}[theorem]{Lemma}
\newtheorem{corollary}[theorem]{Corollary}
\theoremstyle{remark}

\theoremstyle{definition}

\title[Minimal surfaces of density two]{Minimal surfaces with density two at~infinity}
\author[Y. Bi]{Yuchen Bi}
\address[Yuchen Bi]{Mathematical Institute, Department of Pure Mathematics, University of Freiburg, Ernst-Zermelo-Stra{\ss}e 1, D-79104 Freiburg im Breisgau, Germany}
\email{yuchen.bi@math.uni-freiburg.de}
\date{}
\keywords{Minimal surfaces, quadratic area growth, tangent cones at infinity, Alexandrov reflection, catenoid, Scherk surfaces}
\hypersetup{pdftitle={Minimal surfaces with density two at infinity}, pdfauthor={Yuchen Bi},bookmarksnumbered=true}

\begin{document}
\begin{abstract}
We prove that a connected, properly embedded minimal surface without boundary in \(\mathbb R^3\) and with density two at infinity is a catenoid or a singly periodic Scherk surface.
\end{abstract}
\maketitle

\section{Introduction}\label{sec:introduction}

A properly embedded minimal surface without boundary in \(\mathbb R^3\) with density one at infinity is a plane by the monotonicity formula. The catenoid and the singly periodic Scherk surfaces have density two at infinity. We prove the following classification theorem.

\begin{theorem}\label{thm:classification}
Let \(M\subset\mathbb R^3\) be a connected, properly embedded minimal surface without boundary. Suppose that
\begin{equation}\label{eq:density-at-infinity}
\lim_{R\to\infty} \frac{\mathcal H^2(M\cap B_R(0))}{\pi R^2}=2.
\end{equation}
Then \(M\) is a catenoid or a singly periodic Scherk surface.
\end{theorem}

Theorem~\ref{thm:classification} proves the Scherk uniqueness conjecture stated by Meeks in 2003~\cite[Conjecture~10]{Meeks2003}. In the formulation of Meeks and Wolf~\cite[Conjecture~2]{MW}, this conjecture states that a connected, properly immersed minimal surface without boundary in \(\mathbb R^3\) with density at infinity less than three is a plane, a catenoid, or a singly periodic Scherk surface. By \cite[Section~6]{MW}, such a surface is either a plane or has density two at infinity. In the latter case, the surface is embedded by the monotonicity formula.

Meeks and Wolf~\cite{MW} proved the conjecture under the assumption of an infinite symmetry group. Their classification of singly periodic minimal surfaces with four Scherk ends in the quotient removed the genus restrictions in earlier results of Meeks and Rosenberg~\cite{MeeksRosenberg} and Luo~\cite{Luo}. See also \cite{SchoenSymmetry,Collin,Zhou,HHW} for related results on the uniqueness of the catenoid, \cite{LazardHollyMeeks,PerezTraizet} for classifications of periodic minimal surfaces, and P\'erez~\cite{PerezICM} for a survey.

Every tangent cone of \(M\) at infinity is a sum of two planes, counted with multiplicity. If \(M\) has a tangent cone at infinity consisting of two distinct planes, we prove in Sections~\ref{sec:transverse-asymptotics}--\ref{sec:scherk-rigidity} that \(M\) is a singly periodic Scherk surface. Otherwise, every tangent cone at infinity is a plane of multiplicity two, and we prove in Sections~\ref{sec:height-limits}--\ref{sec:catenoid-rigidity} that \(M\) is a catenoid.

Suppose that a tangent cone \(C\) consists of two distinct planes meeting along the \(z\)-axis. In Section~\ref{sec:transverse-asymptotics}, we prove that \(C\) is the unique tangent cone of \(M\) at infinity. Choose \(R_j\to\infty\) such that \(R_j^{-1}M\to C\). On compact subsets away from the axis, \(R_j^{-1}M\) consists of four graphs over the half-planes of \(C\). The graph functions converge smoothly to \(0\) and satisfy the minimal surface equation, whose linearization at \(0\) is the Laplace equation.

Put \(N_z=N\cdot e_3\), where \(N\) is a unit normal. On each graph, interior estimates for the minimal surface equation bound the curvature of \(R_j^{-1}M\cap\{z=t\}\) in terms of the \(L^2\) norm of \(N_z\). If \(M_j\) are translates and dilates of \(M\) with \(N_z\to0\) in \(L^2_{\mathrm{loc}}\), then every varifold limit of \(|M_j|\) is a sum of at most two planes. We use these facts to show that the Laplace equation holds across the axis on each plane. The degree-one homogeneous solutions are therefore linear on each plane and arise from rotations of the two planes of \(C\). As in~\cite{AllardAlmgren,SimonAsymptotics}, the density deficit
\[
e(R)=2-\frac{\mathcal H^2(M\cap B_R)}{\pi R^2}
\]
satisfies \(e(R)=O(R^{-\gamma_0})\) for some \(\gamma_0>0\), which implies uniqueness of \(C\).

In fact, \(e(R)=O(R^{-\gamma})\) for every \(0<\gamma<1\). We argue by contradiction and choose dilations \(M_j=R_j^{-1}M\) such that
\[
\frac{N_z^2\mathcal H^2\llcorner M_j}{e(R_j)}\rightharpoonup\nu\ne0
\]
on \(\mathbb R^3\setminus\{0\}\), with
\[
\nu(B_{2r}\setminus B_r)\le Cr^{2-\gamma}
\]
for \(0<r<1\). The graph functions of \(M_j\), divided by \(\sqrt{e(R_j)}\), have harmonic limits \(U_1,U_2\) on the two planes, with \(U_k(0)=0\) and sublinear growth at infinity. Hence \(U_1=U_2=0\), so \(\nu\) is supported on the axis. Passing to the limit in the first variation formula, we obtain \(\partial_z\nu=0\) on each half-axis and hence \(\nu=c_\pm\,dz\). Since \(\nu(B_{2r}\setminus B_r)=(c_++c_-)r\) and \(\gamma<1\), \(c_+=c_-=0\), contradicting \(\nu\ne0\).

For some \(0<\beta<1/2\), the graph function \(u(r,z)\) of \(M\) over either plane of \(C\) therefore satisfies
\[
|u|+|r||\nabla u|=O(|r|^\beta) \qquad\text{on }|z|\le|r|
\]
as \(|r|\to\infty\). Since \(N_z^2=O(|r|^{2\beta-2})\) and \(2\beta-2<-1\), the integral of \(N_z^2\) over \(M\cap\{|z|\le H\}\) is finite for every \(H>0\). We prove that the limit
\[
c(M)=\lim_{H\to\infty}\frac1{2H}\int_{M\cap\{|z|\le H\}}N_z^2\mathrm d\mathcal H^2
\]
is finite and positive. Moreover, \(c\) is continuous under varifold convergence within the class of surfaces with tangent cone \(C\) at infinity. It is invariant under translations and satisfies \(c(\lambda M)=\lambda c(M)\) for \(\lambda>0\).

For every \(t\in\mathbb R\), there is a point \(p\in M\cap\{z=t\}\) whose normal stays a fixed positive distance from the normals of \(C\). By the positivity and continuity of \(c\), \(M-p\) lies in a compact family of nonflat surfaces with tangent cone \(C\). Thus the graph functions of \(M-p\) are defined for \(|z|\le|r|\) and \(|r|\ge R_0\), with \(R_0\) independent of \(p\). The graph function \(u\) of \(M\) is then defined on
\[
\{(r,z):|r-a(z)|>R_0\}
\]
for some \(a:\mathbb R\to\mathbb R\). With \(u\) extended linearly in \(r\) across \(|r-a(z)|\le R_0\), the identity \mbox{\(\Delta_M(X-ze_3)=0\)} implies
\[
\Delta u=\operatorname{div}\mathcal F,
\]
where \(\mathcal F=(\mathcal F^r,\mathcal F^z)\) satisfies
\[
\|\mathcal F\|_{L^\infty(\mathbb R^2)}+\sup_z\int_{\mathbb R}|\mathcal F(r,z)|\mathrm dr\le C.
\]
By \(\Delta_M z=0\),
\[
\left|\int_I\int_{\mathbb R}\mathcal F^z(r,z)\mathrm dr\mathrm dz\right|\le C(1+|I|)^{2\beta}
\]
for every interval \(I\). Since \(2\beta<1\), an integral representation of \(u\) in terms of \(\mathcal F\) gives \(|u|\le C\). The four graphs cover \(M\) outside a fixed cylinder about the \(z\)-axis. After translating \(M\), the same representation implies
\[
\sup_z|u(R,z)+u(-R,z)|\le CR^{2\beta-1}
\]
for both planes and all sufficiently large \(R\). By Alexandrov reflection, \(M\) is symmetric about the two bisector planes.

With the reflection planes \(x=0\) and \(y=0\), \(M\cap\{x>0,y>0\}\) admits a conformal parametrization \(X\) on \(\mathbb H_+=\{w=z+is:s>0\}\) satisfying
\[
dX=\operatorname{Re}\left[(-\sinh f,-i\cosh f,1)\,dw\right],
\]
where \(f\) is holomorphic, \(0<\operatorname{Im}f<\pi/2\), and \(e^f\) is the stereographic projection of the Gauss map.

For a singly periodic Scherk surface, the parameter \(w\) can be chosen so that
\[
e^{2f(w)}=-\frac{\sin(\pi w/P)}{\sin(\pi(w-a)/P)}
\]
for some \(P>0\) and \(a\in(0,P)\). Then \(f^{-1}(f(w))=w+P\mathbb Z\), so points with the same value of \(f\) have the same \(s\)-coordinate.

For a regular value \(\zeta\), write \(f^{-1}(\zeta)=\{Z_j=x_j+iy_j\}\). For bounded \(h=(h_j)\) with \(h_j\in\mathbb R\), define the weighted discrete Laplacian \(L_\zeta\) by
\[
(L_\zeta h)_j=2\sum_{k\ne j}\frac{h_j-h_k}{|Z_j-Z_k|^2}.
\]
For finitely supported \(h\),
\[
\sum_j h_j(L_\zeta h)_j=\sum_{j\ne k}\frac{(h_j-h_k)^2}{|Z_j-Z_k|^2}.
\]
The surface \(M\) is a singly periodic Scherk surface if and only if \(E(f)=0\), where
\[
E(f)=2\int_{0<\operatorname{Im}\zeta<\pi/2}\cosh(2\operatorname{Re}\zeta)\sum_{j\ne k}\frac{(y_j-y_k)^2}{|Z_j-Z_k|^2}\,\mathrm dA(\zeta).
\]

We prove \(E(f)=0\) by estimating
\[
W_\psi=\int_{0<\operatorname{Im}\zeta<\pi/2}\cosh(2\operatorname{Re}\zeta)\sum_j\psi(x_j)y_j(L_\zeta y)_j\,\mathrm dA(\zeta),
\]
for \(\psi\in C_c^\infty(\mathbb R;\mathbb R)\), with \(y=(y_j)=(s(Z_j))\).
Here
\[
(L_\zeta y)_j=-\operatorname{Im}\mathcal V_f(Z_j)
\]
for every regular value \(\zeta\), where \(\mathcal V_f\) is meromorphic on \(\mathbb H_+\) and satisfies \(\operatorname{Im}\mathcal V_f=0\) on \(\mathbb R\). Put \(X_h=X-ze_3\) and \(\mathcal P=|X_h|^2-s^2\). The product \(\mathcal V_f\mathcal P_{ww}\) is smooth on \(\mathbb H_+\) and satisfies
\[
\partial_{\bar w}^{\,2}(\mathcal V_f\mathcal P_{ww})
=\mathcal V_f\,\partial_{\bar w}^{\,2}\mathcal P_{ww}
=\frac12\mathcal V_f\cosh(2\operatorname{Re}f)|f'|^2.
\]
Integrating by parts, we obtain
\[
\begin{aligned}
W_\psi={}&-\operatorname{Re}\int_{\mathbb H_+}\psi'(z)\mathcal V_f\mathcal P_{ww}\,\mathrm dz\,\mathrm ds\\
&-\frac12\int_0^\infty s\left(\int_{\mathbb R}\psi''(z)\operatorname{Im}(\mathcal V_f\mathcal P_{ww})\,\mathrm dz\right)\mathrm ds.
\end{aligned}
\]

We now choose \(\psi_R\in C_c^\infty(\mathbb R)\) with \(0\le\psi_R\le1\) and \(\psi_R\to1\) on compact sets so that
\[
|W_{\psi_R}|\le C,\qquad E(f)\le2\liminf_{R\to\infty}W_{\psi_R}.
\]
Moreover, \(W_{\psi_R}\to0\) if \(\mathcal V_{f_*}'=0\) for every locally uniform limit \(f_*=\lim_{n\to\infty}f(\cdot+h_n)\) on \(\mathbb H_+\) with \(h_n\in\mathbb R\) and \(|h_n|\to\infty\). Since \(E(f)<\infty\), every such limit satisfies \(E(f_*)=0\) and \(\mathcal V_{f_*}'=0\). It follows that \(W_{\psi_R}\to0\), hence \(E(f)=0\) and \(M\) is a singly periodic Scherk surface.

Suppose instead that every tangent cone at infinity is a plane of multiplicity two. Let \(R_j\to\infty\) and put \(\Sigma_j=R_j^{-1}M\). In Section~\ref{sec:height-limits} we choose coordinates \(X=(x,h,z)\) so that \(\{h=0\}\) minimizes the height excess of \(\Sigma_j\) among planes through the origin. We approximate \(\Sigma_j\) by the graphs of two Lipschitz functions \(u_{j,1}\le u_{j,2}\)~\cite{BDD} and put \(a_j=\|h\|_{L^2(\Sigma_j\cap B_8)}\). After passing to a subsequence, \((u_{j,1},u_{j,2})/a_j\) converges weakly in \(W^{1,2}_{\mathrm{loc}}(\mathbb R^2)\) to \((-v,v)\), where \(v\ge0\) is subharmonic.

For every \(\delta>0\), \(\Sigma_j\cap\{|h|>\delta a_j\}\) consists of smooth graphs on each fixed compact set for all large \(j\). On each side of \(\{h=0\}\), write the smooth graph as \(h=h(x,z)\) and extend \(r_j(x,z)=(|h(x,z)|/a_j-\delta)_+\) by zero to \(\mathbb R^2\). These functions converge weakly to \((v-\delta)_+\) in \(W^{1,2}_{\mathrm{loc}}\). By the minimal surface equation, \(\operatorname{div}(b_j\nabla r_j)\ge0\) with \(b_j\to1\) locally uniformly. The trace-free part of \(\nabla r_j\otimes\nabla r_j\) converges in distributions to that of \(\nabla(v-\delta)_+\otimes\nabla(v-\delta)_+\)~\cite[Proposition~3]{Delort}. For each compact \(K\subset\mathbb R^3\), the equation \(\Delta_{\Sigma_j}h=0\) implies
\[
\int_{\Sigma_j\cap K\cap\{|h|\le\delta a_j\}}|\nabla_{\Sigma_j}h|^2\mathrm d\mathcal H^2\le C_K\delta a_j^2.
\]
Letting \(j\to\infty\) and then \(\delta\downarrow0\) in the first variation formula, we obtain
\[
\operatorname{div}\left(2\nabla v\otimes\nabla v-|\nabla v|^2I\right)=0.
\]
In coordinates \(\xi=(x,z)\), \((v_x-iv_z)^2\) is holomorphic on the whole plane. It is constant since \(v\) has at most linear growth. Thus \(v(\xi)=c+|p\cdot\xi-b|\) with \(c\ge0\). Put \(D_r=\{\xi\in\mathbb R^2:|\xi|<r\}\). The function \(t\mapsto\log\|v\|_{L^2(D_{e^t})}\) is strictly convex unless \(v\) is homogeneous. A three-annulus argument (cf.~\cite[Proposition~3.7 and Corollary~3.8]{DS}; \cite[Propositions~3.9--3.10]{LSS}) then shows that either all limits \(v\) are positive constants or all have the form \(v(\xi)=|p\cdot\xi|\) with \(p\ne0\).

In either case, \(\Delta(v^2)=2|\nabla v|^2\). The identity \(\Delta_{\Sigma_j}h^2=2|\nabla_{\Sigma_j}h|^2\) and the convergence \((u_{j,1},u_{j,2})/a_j\to(-v,v)\) in \(L^2_{\mathrm{loc}}\) imply
\[
\sum_{i=1}^2\int\varphi\left|\nabla\frac{u_{j,i}}{a_j}\right|^2\mathrm d A
\longrightarrow2\int\varphi|\nabla v|^2\mathrm d A
\]
for \(0\le\varphi\in C_c^\infty(\mathbb R^2)\). Hence \((u_{j,1},u_{j,2})/a_j\to(-v,v)\) strongly in \(W^{1,2}_{\mathrm{loc}}\).

If \(v(\xi)=|p\cdot\xi|\) with \(p\ne0\), then \(R^{-1}M\) is approximated by two planes making an angle \(\alpha_R\), where
\[
\alpha_R\to0,\qquad \alpha_R\ge c_\varepsilon R^{-\varepsilon}\quad\text{for every }\varepsilon>0.
\]
In Section~\ref{sec:nonconstant-heights} we rule out this case. Put
\[
F(R)=e(R)+E_h(R)+E_C(R).
\]
Here \(e(R)\), \(E_h(R)\), and \(E_C(R)\) denote the density deficit, height excess, and tilt excess, respectively. We prove
\begin{equation}\label{eq:intro-excess-improvement}
F(R)\le\varepsilon F(\Lambda R)+\frac{C_\varepsilon}{(\Lambda R)^4}\int_{M\cap B_{\Lambda R}}(X\cdot N)^2\mathrm d\mathcal H^2
\end{equation}
for some fixed \(\Lambda>1\), every \(\varepsilon>0\) and all sufficiently large \(R\). The constant \(C_\varepsilon\) is independent of the angle between the planes. By \eqref{eq:intro-excess-improvement} and the monotonicity formula, \(\alpha_R\le CR^{-\gamma}\) for some \(\gamma>0\), contrary to the lower bound \(c_\varepsilon R^{-\varepsilon}\).

If \eqref{eq:intro-excess-improvement} fails for some \(\varepsilon>0\), there are \(R_j\to\infty\) such that
\[
F(\Lambda R_j)=O(\beta_j^2),\qquad
\int_{\Sigma_j\cap B_1}(X\cdot N)^2\mathrm d\mathcal H^2=o(\beta_j^2),
\]
where \(\beta_j^2=F(R_j)>0\) and \(\Sigma_j=(\Lambda R_j)^{-1}M\).
Choose a pair of planes attaining \(E_h(\Lambda R_j)\). After rotating coordinates, their union consists of the graphs of
\[
f_{j,1}(x,z)=m_j-k_j|x-d_j|,\qquad f_{j,2}(x,z)=m_j+k_j|x-d_j|,
\]
where \(k_j>0\) and \(k_j,m_j,d_j\to0\). Let \(u_{j,1}\le u_{j,2}\) be Lipschitz functions on \(D_1\) whose graphs approximate \(\Sigma_j\). The strong convergence \((u_{j,1},u_{j,2})/a_j\to(-v,v)\) gives \(\beta_j=o(k_j)\).
After passing to a subsequence,
\[
\frac{(u_{j,1}-f_{j,1},\,u_{j,2}-f_{j,2})}{\beta_j}
\rightharpoonup
\bigl(b-\operatorname{sgn}(x)\eta,\,b+\operatorname{sgn}(x)\eta\bigr)
\quad\text{in }L^2_{\mathrm{loc}}(D_1),
\]
where \(b\) is harmonic and \(\eta_x-i\eta_z\) extends holomorphically across \(\{x=0\}\). By the bound on \(X\cdot N\), \(b-b(0)\) and \(\eta-c\) are homogeneous of degree one for some constant \(c\). Thus \(b\) and \(\eta\) are affine. Since the pair of planes attains \(E_h(\Lambda R_j)\), we have \(b=\eta=0\). The functions \((u_{j,i}-f_{j,i})/\beta_j\) converge strongly to zero in \mbox{\(W^{1,2}_{\mathrm{loc}}(D_1\setminus\{x=0\})\)}. Thus the measures \(\beta_j^{-2}N_z^2\mu_j\), where \(\mu_j=\mathcal H^2\llcorner\Sigma_j\), can concentrate only on the \(z\)-axis.

To exclude concentration on the axis, we compare \(\Sigma_j\) with the pair \mbox{\(C_j=|\{h=k_jx\}|+|\{h=-k_jx\}|\)} through the origin. By the bound on \(X\cdot N\), \(f_{j,1}=-k_j|x|+o(\beta_j)\) and \(f_{j,2}=k_j|x|+o(\beta_j)\). For \(0\le\phi\in C_c^\infty(B_1)\) radial, constant near \(0\), and nonincreasing in \(|X|\), the first variation identities for \(\Sigma_j\) and \(C_j\) give
\begin{equation}\label{eq:intro-axial-integral}
0\le\frac13\int\phi N_z^2\mathrm d\mu_j=\int\phi\mathrm d(\mu_j-\|C_j\|)+o(\beta_j^2)\le o(\beta_j^2),
\end{equation}
where the last inequality follows from \(\mu_j(B_r)\le2\pi r^2=\|C_j\|(B_r)\). Thus \(\int_{B_r}N_z^2\mathrm d\mu_j=o(\beta_j^2)\) for \(0<r<1\). Together with the bound on \(X\cdot N\), this excludes concentration of \(\beta_j^{-2}|u_{j,i}-f_{j,i}|^2\mathrm d A\) on \(\{x=0\}\), so \((u_{j,i}-f_{j,i})/\beta_j\to0\) in \(L^2_{\mathrm{loc}}(D_1)\) and \(E_h(R_j)=o(\beta_j^2)\). Equation~\eqref{eq:intro-axial-integral} and the monotonicity formula also imply \(e(R_j)=o(\beta_j^2)\).

Put \(g_j=(u_{j,2}-u_{j,1})/2\) and \(\zeta_j=|\partial_xg_j|-k_j\). The tilt excess satisfies
\[
E_C(R_j)\le C\int_{D_{1/2}}\zeta_j^2\mathrm d A+o(\beta_j^2).
\]
To show that \(\zeta_j=o(\beta_j)\) in \(L^2_{\mathrm{loc}}(D_1)\), set
\[
w_j=\frac{(\partial_xg_j)^2-k_j^2}{2}=k_j\zeta_j+\frac12\zeta_j^2\ge-\frac12k_j^2.
\]
By the first variation formula, \(w_j\) satisfies
\[
L_jw_j=\operatorname{div}\operatorname{div}\mathcal R_j,\qquad
|\mathcal R_j|(D)=o(\beta_j^2)
\]
for every \(D\Subset D_1\), where \(L_j=\partial_x^2+(1+k_j^2)\partial_z^2\) and \(\mathcal R_j\) is a matrix of Radon measures.

Fix \(D\Subset\Omega\Subset D_1\). On \(\{w_j\le k_j^2\}\), \(\zeta_j^2\le4w_j^2/k_j^2\). Since \(\zeta_j/\beta_j\rightharpoonup0\) in \(L^2_{\mathrm{loc}}(D_1)\), the weak \((1,1)\) estimate for \(L_j^{-1}\operatorname{div}\operatorname{div}\) gives
\[
\int_{D\cap\{w_j\le k_j^2\}}\zeta_j^2\mathrm d A=o(\beta_j^2).
\]

For a fixed sufficiently small \(\delta>0\), we can choose for each \(j\) a collection \(\mathcal D_j\) of disks with uniformly bounded overlap covering \(D\cap\{w_j>k_j^2\}\), except for a set of \(\zeta_j^2\,\mathrm d A\)-measure \(o(\beta_j^2)\). Each \(D_r(y)\in\mathcal D_j\) satisfies
\[
D_{2r}(y)\subset\Omega,\qquad
\int_{D_r(y)}\zeta_j^2\mathrm d A=\delta k_j^2r^2.
\]
For any choice \(D_{r_j}(y_j)\in\mathcal D_j\), after passing to a subsequence \(k_j^{-2}w_j(y_j+r_j\,\cdot)\) converges in \(L^1(D_1)\) to a harmonic function \(w_\infty\) on \(D_2\) satisfying
\[
\|w_\infty\|_{L^\infty(D_1)}\le C(\sqrt\delta+\delta)<1.
\]
Rescaling back and using \(\zeta_j^2\le2w_j\) on \(\{w_j>k_j^2\}\), we obtain
\[
\int_{D_r(y)\cap\{w_j>k_j^2\}}\zeta_j^2\mathrm d A
=o(1)\int_{D_r(y)}\zeta_j^2\mathrm d A
\]
uniformly for \(D_r(y)\in\mathcal D_j\). Summing and using bounded overlap gives
\[
\int_{D\cap\{w_j>k_j^2\}}\zeta_j^2\mathrm d A=o(\beta_j^2).
\]
Thus \(\zeta_j/\beta_j\to0\) strongly in \(L^2_{\mathrm{loc}}(D_1)\). Hence \(E_C(R_j)=o(\beta_j^2)\) and \(F(R_j)=o(\beta_j^2)\), contrary to \(\beta_j^2=F(R_j)>0\).

Thus \(v\) is constant. In this case, there is a plane \(P_\infty\) such that \(\sup_{X\in M\cap B_R}\operatorname{dist}(X,P_\infty)=O(R^\varepsilon)\) for every \(\varepsilon>0\). In Section~\ref{sec:catenoid-rigidity}, we prove that \(M\) is symmetric about a plane parallel to \(P_\infty\). After translating \(M\), take coordinates \(X=(x,z)\in\mathbb R^2\times\mathbb R\) with \(P_\infty=\{z=0\}\) and put
\[
L(R)=\sup_{\substack{(x,s),(x,t)\in M\\s\ne t,\ |x|\le R}}\frac{|s+t|}{2}.
\]
Suppose \(L\not\equiv0\). For every sufficiently large \(R\) there is \(r\in(R,2R)\) such that \(M\cap(\partial D_r\times\mathbb R)\) consists of two smooth graphs \(u_+>u_-\) with \mbox{\(\sup_{\partial D_r}|u_++u_-|=o(L(8R))\)}. By Schoen's reflection argument~\cite[\mbox{Theorem~1}]{SchoenSymmetry},
\[
L(R)\le\frac12\sup_{\partial D_r}|u_++u_-|=o(L(8R)).
\]
Thus \(L(8R)\ge8L(R)\) for all sufficiently large \(R\), so \(L(R)\ge cR\), contrary to \(L(R)=O(R^\varepsilon)\) for \(0<\varepsilon<1\). Thus \(L\equiv0\) and \(M\) consists of two graphs interchanged by reflection in \(\{z=0\}\).

\begin{samepage}
For \(b>0\), the graph
\[
z=b\operatorname{arcosh}(|x|/b)
\]
over \(\{|x|>b\}\) is the upper half of a catenoid. Under the change of variables \(y=(1+\sqrt{1-b^2/|x|^2})x\), we have \mbox{\(z=U(y)=b\log(|y|/b)\)} on \(\{|y|>b\}\).
\end{samepage}

Following Traizet~\cite[Theorem~10]{TraizetRoof}, by a change of variables we write \(z=U(y)\) on \(M\cap\{z>0\}\), where \(U\) is a positive harmonic function on \(\Omega'\subset\mathbb R^2\), with \(U=0\) and \(|\nabla U|=1\) on \(\partial\Omega'\). It also satisfies \(U(y)=O(|y|^\varepsilon)\) for every \(\varepsilon>0\). With \(U=0\) on \(\mathbb R^2\setminus\Omega'\), we have \(\Delta U=\mathcal H^1\llcorner\partial\Omega'\), and hence \(\mathcal H^1(\partial\Omega'\cap D_R)=O(R^\varepsilon)\). For \(x\in\partial\Omega'\),
\begin{equation}\label{eq:intro-boundary-integral}
\operatorname{p.v.}\int_{\partial\Omega'}\frac{x-y}{|x-y|^2}\,\mathrm d\mathcal H^1(y)=\pi \mathbf{n}(x),
\end{equation}
where \(\mathbf{n}\) is the inward unit normal to \(\Omega'\). If \(K=\mathbb R^2\setminus\Omega'\) were unbounded, the coarea formula and the boundary length bound would imply that a line in \(\Omega'\) separates two components of \(K\). Take the separating line to be \(\{x_1=0\}\). By \eqref{eq:intro-boundary-integral}, the antisymmetry of the kernel and the boundary length bound, we obtain
\[
0<\iint_{\substack{x,y\in\partial K\\x_1>0>y_1}}
\frac{x_1-y_1}{|x-y|^2}\,\mathrm d\mathcal H^1(x)\,\mathrm d\mathcal H^1(y)\le0.
\]
Thus \(K\) is bounded. By \eqref{eq:intro-boundary-integral},
\[
2\pi|K|=\pi\int_{\partial K}x\cdot \mathbf{n}(x)\,\mathrm d\mathcal H^1(x)=\frac12\mathcal H^1(\partial K)^2.
\]
Equality in the isoperimetric inequality implies that \(K\) is a disk, so \(M\) is a catenoid.

\subsection*{AI disclosure}
A publicly available large language model (LLM) was used to develop mathematical arguments for this paper. The LLM was also used to draft and revise the manuscript. The author takes full responsibility for the paper.

\section{Area bounds and varifold limits}\label{sec:varifold-limits}

Throughout the paper, \(M\) satisfies the hypotheses of Theorem~\ref{thm:classification}. Write \(B_R=B_R(0)\), \(\mu=\mathcal H^2\llcorner M\), and let \(|M|\) be the associated integral varifold. A tangent cone of \(M\) at infinity is a varifold limit of \(|R_j^{-1}M|\) on \(\mathbb R^3\), where \(R_j\to\infty\).

The monotonicity formula \cite[Chapter~4, Section~3]{SimonGMT} gives
\begin{equation}\label{eq:area-bound}
\mu(B_r(p))\le2\pi r^2,
\end{equation}
for \(p\in\mathbb R^3\) and \(r>0\), since a fixed translation does not change the density at infinity.

For a varifold \(V\), write \(\mu_V=\|V\|\). Write \(\mu_j=\|V_j\|\) for varifolds and \(\mu_j=\mathcal H^2\llcorner M_j\) for surfaces. For \(p\in\mathbb R^3\) and \(r>0\), put \(\eta_{p,r}(X)=(X-p)/r\) and write \(V-p=(\eta_{p,1})_\#V\).

We recall the following facts about varifold convergence; see~\cite{Allard} and~\cite[Chapter~8, Theorem~5.8 and Remark~5.9]{SimonGMT}.

\begin{lemma}\label{lem:varifold-limits}
Let \(V_j\to V\) as varifolds in an open set \(U\subset\mathbb R^3\), where the \(V_j\) are stationary integral \(2\)-varifolds. Then \(V\) is stationary and integral. The supports converge locally in the Hausdorff topology, and the convergence is smooth near every point of density one of \(V\).
\end{lemma}

\begin{lemma}\label{lem:cones}
Let \(V\) be a varifold limit in \(\mathbb R^3\) of properly embedded minimal surfaces without boundary and with density at infinity at most two. Then \(V\) is a sum of at most two affine planes, counted with multiplicity, or \(V=|S|\) for a connected, properly embedded minimal surface \(S\) of density two at infinity. If \(V\) is a cone, it is a sum of planes through the origin.
\end{lemma}
\begin{proof}
By the monotonicity formula, \(\|V\|(B_r(p))\le2\pi r^2\), and \(V-p\) is a cone whenever \(\Theta(V,p)=2\).

Suppose first that \(V\ne0\) is a cone. A tangent cone \(W\) at \(x\in\operatorname{spt}V\setminus\{0\}\) is invariant under translations parallel to \(x\). By stationarity,
\[
W=\sum_{i=1}^n|\mathbb R_+\xi_i+\mathbb R x|, \qquad \sum_{i=1}^n\xi_i=0, \qquad \Theta(V,x)=\frac n2\le2,
\]
where \(\xi_i\in x^\perp\cap S^2\). If \(n=3\), the half-planes are distinct. The cone \(W\) is also a varifold limit of properly embedded minimal surfaces \(T_j\) without boundary. By Lemma~\ref{lem:varifold-limits}, the circle \(\Gamma=x^\perp\cap\partial B_1\) is transverse to \(T_j\) and \(\#(\Gamma\cap T_j)=3\) for large \(j\). Since \(\Gamma\) is closed, \(\#(\Gamma\cap T_j)\) is even, a contradiction.

Thus \(n\in\{2,4\}\); see also~\cite{White}. After relabeling,
\[
\xi_{2i}=-\xi_{2i-1}, \qquad i=1,\ldots,n/2.
\]
If \(n=4\) for some \(x\), then \(\Theta(V,x)=2\) and \(V-x=W\). Otherwise \(V\) is smooth away from \(0\) by Allard's regularity theorem. Its second fundamental form \(A\) satisfies
\[
A(X,\cdot)=0, \qquad \operatorname{tr}A=0.
\]
Hence \(A=0\) and \(V\) is a plane.

For general \(V\), suppose that \(V\) is not of the form \(\sum_{i=1}^m|P_i|\), with \(0\le m\le2\) and affine planes \(P_i\). Since \(\Theta(V,p)\in\{1,2\}\) and \(V-p\) is not a cone,
\[
\Theta(V,p)=1\qquad\text{for every }p\in\operatorname{spt}V.
\]
By Allard's regularity theorem, \(V=|S|\) for a properly embedded minimal surface \(S\). Let \(S_0\) be a connected component of \(S\) that is not a plane. Then
\[
\lim_{R\to\infty}\frac{\mathcal H^2(S_0\cap B_R)}{\pi R^2}=2.
\]
Since \(\mathcal H^2(S\cap B_R)\le2\pi R^2\) for every \(R>0\), we have \(S=S_0\) by the monotonicity formula.
\end{proof}

\begin{lemma}\label{lem:constant-direction}
Let \(V\) satisfy the hypotheses of Lemma~\ref{lem:cones}. If \(e_3\in T_xV\) for \(\|V\|\)-almost every \(x\), then \(V\) is a sum of at most two affine planes parallel to \(e_3\), counted with multiplicity.
\end{lemma}
\begin{proof}
Suppose \(V\) is not a sum of planes and write \(V=|S|\) by Lemma~\ref{lem:cones}. Since \(e_3\) is tangent to \(S\) everywhere, its second fundamental form satisfies \(A(e_3,\cdot)=0=\operatorname{tr}A\). Hence \(A=0\), a contradiction.
\end{proof}

For \(R>0\), recall that
\begin{equation}\label{eq:deficit}
e(R)=2-\frac{\mu(B_R)}{\pi R^2}.
\end{equation}
Then \(e(R)\downarrow0\), and for \(0<r<s\),
\begin{equation}\label{eq:monotonicity}
e(r)-e(s)=\frac1\pi\int_{M\cap(B_s\setminus B_r)} \frac{|X\cdot N|^2}{|X|^4} \mathrm d\mu,
\end{equation}
where \(N\) is a unit normal to \(M\).

\section{Transverse tangent cones}\label{sec:transverse-asymptotics}

\begin{theorem}\label{thm:transverse-scherk-classification}
Suppose that a tangent cone of \(M\) at infinity is \(C=|P_1|+|P_2|\), with \(P_1\ne P_2\). Then \(M\) is a singly periodic Scherk surface.
\end{theorem}

Throughout Sections~\ref{sec:transverse-asymptotics}--\ref{sec:scherk-rigidity}, assume that \(C\) is a tangent cone of \(M\) at infinity, and take \(e_3\) parallel to \(P_1\cap P_2\). Put \(z=X\cdot e_3\), \(N_z=N\cdot e_3\), and \(X_h=X-ze_3\).

\subsection{Graphs and height estimates}\label{sec:transverse-graphs}

\begin{lemma}\label{lem:intermediate}
Let \(V_j\) satisfy the hypotheses of Lemma~\ref{lem:cones}. Fix \(0<a<b\). Suppose \(s_j\downarrow0\), \(|t_j|\) stays in a compact subinterval of \((a,b)\), and
\[
\delta_j:=\frac{\mu_j(B_b)}{\pi b^2} -\frac{\mu_j(B_a)}{\pi a^2}=o(s_j^2).
\]
Then every varifold limit of \((\eta_{t_je_3,s_j})_\#V_j\) is a sum of at most two affine planes parallel to \(e_3\), counted with multiplicity.
\end{lemma}
\begin{proof}
For each fixed \(L>0\) and all sufficiently large \(j\), \(a<|X|<b\), \(|z|>a/2\), and \(|X_h|\le Ls_j\) on \(B_{Ls_j}(t_je_3)\). Since \(zN_z=X\cdot N-X_h\cdot N\),
\begin{equation}\label{eq:verticalenergy}
s_j^{-2}\int_{B_{Ls_j}(t_je_3)}N_z^2 \mathrm{d}\mu_j \le C\frac{\delta_j}{s_j^2}+CL^4s_j^2\longrightarrow0.
\end{equation}
For every limit \(W\), \eqref{eq:verticalenergy} implies \(e_3\in T_xW\) for \(\|W\|\)-almost every \(x\). Lemma~\ref{lem:constant-direction} applies.
\end{proof}

Write \(x=(r,z)\) and \(\mathcal{E}=r\partial_r+z\partial_z\).

\begin{lemma}\label{lem:coercive}
Let \(D\subset\mathbb R^2\) be a bounded Lipschitz domain with \(\inf_D|x|>0\). There are \(\varepsilon_0,C>0\), depending only on \(D\), such that every \(a\in L^\infty(D;\mathbb R)\) with \(\|a-1\|_{L^\infty(D)}\le\varepsilon_0\) and every \(v\in H^1(D)\) satisfying \(\operatorname{div}(a\nabla v)=0\) obey
\begin{equation}\label{eq:coercive}
\|\nabla v\|_{L^2(D)}\le C\left(\|v\|_{L^2(D)}+\|\mathcal{E}v-v\|_{L^2(D)}\right).
\end{equation}
\end{lemma}
\begin{proof}
By \cite[Chapter~3, Lemma~7.1]{Necas},
\begin{equation}\label{eq:gradient-negative-norm}
\|\nabla v\|_{L^2(D)}\le C\left(\|v\|_{L^2(D)}+\|\nabla^2 v\|_{H^{-1}(D)}\right).
\end{equation}
For \(x\ne0\), the linear map \(H\mapsto(Hx,\operatorname{tr}H)\) on symmetric \(2\times2\) matrices is invertible. Its inverse has uniformly bounded \(C^1\) coefficients on \(D\).
The identities
\[
(\nabla^2 v)x=\nabla (\mathcal Ev-v), \qquad \operatorname{tr}\nabla^2 v=-\operatorname{div}((a-1)\nabla v)
\]
therefore give
\[
(1-C\varepsilon_0)\|\nabla v\|_{L^2(D)} \le C\left(\|v\|_{L^2(D)}+\|\mathcal Ev-v\|_{L^2(D)}\right).\qedhere
\]
\end{proof}

For fixed \(0<a<b\), set
\[
A_\sigma(a,b)=\{r>\sigma, a<\sqrt{r^2+z^2}<b\}, \qquad 0\le\sigma\le a/2.
\]
In polar coordinates \(r=\rho\cos\theta,z=\rho\sin\theta\), the map
\[
\Phi_\sigma(\rho,\theta)= \left(\rho,(2/\pi)\arccos(\sigma/\rho)\theta\right)
\]
takes \(A_0(a,b)\) onto \(A_\sigma(a,b)\). In the \((r,z)\)-coordinates, \(\Phi_\sigma\) and its inverse have uniformly bounded \(C^2\) norms and converge in \(C^2\) to the identity as \(\sigma\to0\). Thus \(\varepsilon_0,C\) in Lemma~\ref{lem:coercive} are independent of \(\sigma\).

Fix \(r_0>0\), an interval \(I\), and \(\lambda\ge0\), and set
\begin{equation}\label{eq:graphdomain}
\begin{gathered}
\Omega_\rho=\{z\in I, \rho(z)<r<r_0\},\\
0\le\rho\le s_0<r_0/8, \qquad \operatorname{Lip}\rho\le\lambda, \qquad |I|\le\frac{r_0}{8(1+\lambda)}.
\end{gathered}
\end{equation}
Each \(\Omega_\rho\) is star-shaped with respect to every point of \((r_0/2,r_0)\times I\), so \eqref{eq:gradient-negative-norm} holds with \(D=\Omega_\rho\) and \(C\) independent of \(\rho\).

For \(u\in H^1(\Omega_\rho)\), set \(\widetilde u(r,z)=u(\rho(z)+|r-\rho(z)|,z)\) on \((0,r_0)\times I\). Then
\begin{equation}\label{eq:trace}
\|\widetilde u\|_{H^1((0,r_0)\times I)}^2 +\sup_{0<r<r_0}\int_I\widetilde u(r,z)^2 \mathrm{d} z \le C\|u\|_{H^1(\Omega_\rho)}^2.
\end{equation}
In particular, \(\int_{\Omega_\rho\cap\{r<s_0\}}u^2\le Cs_0\|u\|_{H^1(\Omega_\rho)}^2\). If \(u\) satisfies the minimal surface equation, \(\inf_{\Omega_\rho}|x|\ge c_0>0\), and \(\|\nabla u\|_{L^\infty(\Omega_\rho)}+s_0\) is sufficiently small, Lemma~\ref{lem:coercive} and \eqref{eq:trace} give
\begin{equation}\label{eq:absorb}
\|u\|_{H^1(\Omega_\rho)}^2 \le C\left(\int_{\Omega_\rho\cap\{r>s_0\}}u^2 +\|\mathcal{E}u-u\|_{L^2(\Omega_\rho)}^2\right).
\end{equation}
Let \(M_j\) satisfy the hypotheses of Theorem~\ref{thm:classification}. Assume that \(|M_j|\to C=|P_1|+|P_2|\) as varifolds, where \(P_1\ne P_2\). Rotate so that \(P_1\cap P_2=\mathbb{R} e_3\). Let \(\alpha=\angle(P_1,P_2)>0\).

For each \(k\), choose a unit vector \(t_{k,+}\in P_k\cap e_3^\perp\) and a unit normal \(\mathbf{n}_k\), and put \(t_{k,-}=-t_{k,+}\), \(\mathbf{n}_{k,\pm}=\mathbf{n}_k\). Index the four half-planes by \(i=(k,\sigma)\), with \(k\in\{1,2\}\) and \(\sigma\in\{+,-\}\).
Fix \(1<c<d<2\) and write \(A_\sigma=A_\sigma(c,d)\).

By Lemma~\ref{lem:varifold-limits}, \(M_j\to P_1\cup P_2\) smoothly away from \(\mathbb Re_3\). There are \(\sigma_j\downarrow0\) and four smooth functions \(u_{j,i}\) defined on a neighborhood of \(\overline{A_{\sigma_j}}\) such that
\[
X_{j,i}(r,z)=rt_i+ze_3+u_{j,i}(r,z)\mathbf{n}_i\in M_j, \qquad (r,z)\in A_{\sigma_j},
\]
and
\begin{equation}\label{eq:initial-graph-smallness}
\max_i\|\nabla u_{j,i}\|_{L^\infty(A_{\sigma_j})}\to0, \qquad \max_i\sup\frac{|u_{j,i}|}{r}\to0.
\end{equation}
For every fixed annulus \(A'\Subset\{c<|X|<d\}\), these graphs cover \(M_j\cap A'\cap\{|X_h|>2\sigma_j\}\) for all sufficiently large \(j\). Set
\begin{equation}\label{eq:graph-error-density-drop}
H_j=\sum_i\int_{A_{\sigma_j}}u_{j,i}^2, \qquad \delta_j=\frac1\pi\int_{M_j\cap\{1<|X|<2\}} \frac{(X\cdot N)^2}{|X|^4} \mathrm{d}\mu_j,
\end{equation}
and \(a_j=(H_j+\delta_j)^{1/2}>0\). By \eqref{eq:initial-graph-smallness} and varifold convergence, \(H_j,\delta_j\to0\), so \(a_j\to0\).

\begin{lemma}\label{lem:transverse-initial-graphs}
For all sufficiently large \(j\),
\begin{equation}\label{eq:initial-graph-h1-bound}
\sum_i\|u_{j,i}\|_{H^1(A_{\sigma_j})}^2\le C a_j^2.
\end{equation}
\end{lemma}
\begin{proof}
Write \(u=u_{j,i}\), and choose the unit normal \(N\) with \(N\cdot \mathbf{n}_i>0\). Then
\begin{equation}\label{eq:eulerexact}
X\cdot N=\frac{u-ru_r-zu_z}{\sqrt{1+|\nabla u|^2}}.
\end{equation}
For large \(j\), \(X_{j,i}(A_{\sigma_j})\subset\{1<|X|<2\}\), so \(\sum_i\|\mathcal Eu_{j,i}-u_{j,i}\|_{L^2(A_{\sigma_j})}^2\le C\delta_j\). Lemma~\ref{lem:coercive} proves \eqref{eq:initial-graph-h1-bound}.
\end{proof}

Fix \(\lambda=1/16\). Choose \(r_0>0\) and \(I_1\Subset I_2\Subset I_3\) so that
\[
\{|X_h|\le2r_0, z\in\overline{I_3}\} \Subset\{c<|X|<d\}, \qquad |I_3|\le\frac{r_0}{8(1+\lambda)}.
\]
Choose \(0<\eta<1\) so that the four sets \(\{rt_i+ze_3+h\mathbf{n}_i:r>0, z\in\mathbb{R}, |h|<\eta r\}\) are pairwise disjoint and Lemma~\ref{lem:coercive} applies to minimal graphs \(u\) with \(|\nabla u|\le\eta\).

\begin{proposition}\label{prop:transverse-graphs}
There are functions \(\rho_j:I_3\to[8a_j,\infty)\), with \(\operatorname{Lip}\rho_j\le\lambda\) and \(\sup_{I_3}\rho_j\to0\), such that, for large \(j\), \(X_{j,i}\) extend smoothly as maps into \(M_j\) over \(A_{\sigma_j}\cup \Omega_{\rho_j}\), where
\[
\Omega_{\rho_j}=\{z\in I_2, \rho_j(z)<r<r_0\}.
\]
On \(\Omega_{\rho_j}\), \(|u_{j,i}|<\eta r/2\) and \(|\nabla u_{j,i}|<\eta\), and
\[
\begin{aligned}
&M_j\cap\{|X_h|<r_0/2, z\in I_1, |X_h|>2\rho_j(z)\} \\
\subset{}&\bigcup_i X_{j,i}(\Omega_{\rho_j}).
\end{aligned}
\]
The functions \(\widetilde u_{j,i}\) on \((0,r_0)\times I_2\) satisfy
\begin{equation}\label{eq:general-Hone}
\sum_i\|\widetilde u_{j,i}\|_{H^1((0,r_0)\times I_2)}^2 +\sup_{0<r<r_0}\sum_i\int_{I_2}\widetilde u_{j,i}(r,z)^2 \mathrm{d} z\le Ca_j^2.
\end{equation}
Moreover, for \(q\ge Ca_j\),
\begin{equation}\label{eq:general-tail}
|\{z\in I_1:\rho_j(z)>q\}| \le\frac{Ca_j^2}{q^2}.
\end{equation}
\begin{equation}\label{eq:general-rho}
\int_{I_1}\rho_j\le Ca_j, \qquad \sup_{I_1}\rho_j\le Ca_j^{2/3}, \qquad \int_{I_1}\rho_j^3\le Ca_j^{8/3}.
\end{equation}
\end{proposition}
\begin{proof}
For \(t\in I_3\) and \(s>0\), consider four smooth functions \(u_i\) on \((s/8,8s)\times(t-8s,t+8s)\), with \(|u_i|<\eta r/2\) and \(|\nabla u_i|<\eta\), such that
\[
\{h\in(-\eta r,\eta r):rt_i+ze_3+h\mathbf{n}_i\in M_j\}=\{u_i(r,z)\},
\]
and whose graphs cover
\[
M_j\cap\{s/4<|X_h|<4s, |z-t|<4s\}.
\]
By \eqref{eq:initial-graph-smallness}, such functions exist for \(C\sigma_j\le s\le s_*\) and large \(j\), where \(s_*>0\) is fixed and sufficiently small.

Let \(E_j\) consist of the pairs \((t,s)\) with \(s=s_*2^{-m}\ge a_j\), \(m=0,1,\ldots\), and \(t\in s\mathbb Z\cap I_3\) for which no such four graphs exist. Define
\begin{equation}\label{eq:graph-extension-radius}
\rho_j(z)=\max\left(\{8a_j\}\cup \{16s-\lambda|z-t|:(t,s)\in E_j\}\right).
\end{equation}
Then \(\operatorname{Lip}\rho_j\le\lambda\), and
\begin{equation}\label{eq:maxwidth}
\sup\rho_j\le C(\sigma_j+a_j)\longrightarrow0.
\end{equation}
If \((r,z)\in \Omega_{\rho_j}\) and \(r\le s_*\), choose a dyadic scale \(s\le r<2s\) and \(t\in s\mathbb Z\cap I_3\) with \(|t-z|\le s/2\). Then \(s>4a_j\) and \(16s-\lambda|z-t|>r>\rho_j(z)\), so \((t,s)\notin E_j\). For \(r>s_*\), use \(u_{j,i}\) on \(A_{\sigma_j}\). Thus each \(u_{j,i}\) extends smoothly to \(A_{\sigma_j}\cup\Omega_{\rho_j}\).

For \(X\in M_j\) with \(z\in I_1\) and \(2\rho_j(z)<|X_h|<r_0/2\), apply the preceding argument with \(|X_h|\) in place of \(r\). Then \(X=rt_i+ze_3+u_i(r,z)\mathbf n_i\) for some \(i\), with \(r\le|X_h|<2r\). Hence \((r,z)\in\Omega_{\rho_j}\).

With \(s_0=C(\sigma_j+a_j)\), \(\Omega_{\rho_j}\cap\{r>s_0\}\subset A_{\sigma_j}\). On \(\Omega_{\rho_j}\), equation~\eqref{eq:eulerexact} gives \(\sum_i\|\mathcal Eu_{j,i}-u_{j,i}\|_{L^2(\Omega_{\rho_j})}^2\le C\delta_j\). Equations~\eqref{eq:absorb} and \eqref{eq:trace} prove \eqref{eq:general-Hone}.

Consider \((t,s)\in E_j\) for which
\begin{equation}\label{eq:rho-maximizer}
\rho_j(z)=16s-\lambda|z-t|>8a_j \quad \text{ for some }z\in I_1.
\end{equation}
Thus \(16s\le\rho_j(t)\le\rho_j(z)+\lambda|z-t|=16s\). Put \(K=64\). For large \(j\),
\[
R_{t,s}:=\{Ks<r<2Ks, |z-t|<Ks\}\subset \Omega_{\rho_j}.
\]
Indeed, \(s\to0\) uniformly over \(E_j\) and \(\operatorname{dist}(t,I_1)<16s/\lambda\), so \(R_{t,s}\subset(0,r_0)\times I_2\); on this rectangle, \(\rho_j(z)\le16s+\lambda|z-t|<(16+\lambda K)s<Ks\).

For all sufficiently large \(j\) and every pair satisfying \eqref{eq:rho-maximizer},
\begin{equation}\label{eq:height-defect-charge}
s^{-4}\sum_i\int_{R_{t,s}}u_{j,i}^2 \mathrm{d} r \mathrm{d} z +\frac{\delta_j}{s^2}\ge C^{-1}.
\end{equation}
Otherwise, after passing to a subsequence we can choose \((t_j,s_j)\) for which both terms tend to zero. The functions \(s_j^{-1}u_{j,i}(s_jr,t_j+s_jz)\) have Lipschitz constant at most \(\eta\) and tend to zero in \(L^2((K,2K)\times(-K,K))\), hence uniformly. After passing to a further subsequence, let \(|s_j^{-1}(M_j-t_je_3)|\to W\). By Lemma~\ref{lem:varifold-limits}, \(rt_i+ze_3\in\operatorname{spt}W\) for \(K<r<2K\), \(|z|<K\), and each \(i\). Since \(\delta_j/s_j^2\to0\), Lemma~\ref{lem:intermediate}, with \((a,b)=(1,2)\), implies \(W=C\). The convergence is then smooth off \(\mathbb Re_3\), so \((t_j,s_j)\notin E_j\) for large \(j\), a contradiction.

Set \(g_j(z)=\sup_{0<r<r_0}\sum_i\widetilde u_{j,i}(r,z)^2\) on \(I_2\), and zero elsewhere. By Sobolev's inequality in \(r\) and \eqref{eq:general-Hone}, \(\|g_j\|_{L^1(\mathbb R)}\le Ca_j^2\). If \(z\in I_1\) and \(\rho_j(z)>q\ge Ca_j\), choose \((t,s)\) realizing the maximum in \eqref{eq:graph-extension-radius}. Then
\[
s>q/16, \qquad |z-t|<16s/\lambda, \qquad \frac{\delta_j}{s^2}\le\frac{256a_j^2}{q^2}.
\]
For \(q\ge Ca_j\) with \(C\) sufficiently large, \eqref{eq:height-defect-charge} gives
\[
C^{-1}s^4\le\sum_i\int_{R_{t,s}}u_{j,i}^2 \le Ks\int_{t-Ks}^{t+Ks}g_j.
\]
Since \(|z-t|<16s/\lambda\),
\[
\mathcal M g_j(z)\ge\frac1{Cs}\int_{t-Ks}^{t+Ks}g_j\ge C^{-1}s^2\ge C^{-1}q^2,
\]
where \(\mathcal M\) is the maximal operator. Its weak \((1,1)\) estimate proves \eqref{eq:general-tail}.

For \(h=\sup_{I_1}\rho_j\ge Ca_j\), the bound \(\operatorname{Lip}\rho_j\le\lambda\) implies \(|\{z\in I_1:\rho_j(z)>h/2\}|\ge C^{-1}h\). By \eqref{eq:general-tail}, \(h^3\le Ca_j^2\). For \(p=1,3\),
\[
\int_{I_1}\rho_j^p\mathrm dz
=p\int_0^h q^{p-1}|\{z\in I_1:\rho_j(z)>q\}|\mathrm dq.
\]
Together with \eqref{eq:general-tail}, this proves \eqref{eq:general-rho}.
\end{proof}

Fix \(I\Subset I'\Subset I_1\) and set \(Z=\{|X_h|<r_0/2, z\in I'\}\). Let \(T_j\) be the part of \(M_j\cap\{|X_h|<r_0/2, z\in I\}\) outside \(\bigcup_i X_{j,i}(\Omega_{\rho_j})\), and put
\[
F_j=\bigcup_i\{rt_i+ze_3:z\in I, 0<r\le\rho_j(z)\}.
\]

\begin{corollary}\label{cor:transverse-remainder}
For all sufficiently large \(j\),
\begin{equation}\label{eq:coremoments}
\int_{T_j}|X_h|^p \mathrm{d}\mu_j +\int_{F_j}|X_h|^p \mathrm{d}\mu_C \le C\int_{I_1}\rho_j^{p+1}, \qquad p=0,1,2.
\end{equation}
\begin{equation}\label{eq:finite-tilt}
\int_{M_j\cap Z}N_z^2 \mathrm{d}\mu_j\le Ca_j^2.
\end{equation}
\begin{equation}\label{eq:horizontal-height}
\int_{M_j\cap\{|X_h|<r_0/2, z\in I\}}(X_h\cdot N)^2 \mathrm d\mu_j \le C\sum_i\|u_{j,i}\|_{L^2(\Omega_{\rho_j})}^2+Ca_j^{8/3}.
\end{equation}
In particular,
\begin{equation}\label{eq:general-core}
\mu_j(T_j)+\mu_C(F_j)\le Ca_j, \qquad \int_{T_j}|X_h|^2 \mathrm{d}\mu_j=o(a_j^2).
\end{equation}
\end{corollary}
\begin{proof}
For \(z\in I'\) and \(|t-z|<\rho_j(z)\), we have \(t\in I_1\) and \(\rho_j(z)/2\le\rho_j(t)\le2\rho_j(z)\) for large \(j\). Thus, for \(p=0,1,2\),
\[
\begin{aligned}
\int_{M_j\cap Z\cap\{|X_h|\le2\rho_j(z)\}} |X_h|^p \mathrm d\mu_j &\le C\int_{I_1}\rho_j(t)^{p-1} \mu_j\left(B_{C\rho_j(t)}(te_3)\right) \mathrm dt \\
&\le C\int_{I_1}\rho_j^{p+1},
\end{aligned}
\]
and \(\int_{F_j}|X_h|^p \mathrm d\mu_C =\frac4{p+1}\int_I\rho_j^{p+1}\). The cases \(p=0,2\), together with \eqref{eq:general-rho}, give \eqref{eq:general-core}.

The images \(X_{j,i}(\Omega_{\rho_j})\) cover \(M_j\cap Z\cap\{|X_h|>2\rho_j(z)\}\), so \eqref{eq:general-Hone} gives
\[
\int_{M_j\cap Z\cap\{|X_h|>2\rho_j(z)\}}N_z^2\mathrm d\mu_j
\le\sum_i\int_{\Omega_{\rho_j}}|\partial_z u_{j,i}|^2\mathrm dr\mathrm dz\le Ca_j^2.
\]
On \(M_j\cap Z\cap\{|X_h|\le2\rho_j(z)\}\), \(|z|>c\) and \(zN_z=X\cdot N-X_h\cdot N\), so \(N_z^2\le C((X\cdot N)^2+|X_h|^2)\). Together with the bound for \(p=2\), \eqref{eq:general-rho}, and \(\delta_j\le a_j^2\), this proves \eqref{eq:finite-tilt}.

Set \(G(p)=p/\sqrt{1+|p|^2}\) for \(p\in\mathbb R^2\). Choose \(\zeta\in C_c^\infty((-r_0,r_0)\times I_1)\) with \(\zeta=1\) on \((-r_0/2,r_0/2)\times I\). For \(u=u_{j,i}\), multiplying \(\operatorname{div}G(\nabla u)=0\) by \((r-\rho_j(z))^2\zeta^2u\) and integrating by parts gives
\[
\int_{\Omega_{\rho_j}}(r-\rho_j)^2\zeta^2|\nabla u|^2 \le C\int_{\Omega_{\rho_j}}u^2.
\]
On \(r\ge2\rho_j\), \(r^2\le4(r-\rho_j)^2\); on \(\rho_j<r<2\rho_j\), \((u-ru_r)^2\le Cr^2\). Since on a graph
\[
(X_h\cdot N)^2 \mathrm d\mu_j =\frac{(u-ru_r)^2}{\sqrt{1+|\nabla u|^2}} \mathrm dr \mathrm dz,
\]
equations~\eqref{eq:coremoments} and \eqref{eq:general-rho} give \eqref{eq:horizontal-height}.
\end{proof}

\subsection{Harmonic limits across the axis}\label{sec:transverse-harmonic}

In Lemmas~\ref{lem:four-graph-extension} and~\ref{lem:tilt-iteration}, \(V\) satisfies the hypotheses of Lemma~\ref{lem:cones}.

\begin{lemma}\label{lem:four-graph-extension}
For each \(\alpha_0>0\), there are \(C,\eta_0>0\) such that for every \(0<\eta\le\eta_0\) there is \(\delta_0>0\) with the following property. Suppose \(P_1\cap P_2=\mathbb Re_3\), \(\angle(P_1,P_2)\ge\alpha_0\), and \(V\) has a smooth graph \(u_k\) over one half-plane of each \(P_k\), on \(1<r<2, |z|<2\), with
\[
\begin{gathered}
|u_k|\le\eta,\qquad |\nabla u_k|\le\eta, \\
\int_{B_{32}}N_z^2 \mathrm d\|V\|\le\delta_0.
\end{gathered}
\]
Then there are \(p'\in e_3^\perp\), affine planes \(P'_1,P'_2\) with \(P'_1\cap P'_2=p'+\mathbb Re_3\), and graphs \(u'_i\) over the four components of \((P'_1\cup P'_2)\cap\{1/4<|X_h-p'|<4, |z|<4\}\), with
\[
\begin{gathered}
|p'|+\max_k\angle(P'_k,P_k)\le C\eta, \\
|u'_i|\le\frac{\eta}{4},\qquad |\nabla u'_i|\le\frac{\eta}{4}, \\
\operatorname{spt}V\cap\{1/2<|X_h-p'|<2, |z|<2\} \subset\bigcup_i\operatorname{graph}u'_i.
\end{gathered}
\]
\end{lemma}
\begin{proof}
We argue by contradiction and suppose that \(V_j\) satisfy the hypotheses for fixed \(\alpha_0,\eta\), with \(\int_{B_{32}}N_z^2 \mathrm d\|V_j\|\to0\), but admit no four graphs \(u'_i\) as above. After passing to a subsequence, \(P_{j,k}\to P_k\) and \(V_j\to W\), with \(e_3\in T_xW\) almost everywhere in \(B_{32}\). If \(W\) is not a sum of planes, write \(W=|S|\) with \(S\) connected, by Lemma~\ref{lem:cones}. On \(S\cap B_{32}\), \(A(e_3,\cdot)=0=\operatorname{tr}A\), hence \(A=0\). Thus \(S\) is a plane by unique continuation, a contradiction. The limits of the two graphs lie in distinct planes, so \(W=|P'_1|+|P'_2|\), where \(P'_1\cap P'_2=p'+\mathbb Re_3\). For large \(j\),
\[
|p'|+\max_k\angle(P'_k,P_{j,k})\le C\eta<1.
\]
By Lemma~\ref{lem:varifold-limits}, \(V_j\to W\) smoothly off \(P'_1\cap P'_2\).
\end{proof}

Suppose \(V\) has four smooth graphs \(u_i\) on \(h/2<r<2h\), \(|z-t|<2h\), over the half-planes of two planes \(P_1,P_2\) meeting along \(p+\mathbb Re_3\), \(p\in e_3^\perp\), at an angle at least \(\alpha>0\). Write the graphs as \(p+rt_i+ze_3+u_i(r,z)\mathbf{n}_i\), where \(t_i,\mathbf{n}_i\in e_3^\perp\cap S^2\) and \(t_i\perp \mathbf{n}_i\). After translation by \(-te_3\), the tangent line to the section \(z=t\) at \(r=h\) is
\[
\mathcal L_i(t)=p+ht_i+u_i(h,t)\mathbf{n}_i +\mathbb R\left(t_i+u_{i,r}(h,t)\mathbf{n}_i\right).
\]

\begin{lemma}\label{lem:tilt-iteration}
There are constants \(\eta>0\) and \(A_\alpha,C<\infty\), depending only on \(\alpha\), such that, if \(h\ge A_\alpha\kappa>0\),
\[
|u_i|\le\eta h,\qquad |\nabla u_i|\le\eta
\]
for \(h/2<r<2h\) and \(|z-t|<2h\), and
\begin{equation}\label{eq:iteration-tilt}
\int_{B_{64h}(p+te_3)\cap\{|z-t|<a\}}N_z^2 \mathrm{d}\|V\| \le2\kappa a
\end{equation}
for every \(0<a\le64h\), then there is \(p_t\in e_3^\perp\) with
\begin{equation}\label{eq:tangent-concurrency}
\operatorname{dist}(p_t,\mathcal L_i(t))\le C\sqrt{\kappa h}, \qquad |p_t-p|\le Ch.
\end{equation}
There are also \(s\in[A_\alpha\kappa,2A_\alpha\kappa)\) and planes \(P'_1,P'_2\) with \(P'_1\cap P'_2=p_t+\mathbb Re_3\) and \(\angle(P'_1,P'_2)\ge\alpha/2\). Over each half-plane, \(V\) has a smooth graph \(v_i\) on \(s<r<2s\), \(|z-t|<2s\), with
\[
|v_i|\le\eta s,\qquad |\nabla v_i|\le\eta.
\]
\end{lemma}
\begin{proof}
Translate so that \(p=0\). If \(u\) is a graph of \(V\) over \(B_\sigma((r,t))\), with image in \(B_{64h}(te_3)\) and \(|\nabla u|\le1\), \eqref{eq:iteration-tilt} gives \(\|u_z\|_{L^2(B_\sigma((r,t)))}^2\le C\kappa\sigma\). The minimal surface equation and interior estimates for \(u_z\) give on \(B_{\sigma/2}\)
\begin{equation}\label{eq:horizontal-curvature}
\frac{|u_{rr}|}{(1+u_r^2)^{3/2}} \le C|\nabla u_z|\le C\sigma^{-2}\|u_z\|_{L^2(B_\sigma((r,t)))} \le C\sqrt\kappa \sigma^{-3/2}.
\end{equation}

Set \(s_n=2^{-n}h\), and choose \(N\ge0\) with \(A_\alpha\kappa\le s_N<2A_\alpha\kappa\). Put \(p_0=0\), \(P_{0,k}=P_k\). Suppose \(n<N\), \(P_{n,1}\cap P_{n,2}=p_n+\mathbb Re_3\), and \(\angle(P_{n,1},P_{n,2})\ge\alpha/2\). Over their four half-planes, \(V\) has graphs \(u_{n,i}\) on \(s_n<r<2s_n\), \(|z-t|<2s_n\), with
\[
|u_{n,i}|\le\eta s_n,\qquad |\nabla u_{n,i}|\le\eta.
\]
Assume also \(|p_{m+1}-p_m|\le C\eta s_m\) for \(m<n\). Then \(|p_n|\le C\eta\sum_{m<n}s_m\le C\eta h\), so \(B_{32s_n}(p_n+te_3)\subset B_{64h}(te_3)\) and
\begin{equation}\label{eq:local-Morrey-ball}
s_n^{-2}\int_{B_{32s_n}(p_n+te_3)}N_z^2 \mathrm{d}\|V\| \le64\kappa/s_n\le64/A_\alpha\le\delta_0(\alpha/2,\eta).
\end{equation}
Apply Lemma~\ref{lem:four-graph-extension} to \((\eta_{p_n+te_3,s_n})_\#V\). Then \(|p_{n+1}-p_n|\le C\eta s_n\).

In \(\operatorname{spt}V\cap\{z=t\}\), join the point \(r=s_n\) on the graph of \(u_{n,i}\) to the point \(r=s_{n+1}\) on the graph of \(u_{n+1,i}\). The arc has length at most \(Cs_n\). By \eqref{eq:horizontal-curvature} with \(\sigma=cs_{n+1}\), its total curvature is at most \(C\sqrt{\kappa/s_n}\). Hence
\[
\angle(P_{n+1,1},P_{n+1,2}) \ge\alpha-C\eta-C\sum_{m=0}^n\sqrt{\kappa/s_m} \ge\alpha-C\eta-C/\sqrt{A_\alpha}\ge\alpha/2.
\]

Take \(p_t=p_N\) and \(s=s_N\). Since \(\sum_{m=n}^{N-1}s_m\le2s_n\), the total arc length from \(r=s_n\) to the endpoint at \(r=s_N\) is at most \(Cs_n\). The endpoint is within \(Cs_N\) of \(p_N\), so
\[
\operatorname{dist}(p_t,\mathcal L_i(t)) \le Cs_N+C\sum_{n<N}s_n\sqrt{\kappa/s_n} \le C\sqrt{\kappa h}.\qedhere
\]
\end{proof}

Choose \(\eta\) in Proposition~\ref{prop:transverse-graphs} no larger than the constant in Lemma~\ref{lem:tilt-iteration}. For \(k=1,2\) and \((q,z)\in(-r_0/2,r_0/2)\times I_1\), set
\[
\widehat U_{j,k}(q,z)= \begin{cases} \widetilde u_{j,k,+}(q,z),&q>0,\\ \widetilde u_{j,k,-}(-q,z),&q<0. \end{cases}
\]

\begin{proposition}\label{prop:mean-equation}
For every rectangle \(D\Subset(-r_0/2,r_0/2)\times I_1\), all sufficiently large \(j\), \(k=1,2\), and \(\varphi\in C_c^\infty(D)\),
\begin{equation}\label{eq:full-weak-equation}
\left|\int_D\widehat U_{j,k}\Delta\varphi\right| \le Ca_j^{6/5}\|\varphi\|_{C^3(D)}.
\end{equation}
\end{proposition}
\begin{proof}
Fix \(j,k\), suppress \(j\), and choose \(I\Subset I_1\) with \(D\subset\mathbb R\times I\). Extend \(\varphi\in C_c^\infty(D)\) by zero and put \(\varphi_\sigma(r,z)=\varphi(\sigma r,z)\). Integration by parts gives
\[
\begin{aligned}
&\int_D\widehat U_k\Delta\varphi \\
={}& -\sum_{\sigma=\pm}\int_{(0,r_0)\times I} \nabla \widetilde u_{k,\sigma}\cdot \nabla \varphi_\sigma -\int_I(\widetilde u_{k,+}(0,z)-\widetilde u_{k,-}(0,z)) \varphi_q(0,z) dz.
\end{aligned}
\]
The functions \(\widetilde u_{k,\pm}\) satisfy
\begin{equation}\label{eq:extension-weak-gradient}
\left|\sum_{\sigma=\pm}\int_{(0,r_0)\times I} \nabla \widetilde u_{k,\sigma}\cdot \nabla \varphi_\sigma\right| \le Ca^{3/2}\|\varphi\|_{C^3(\mathbb R^2)}.
\end{equation}
Equations~\eqref{eq:general-Hone} and \eqref{eq:general-rho} give
\begin{equation}\label{eq:extension-gradient-lone}
\sum_{\sigma=\pm}\int_I\int_0^{2\rho(z)} |\nabla \widetilde u_{k,\sigma}|\mathrm d r\mathrm d z \le Ca\left(\int_I\rho\right)^{1/2}\le Ca^{3/2}.
\end{equation}
If \(\varphi(0,z)=0\), put \(u=u_{k,\sigma}\) and let \(\chi\) be smooth, with \(\chi=0\) on \((-\infty,1]\) and \(\chi=1\) on \([2,\infty)\). Multiply \(\operatorname{div}G(\nabla u)=0\) by \(\chi(r/\rho(z))\varphi_\sigma\) and integrate by parts over \(\Omega_\rho\). Since \(\operatorname{Lip}\rho\le\lambda\), \(|\nabla(\chi(r/\rho(z))\varphi_\sigma)|\le C\|\varphi\|_{C^1(\mathbb R^2)}\). By \eqref{eq:general-Hone}, \(\|G(\nabla u)-\nabla u\|_{L^1(\Omega_\rho)}\le Ca^2\). Together with \eqref{eq:extension-gradient-lone}, this proves \eqref{eq:extension-weak-gradient}.

Next let \(\zeta\in C_c^\infty((-r_0^2/4,r_0^2/4))\) with \(\zeta=1\) near \(0\) and take \(\varphi(q,z)=\zeta(q^2)g(z)\), \(g\in C_c^\infty(I)\). Choose \(b\in e_3^\perp\) with \(b\cdot \mathbf{n}_\ell=\delta_{k\ell}\) and set \(Y(X)=b\zeta(|X_h|^2)g(z)\). On \(T,F\) from Corollary~\ref{cor:transverse-remainder}, we have \(Y=bg(z)\) for large \(j\), since \(\sup\rho_j\to0\). Hence \(\operatorname{div}_M Y=-(b\cdot N)N_zg'\) on \(T\) and \(\operatorname{div}_C Y=0\) on \(F\). By \eqref{eq:finite-tilt} and \eqref{eq:general-core},
\[
\int_T|\operatorname{div}_M Y| \mathrm{d}\mu \le C\|g'\|_{L^\infty(I)}\mu(T)^{1/2} \left(\int_TN_z^2 \mathrm{d}\mu\right)^{1/2} \le Ca^{3/2}\|\varphi\|_{C^1(\mathbb R^2)}.
\]
For \(u=u_{\ell,\sigma}\) and \(x=rt_{\ell,\sigma}+ze_3\),
\[
\begin{aligned}
&(\operatorname{div}_M Y)(x+u\mathbf{n}_\ell)\sqrt{1+|\nabla u|^2} -(\operatorname{div}_{P_\ell}Y)(x) \\
={}& (b\cdot \mathbf{n}_\ell)\nabla u\cdot \nabla \varphi_\sigma +O\left((u^2+|\nabla u|^2)\|\varphi\|_{C^3(\mathbb R^2)}\right).
\end{aligned}
\]
Since \(\int_M\operatorname{div}_M Y \mathrm d\mu =\int_C\operatorname{div}_C Y \mathrm d\mu_C=0\), \eqref{eq:general-Hone} and \eqref{eq:extension-gradient-lone} imply \eqref{eq:extension-weak-gradient}. For arbitrary \(\varphi\), \eqref{eq:extension-weak-gradient} holds for both \(\varphi-\zeta(q^2)\varphi(0,z)\) and \(\zeta(q^2)\varphi(0,z)\), hence for \(\varphi\).

To estimate \(\widetilde u_{k,+}(0,\cdot)-\widetilde u_{k,-}(0,\cdot)\), put \(h_0=a^{4/5}\) and define
\[
\nu(J)=\int_{M\cap Z\cap\{z\in J\}}N_z^2 \mathrm{d}\mu, \qquad \nu(\mathbb R)\le Ca^2,
\]
and
\[
E=\{t\in I:\rho(t)>h_0/4\} \cup\left\{t\in I: \sup_{s>0}\frac{\nu((t-s,t+s))}{2s}>h_0^2\right\}.
\]
By \eqref{eq:general-tail} and the weak \((1,1)\) estimate for the maximal function,
\begin{equation}\label{eq:bad-height-size}
|E|\le Ca^2/h_0^2=Ca^{2/5}.
\end{equation}
For \(t\in I\setminus E\), \(h\in[h_0,2h_0]\), and \(|z-t|<2h\), \(\rho(z)\le h_0/4+2\lambda h<h/2\). For small \(h_0\), we also have \(h\ge A_\alpha h_0^2\) and \(B_{64h}(te_3)\subset Z\). Thus Lemma~\ref{lem:tilt-iteration} applies with \(\kappa=h_0^2\). By \eqref{eq:tangent-concurrency},
\[
|u_{k,+}(h,t)-u_{k,-}(h,t)| \le Ch_0^{3/2}+Ch\sum_{\sigma=\pm}|\partial_r u_{k,\sigma}(h,t)|.
\]
Integrating \(\partial_r\widetilde u_{k,\sigma}\) from \(0\) to \(h\) and averaging over \(h\in[h_0,2h_0]\), we obtain
\[
|\widetilde u_{k,+}(0,t)-\widetilde u_{k,-}(0,t)| \le Ch_0^{3/2} +C\sum_{\sigma=\pm}\int_0^{2h_0}|\partial_r\widetilde u_{k,\sigma}(r,t)| dr.
\]
Integrating over \(I\setminus E\) and using \eqref{eq:general-Hone} and \eqref{eq:bad-height-size} on \(E\), we obtain
\begin{equation}\label{eq:power-trace}
\|\widetilde u_{k,+}(0,\cdot)-\widetilde u_{k,-}(0,\cdot)\|_{L^1(I)} \le C\left(\frac{a^2}{h_0}+h_0^{3/2}+a\sqrt{h_0}\right) \le Ca^{6/5}. \qedhere
\end{equation}
\end{proof}

\begin{proposition}\label{cor:linear-matching}\label{lem:linear}
Let \(\Phi_{\sigma_j}:A_0\to A_{\sigma_j}\) be the maps from Section~\ref{sec:transverse-graphs}, with \(\Phi_{\sigma_j}\to\mathrm{id}\) in \(C^2\). If \(\delta_j/H_j\to0\), then, after passing to a subsequence, there are constants \(A_k,B_k\), \(k=1,2\), such that
\begin{equation}\label{eq:linear}
\frac{u_{j,k,\sigma}\circ\Phi_{\sigma_j}}{a_j} \longrightarrow\ell_{k,\sigma}(r,z)=\sigma A_kr+B_kz \quad \text{ in }L^2(A_0),\quad \sum_i\int_{A_0}\ell_i^2=1,
\end{equation}
and
\begin{equation}\label{eq:extended-graph-strong-convergence}
\widetilde u_{j,i}/a_j\longrightarrow\ell_i \quad \text{ in }L^2((0,r_0)\times I_2).
\end{equation}
\end{proposition}
\begin{proof}
By \eqref{eq:initial-graph-h1-bound}, \eqref{eq:general-Hone}, and Rellich's theorem, after passing to a subsequence we have
\[
\begin{aligned}
u_{j,i}\circ\Phi_{\sigma_j}/a_j&\longrightarrow\ell_i &&\text{in }L^2(A_0),\\
\widetilde u_{j,i}/a_j&\longrightarrow\ell_i &&\text{in }L^2((0,r_0)\times I_2).
\end{aligned}
\]

On \(\Omega=\{(q,z):c<\sqrt{q^2+z^2}<d\}\), define, for \(\sigma q>0\),
\[
U_k(q,z)=\ell_{k,\sigma}(\sigma q,z).
\]
Since \(\operatorname{div}G(\nabla u_{j,i})=0\) and \(\|\nabla u_{j,i}\|_{L^\infty(A_{\sigma_j})}\to0\), \(\Delta U_k=0\) in \(\Omega\setminus\{q=0\}\). Dividing \eqref{eq:full-weak-equation} by \(a_j\) and passing to the limit gives \(\Delta U_k=0\) on \(\Omega\). By \eqref{eq:eulerexact},
\[
\frac{\|\mathcal E u_{j,i}-u_{j,i}\|_{L^2(A_{\sigma_j})}}{a_j}
\le C\left(\frac{\delta_j}{H_j+\delta_j}\right)^{1/2}\longrightarrow0,
\]
hence \(q\partial_q U_k+z\partial_z U_k=U_k\). Thus \(U_k=A_kq+B_kz\), and
\[
\sum_i\int_{A_0}\ell_i^2 =\lim_{j\to\infty}\frac{H_j}{a_j^2}=1.\qedhere
\]
\end{proof}

\subsection{Uniqueness of the tangent cone}\label{sec:transverse-decay}

Recall \(e(R)=2-\mu(B_R)/(\pi R^2)\), and put \(d(R)=e(R)-e(2R)\).

\begin{theorem}\label{thm:transverse-tangent-cone-uniqueness}
\(|R^{-1}M|\to C\) as varifolds as \(R\to\infty\). Moreover, \(e(R)\le Kd(R)\) for all sufficiently large \(R\) and some \(K>2\).
\end{theorem}

For related uniqueness and regularity results at singular points, see Simon~\cite{SimonCylindrical} and Becker-Kahn~\cite{BeckerKahn}.

Fix radii \(0<a_0<c<a_1<a_2<d<a_3\), and write \(A(a,b)=\{a<|X|<b\}\). Let \(\mathcal K\) be the compact space of cones \(C=|P_1|+|P_2|\), with the planes unordered and possibly coincident. For \(C'=|P'_1|+|P'_2|\), use the metric
\[
d_{\mathcal K}(C,C')= \min_{\sigma\in S_2}\max_{k=1,2}\angle(P_k,P'_{\sigma(k)}).
\]
Write \(e_T,d_T\) for \(e,d\) with \(M\) replaced by \(T\). For a sequence as in Section~\ref{sec:transverse-graphs}, set \(e_j=e_{M_j}\) and define
\begin{equation}\label{eq:best-two-plane-approximation}
\mathcal F_j(Q)=\int_{M_j\cap A(a_1,a_2)}\operatorname{dist}(X,\operatorname{spt} Q)^2 \mathrm{d}\mu_j, \qquad C_j\in\operatorname*{argmin}_{Q\in\mathcal K}\mathcal F_j(Q).
\end{equation}
Since \(\mu_C(\partial A(a_1,a_2))=0\), \(\mathcal F_j(C_j)\le\mathcal F_j(C)\to0\). Every subsequential limit of \(C_j\) therefore contains \(\operatorname{spt}C\cap A(a_1,a_2)\) and equals \(C\). Thus \(C_j\to C\).

After applying rotations \(Q_j\to I\) to \(M_j\) and \(C_j\), assume that the planes in \(C_j\) intersect along \(\mathbb R e_3\). Write \(C_j=|P_{j,1}|+|P_{j,2}|\), and choose unit vectors \(t_{j,k,+}\in P_{j,k}\cap e_3^\perp\) and \(\mathbf n_{j,k}\perp P_{j,k}\). Let \(u_{j,i}\) be the four graph functions over \(C_j\) on \(A_{\sigma_j}(c,d)\) satisfying \eqref{eq:initial-graph-smallness}, and set
\[
H_j=\sum_i\int_{A_{\sigma_j}(c,d)}u_{j,i}^2.
\]

\begin{proposition}\label{prop:graph-error-density-drop}
For all sufficiently large \(j\),
\[
H_j+\mathcal F_j(C_j)\le C\left(e_j(a_0)-e_j(a_3)\right).
\]
The constant depends only on \(a_i,c,d\) and a positive lower bound for \(\alpha\).
\end{proposition}
\begin{proof}
Put \(\Delta_j=e_j(a_0)-e_j(a_3)\) and \(a_j=(H_j+\Delta_j)^{1/2}\). Apply Proposition~\ref{prop:transverse-graphs} near \(\overline{A(a_1,a_2)}\cap\mathbb R e_3\), using \((a,b)=(a_0,a_3)\) in Lemma~\ref{lem:intermediate}. Let \(T_j\) be the part of \(M_j\cap A(a_1,a_2)\) not covered by these graphs. Since \(\operatorname{dist}(X_{j,i},\operatorname{spt}C_j)\le|u_{j,i}|\), \eqref{eq:initial-graph-h1-bound}, \eqref{eq:general-Hone}, and Corollary~\ref{cor:transverse-remainder} imply \(\mathcal F_j(C_j)\le Ca_j^2\).

It therefore suffices to prove \(H_j\le C\Delta_j\). We argue by contradiction and assume that \(\Delta_j/H_j\to0\). By Proposition~\ref{cor:linear-matching}, after passing to a subsequence we have \(u_{j,i}\circ\Phi_{\sigma_j}/a_j\to\ell_i\) in \(L^2(A_0(c,d))\), where \(\ell_{k,\sigma}=\sigma A_kr+B_kz\). Define
\begin{equation}\label{eq:competitor}
P'_{j,k} =\{q t_{j,k,+}+ze_3+a_j(A_kq+B_kz)\mathbf{n}_{j,k}:q,z\in\mathbb{R}\},\quad C'_j =|P'_{j,1}|+|P'_{j,2}|.
\end{equation}
For \(i=(k,\sigma)\),
\[
\operatorname{dist}(X_{j,i}(r,z),\operatorname{spt}C'_j) \le |u_{j,i}(r,z)-a_j\ell_i(r,z)|.
\]
On \(T_j\), \(\operatorname{dist}(X,\operatorname{spt}C'_j)\le |X_h|+Ca_j\). By \eqref{eq:linear}, \eqref{eq:extended-graph-strong-convergence}, and \eqref{eq:general-core},
\[
\mathcal F_j(C_j)\le\mathcal F_j(C'_j) \le C\int_{T_j}|X_h|^2 \mathrm d\mu_j+o(a_j^2)=o(a_j^2).
\]
On each compact subset of \(A_0(a_1,a_2)\), \(\operatorname{dist}(X_{j,i},\operatorname{spt}C_j)=|u_{j,i}|\) for large \(j\). Hence \(\ell_i\equiv0\) for every \(i\), contradicting \(\sum_i\int_{A_0(c,d)}\ell_i^2=1\). Therefore \(H_j\le C\Delta_j\).
\end{proof}

Set \(\delta_j=e_j(1)-e_j(2)\).
\begin{proposition}\label{prop:density-deficit-drop}
There is \(K>2\) such that, for all sufficiently large \(j\),
\begin{equation}\label{eq:density-deficit-drop-bound}
e_j(1)\le K\delta_j.
\end{equation}
\end{proposition}
\begin{proof}
By Proposition~\ref{prop:graph-error-density-drop} with \(a_0=1\) and \(a_3=2\), we have \(H_j\le C\delta_j\). Let \(u_{j,i},T_j,F_j\) be as in Proposition~\ref{prop:transverse-graphs} and Corollary~\ref{cor:transverse-remainder}, applied to \(M_j,C_j\) near \(\overline{A(a_1,a_2)}\cap\mathbb R e_3\). Then
\begin{equation}\label{eq:graphsmall}
\sum_i\int_{\Omega_{\rho_j}}(u_{j,i}^2+|\nabla u_{j,i}|^2) +\int_{T_j}(|X_h|^2+N_z^2) \mathrm d\mu_j +\int_{F_j}|X_h|^2 \mathrm d\mu_{C_j} \le C\delta_j.
\end{equation}

Fix \(0\le\phi\in C_c^\infty((a_1,a_2))\) with \(2\pi\int s\phi(s) \mathrm ds=1\), and put \(f(X)=\phi'(|X|)|X_h|^2/|X|\). For \(Y(X)=\phi(|X|)X_h\), the stationarity of \(M_j\) and \(C_j\) gives
\[
\begin{aligned}
&\int\phi(|X|) \mathrm d(\mu_j-\mu_{C_j}) \\
={}& -\int f \mathrm d(\mu_j-\mu_{C_j}) +\int\left(\frac{\phi'(|X|)}{|X|}(X_h\cdot N)(X\cdot N)-\phi(|X|)N_z^2\right) \mathrm d\mu_j.
\end{aligned}
\]
For \(i=(k,\sigma)\), put \(x=\sigma r t_{j,k,+}+ze_3\). Since \(\nabla f(x)\cdot\mathbf n_{j,k}=0\),
\[
\left|f(X_{j,i})\sqrt{1+|\nabla u_{j,i}|^2}-f(x)\right| \le C(u_{j,i}^2+|\nabla u_{j,i}|^2).
\]
Since \(|f(X)|\le C|X_h|^2\) and \(X_h\cdot N=X\cdot N-zN_z\), \eqref{eq:graph-error-density-drop}, \eqref{eq:graphsmall}, \eqref{eq:initial-graph-h1-bound} give
\begin{equation}\label{eq:radialmass}
\left|\int\phi(|X|) \mathrm d(\mu_j-\mu_{C_j})\right| \le C_\phi\delta_j.
\end{equation}

On \(\operatorname{spt}\phi\), \(|e_j(s)-e_j(1)|\le\delta_j\). Since \(\mu_j(B_s)-\mu_{C_j}(B_s)=-\pi s^2e_j(s)\), integration by parts and \eqref{eq:radialmass} give
\[
\begin{aligned}[b]
e_j(1) &=-\pi\int\phi'(s)s^2e_j(s) \mathrm ds+O_\phi(\delta_j) \\
&=-\int\phi(|X|) \mathrm d(\mu_j-\mu_{C_j}) +O_\phi(\delta_j)\le C_\phi\delta_j.
\end{aligned}\qedhere
\]
\end{proof}

Fix \(0\le\chi\in C_c^\infty((1,2))\) with \(\int_1^2s\chi(s) \mathrm{d}s=1/2\). For a surface \(T\), define a measure on \(S^2\) by
\[
\Lambda_T(R)[f]=R^{-2}\int_T\chi(|X|/R)f(X/|X|) \mathrm{d}\mathcal H^2, \qquad f\in C(S^2),
\]
and write \(\Lambda(R)=\Lambda_M(R)\). For \(C=|P_1|+|P_2|\in\mathcal K\),
\[
\sum_{k=1}^2\Lambda_{P_k}(R)=\tfrac12\sum_{k=1}^2\mathcal H^1\llcorner(P_k\cap S^2)=:\Lambda_C.
\]
For a finite signed measure \(\nu\) on \(S^2\), put
\[
\|\nu\|_*= \sup_{\substack{f\in C^1(S^2)\\ \|f\|_{L^\infty(S^2)}+\|\nabla f\|_{L^\infty(S^2)}\le1}}|\nu[f]|.
\]
We have \(\|\Lambda_{T_j}(1)-\Lambda_C\|_*\to0\) whenever \(|T_j|\to C\) as varifolds. The map \(C\mapsto\Lambda_C\) is continuous and injective on \(\mathcal K\). Also, \(\Lambda_{R^{-1}T}(s)=\Lambda_T(Rs)\).

For \(f\in C^1(S^2)\), let \(\widetilde f(X)=f(X/|X|)\).
As \(\int_M\operatorname{div}_M(\chi(|X|/R)\widetilde f X) \mathrm d\mu=0\),
\[
\begin{aligned}
&R\frac{\mathrm{d}}{\mathrm{d}R}\Lambda(R)[f] \\
={}& -R^{-2}\int_M\chi(|X|/R)(X\cdot N)(\nabla\widetilde f\cdot N) \mathrm{d}\mu -R^{-3}\int_M\chi'(|X|/R)\widetilde f\frac{(X\cdot N)^2}{|X|} \mathrm{d}\mu.
\end{aligned}
\]
By \eqref{eq:monotonicity},
\[
\left|R\frac{\mathrm{d}}{\mathrm{d}R}\Lambda(R)[f]\right| \le C\|f\|_{C^1(S^2)}\left(\sqrt{d(R)}+d(R)\right) \le C\|f\|_{C^1(S^2)}\sqrt{e(R)}.
\]
Monotonicity of \(e\) and integration over \([R,s]\subset[R,2R]\) give
\begin{equation}\label{eq:Qvariation}
\sup_{R\le s\le2R}\|\Lambda(s)-\Lambda(R)\|_* \le C\sqrt{e(R)}.
\end{equation}

\begin{proof}[Proof of Theorem~\ref{thm:transverse-tangent-cone-uniqueness}]
Choose \(\eta>0\) with \(\|\Lambda_C-\Lambda_{2|P|}\|_*>2\eta\) for every plane \(P\) through the origin. There are \(K>2\) and \(R_*>0\) such that
\begin{equation}\label{eq:neighborhood}
e(R)\le Kd(R)
\end{equation}
whenever \(R\ge R_*\) and \(\|\Lambda(R)-\Lambda_C\|_*\le\eta\).
Otherwise there are \(R_j\to\infty\) with \(\|\Lambda(R_j)-\Lambda_C\|_*\le\eta\) and \(e(R_j)>j d(R_j)\). Since \(e(sR_j)\to0\) for every \(s>0\), after passing to a subsequence we have \(|R_j^{-1}M|\to C'\), where \(C'\) is a cone of density two with \(\|\Lambda_{C'}-\Lambda_C\|_*\le\eta\). By Lemma~\ref{lem:cones}, \(C'=|P'_1|+|P'_2|\) with \(P'_1\ne P'_2\), contradicting Proposition~\ref{prop:density-deficit-drop}.

Put \(q=1-1/K\). Choose \(R_j\to\infty\) with \(|R_j^{-1}M|\to C\), and take \(R_0=R_j\ge R_*\) for \(j\) large enough that
\begin{equation}\label{eq:initial-radius}
\|\Lambda(R_0)-\Lambda_C\|_* +\frac{C\sqrt{e(R_0)}}{1-\sqrt q}<\eta.
\end{equation}
Equations~\eqref{eq:neighborhood} and \eqref{eq:Qvariation} give inductively, for every \(n\ge0\),
\[
e(2^nR_0)\le q^n e(R_0),\quad \sup_{R_0\le R\le2^{n+1}R_0}\|\Lambda(R)-\Lambda_C\|_*<\eta.
\]
Thus \eqref{eq:neighborhood} applies at every \(R\ge R_0\). With \(\gamma_0=-\log_2 q\in(0,1)\), the monotonicity of \(e\) gives
\begin{equation}\label{eq:uf-capture}
e(R)\le Kd(R), \qquad e(R)\le AR^{-\gamma_0}.
\end{equation}
Since \(e(2R)\le qe(R)\) for \(R\ge R_0\), \eqref{eq:Qvariation} and \(\Lambda(R_j)\to\Lambda_C\) imply
\[
\|\Lambda(R)-\Lambda_C\|_*\le C\sum_{k=0}^{\infty}\sqrt{e(2^kR)}\le\frac{C\sqrt{e(R)}}{1-\sqrt q}\longrightarrow0.
\]
Every tangent cone \(C'\) of \(M\) at infinity satisfies \(\Lambda_{C'}=\Lambda_C\), hence \(C'=C\).
\end{proof}

Fix a transverse cone \(C=|P_1|+|P_2|\) with intersection axis \(\mathbb R e_3\), and let \(\mathcal C(C)\) be the class of connected, properly embedded minimal surfaces without boundary, of density two at infinity, whose unique tangent cone at infinity is \(C\).

\begin{corollary}\label{cor:cone-limit}
Let \(T_j\in\mathcal C(C)\) and \(|T_j|\to V\) as varifolds. If \(V\) has density two at infinity, then \((\eta_{0,R})_\#V\to C\) as varifolds as \(R\to\infty\).
\end{corollary}
\begin{proof}
Suppose that \(V\) has a tangent cone \(D\ne C\) at infinity. By Lemma~\ref{lem:cones}, \(D\in\mathcal K\). After passing to a subsequence, there are \(s_j\to\infty\) with \(S_j=s_j^{-1}T_j\in\mathcal C(C)\) and \(|S_j|\to D\). Thus \(e_{S_j}(1)\to0\).

Choose \(\delta>0\) smaller than \(\|\Lambda_D-\Lambda_C\|_*\) and \(\|\Lambda_{2|P|}-\Lambda_C\|_*\) for every plane \(P\) through \(0\). For each fixed \(r\ge1\), \(\|\Lambda_{S_j}(r)-\Lambda_C\|_*\to\|\Lambda_D-\Lambda_C\|_*>\delta\). Since \(\Lambda_{S_j}(r)\to\Lambda_C\) as \(r\to\infty\) for each \(j\), the radius
\[
R_j=\max\{r\ge1: \|\Lambda_{S_j}(r)-\Lambda_C\|_*=\delta\}
\]
is well defined for large \(j\), tends to infinity, and satisfies \(\|\Lambda_{S_j}(r)-\Lambda_C\|_*\le\delta\) for \(r\ge R_j\).

For every fixed \(s>0\) and large \(j\), \(0\le e_{S_j}(sR_j)\le e_{S_j}(1)\to0\). Every varifold limit \(C'\) of \(|r_j^{-1}S_j|\) with \(r_j\ge R_j\) is a cone of density two satisfying \(\|\Lambda_{C'}-\Lambda_C\|_*\le\delta\). Thus \(C'=|P'_1|+|P'_2|\), where \(\angle(P'_1,P'_2)\ge c(C,\delta)>0\). Proposition~\ref{prop:density-deficit-drop} gives \(K>2\), independent of \(j\), such that \(e_{S_j}(r)\le Kd_{S_j}(r)\) for all \(r\ge R_j\) and large \(j\). Summing \eqref{eq:Qvariation} on dyadic scales gives
\[
\delta=\|\Lambda_{S_j}(R_j)-\Lambda_C\|_* \le C\sqrt{e_{S_j}(R_j)}\longrightarrow0.\qedhere
\]
\end{proof}

\subsection{Uniform decay of the density deficit}\label{sec:uniform-transverse-decay}

Let \(S_j,S\in\mathcal C(C)\) satisfy \(|S_j|\to|S|\) as varifolds, and put \(\mathscr F=\{S,S_1,S_2,\ldots\}\).

\begin{proposition}\label{prop:uniform-decay}
For every \(0<\gamma<1\), there are constants \(A_\gamma,R_\gamma<\infty\) such that, for \(R\ge R_\gamma\) and \(T\in\mathscr F\),
\begin{equation}\label{eq:uf-decay}
e_T(R)\le A_\gamma R^{-\gamma}.
\end{equation}
\end{proposition}

We first obtain \(e_T(R)\le AR^{-\gamma_0}\) for some \(\gamma_0\in(0,1)\), uniformly in \(T\in\mathscr F\).

\begin{lemma}\label{lem:uniform-transverse-estimates}
There are \(A,R_0>0\), \(K>2\), and \(\gamma_0\in(0,1)\) such that, for \(T\in\mathscr F\) and \(R\ge R_0\),
\begin{equation}\label{eq:baseline-drop}
e_T(2R)\le2^{-\gamma_0}e_T(R), \qquad e_T(R)\le Kd_T(R),\quad e_T(R)\le AR^{-\gamma_0}.
\end{equation}
Moreover, \(|R^{-1}T|\to C\) uniformly in \(T\in\mathscr F\). For some \(R_*\ge8R_0\) and \(\sigma_R\downarrow0\), every \(R^{-1}T\), \(R\ge R_*\), has four graphs \(u_{T,R,i}\) over \(C\) on \(A_{\sigma_R}(3/4,3)\). Put
\[
Q_{T,R}=\sum_i\int_{A_{\sigma_R}(3/4,3)}u_{T,R,i}^2 +e_T(R/2)-e_T(4R).
\]
With \(\widetilde u_{T,R,i}\) as in Proposition~\ref{prop:transverse-graphs},
\begin{gather}
Q_{T,R}+\sum_i\|\widetilde u_{T,R,i}\|_{H^1((0,r_0)\times I_2)}^2 \le Ce_T(R/8),\label{eq:uf-budget} \\
\int_{T\cap A(R,2R)}N_z^2 \mathrm{d}\mu_T \le CR^2e_T(R/8).\label{eq:physical-budget}
\end{gather}
\end{lemma}
\begin{proof}
Since \(\sup_{T\in\mathscr F}e_T(R)\to0\) by local convergence and the monotonicity of \(e_T\), the constants \(K,R_*\) in \eqref{eq:neighborhood} can be chosen independently of \(T\in\mathscr F\). We may take \(R_0\ge R_*\) so that \eqref{eq:initial-radius} holds for \(S\); by local convergence, it also holds for \(S_j\) with large \(j\). The proof of Theorem~\ref{thm:transverse-tangent-cone-uniqueness} then gives \eqref{eq:baseline-drop} for \(S\) and these \(S_j\). Since only finitely many \(S_j\) remain, \eqref{eq:baseline-drop} holds with the same constants \(A,R_0\) for all \(T\in\mathscr F\).

By \eqref{eq:Qvariation} and \eqref{eq:baseline-drop},
\[
\sup_{T\in\mathscr F}\|\Lambda_T(R)-\Lambda_C\|_*\le CR^{-\gamma_0/2}.
\]
Every varifold limit of \(|R_j^{-1}T_j|\), with \(T_j\in\mathscr F\) and \(R_j\to\infty\), is a cone of density two by \eqref{eq:baseline-drop} and equals \(C\) since \(\Lambda_{T_j}(R_j)\to\Lambda_C\).

Write \(e=e_T\) and \(\mu=\mu_T\). For each sufficiently large \(R\), choose
\[
C_R\in\operatorname*{argmin}_{Q\in\mathcal K} \int_{T\cap A(R/2,4R)}\operatorname{dist}(X,\operatorname{spt} Q)^2 \mathrm{d}\mu.
\]
As in Section~\ref{sec:transverse-decay}, \(C_R\to C\) uniformly in \(T\). By Proposition~\ref{prop:graph-error-density-drop},
\begin{equation}\label{eq:fitted-annulus}
R^{-4}\int_{T\cap A(R/2,4R)}\operatorname{dist}(X,\operatorname{spt} C_R)^2 \mathrm{d}\mu \le Ce(R/8).
\end{equation}
On fixed rectangles compactly contained in \(P_k\cap A(1,2)\setminus\mathbb R e_3\), \(R^{-1}T\) and the planes in \(C_R\) and \(C_{2R}\) that converge to \(P_k\) are graphs over \(P_k\), so
\[
\begin{aligned}
d_{\mathcal K}(C_R,C_{2R})^2 &\le CR^{-4}\int_{T\cap A(R,2R)} \left(\operatorname{dist}(X,\operatorname{spt} C_R)^2+\operatorname{dist}(X,\operatorname{spt} C_{2R})^2\right) \mathrm{d}\mu \\
&\le Ce(R/8).
\end{aligned}
\]
By \eqref{eq:baseline-drop},
\begin{equation}\label{eq:fixedcone-rate}
d_{\mathcal K}(C_R,C)\le C\sum_{m\ge0}\sqrt{e(2^mR/8)} \le C\sqrt{e(R/8)}.
\end{equation}
By \eqref{eq:area-bound} and \eqref{eq:fitted-annulus}--\eqref{eq:fixedcone-rate},
\begin{equation}\label{eq:fixedcone-height}
R^{-4}\int_{T\cap A(R/2,4R)}\operatorname{dist}(X,\operatorname{spt} C)^2 \mathrm{d}\mu \le Ce(R/8).
\end{equation}

By smooth convergence away from the axis, the four graphs can be chosen with \(\sigma_R\downarrow0\) independent of \(T\in\mathscr F\). On the graph of \(u_{T,R,(k,\sigma)}\), we have \(1/2<|X|<4\) and \(\operatorname{dist}(X,\operatorname{spt}C)=\operatorname{dist}(X,P_k)\). Thus \(Q_{T,R}\le Ce_T(R/8)\) by \eqref{eq:fixedcone-height} and the monotonicity of \(e_T\). Apply Proposition~\ref{prop:transverse-graphs} and Corollary~\ref{cor:transverse-remainder} to \(R^{-1}T\) near \(\{1\le|X|\le2\}\cap\mathbb R e_3\). Equation~\eqref{eq:general-Hone} proves \eqref{eq:uf-budget}; \eqref{eq:initial-graph-h1-bound} and \eqref{eq:finite-tilt} prove \eqref{eq:physical-budget}.
\end{proof}

Put \(\tau_k=t_{k,+}\) and \(\eta_k=\mathbf n_k\), \(k=1,2\). For \(T\in\mathscr F\), extend \(u_{T,R,i}\) as in Proposition~\ref{prop:transverse-graphs}. For each \(i\), the functions \(R\,u_{T,R,i}(\cdot/R)\), \(R\ge R_*\), agree on overlaps. Define \(u_{T,k}\) on \(P_k\) by
\[
u_{T,k}(\pm Rr,Rz)=R\,u_{T,R,(k,\pm)}(r,z)
\]
where the right-hand side is defined for some \(R\ge R_*\), and by zero elsewhere on \(P_k\).

Put \(A_R=\{R<|x|<2R\}\).

\begin{lemma}\label{lem:uniform-harmonic-approximation}
There are a universal \(\theta>1/2\) and a finite cover of \(\overline{A_1}\) by disks \(D^{(\ell)}\Subset\mathbb R^2\setminus\{0\}\) such that, for \(T\in\mathscr F\), \(R\ge R_*\), and \(k=1,2\), we have
\[
R^{-1}\|u_{T,k}\|_{L^2(A_R)}\le CR e_T(R/8)^{1/2},
\]
and on each \(RD^{(\ell)}\) there is a harmonic function \(h\) satisfying
\begin{equation}\label{eq:height-harmonic}
R^{-1}\|u_{T,k}-h\|_{L^2(RD^{(\ell)})}\le CR e_T(R/8)^\theta.
\end{equation}
\end{lemma}
\begin{proof}
Write \(w=R^{-1}u_{T,k}(R \cdot)\) and \(Q=Q_{T,R}\). Let \(\widehat U\) be as in Proposition~\ref{prop:mean-equation} applied to \(R^{-1}T\), and choose rectangles \(\Omega\Subset D\Subset(-r_0/2,r_0/2)\times I_1\). Since \(w=\widehat U\) for \(|r|>\rho(z)\) and \(|w|+|\widehat U|\le C\rho(z)\) otherwise, \eqref{eq:general-rho} gives
\[
\|w-\widehat U\|_{L^2(D)}^2\le C\int\rho^3\le CQ^{4/3}.
\]
For \(|y|<\min\{1,\operatorname{dist}(\Omega,\partial D)\}\), \eqref{eq:general-Hone} gives
\[
\begin{aligned}
\|\widehat U(\cdot-y)-\widehat U\|_{L^2(\Omega)}^2\le{}& C|y|^2\int_{D\setminus\{r=0\}}|\nabla\widehat U|^2\mathrm dx \\
&+C|y|\sup_{|r|<r_0/2}\int_{I_1}|\widehat U(r,z)|^2\mathrm dz\le CQ|y|.
\end{aligned}
\]
Put \(U_\varepsilon=\varphi_\varepsilon*\widehat U\), with \(\varphi_\varepsilon\) an even mollifier and \(0<\varepsilon<\min\{1,\operatorname{dist}(\Omega,\partial D)\}\). Then
\[
\|U_\varepsilon-\widehat U\|_{L^2(\Omega)}\le C\sqrt Q \varepsilon^{1/2}, \qquad \|\Delta U_\varepsilon\|_{L^\infty(\Omega)}\le CQ^{3/5}\varepsilon^{-5},
\]
where the second estimate follows from \eqref{eq:full-weak-equation}.
Let \(h\) be harmonic on \(\Omega\) with boundary values \(U_\varepsilon\). Then
\[
\|w-h\|_{L^2(\Omega)} \le C\left(Q^{2/3}+Q^{1/2}\varepsilon^{1/2} +Q^{3/5}\varepsilon^{-5}\right).
\]
Choose a universal \(\theta>1/2\) sufficiently close to \(1/2\), and set \(\varepsilon=Q^{2\theta-1}\). Then \(\|w-h\|_{L^2(\Omega)}\le CQ^\theta\).

On a disk \(D\) away from the axis, the harmonic function \(h\) with boundary values \(w\) satisfies
\[
\|w-h\|_{L^2(D)} \le C\left\|(1-(1+|\nabla w|^2)^{-1/2})\nabla w\right\|_{L^2(D)} \le CQ^{3/2}.
\]
Choose disks \(D^{(\ell)}\) covering \(\overline{A_1}\), each contained in a rectangle \(\Omega\) near the axis or a disk \(D\) away from it. By \eqref{eq:uf-budget}, \(Q\le Ce_T(R/8)\), so \(Rh(\cdot/R)\) satisfies \eqref{eq:height-harmonic} on \(RD^{(\ell)}\).
\end{proof}

\begin{proof}[Proof of Proposition~\ref{prop:uniform-decay}]
Fix \(0<\gamma<1\). We argue by contradiction and choose \(T_j\in\mathscr F\) and \(R_j\to\infty\) such that for \(R_*\le r\le R_j\),
\begin{equation}\label{eq:record-scales}
R_j^\gamma e_{T_j}(R_j)\to\infty, \qquad r^\gamma e_{T_j}(r)\le R_j^\gamma e_{T_j}(R_j).
\end{equation}
Set \(M_j=R_j^{-1}T_j\) and \(a_j^2=e_{M_j}(1)>0\). Then \(|M_j|\to C\) as varifolds and \(a_j\to0\). Equations~\eqref{eq:record-scales} and \eqref{eq:baseline-drop} give
\begin{equation}\label{eq:record-ratios}
\frac{e_{M_j}(s)}{a_j^2}\le \begin{cases} s^{-\gamma},&R_*/R_j\le s\le1,\\ Cs^{-\gamma_0},&s>1. \end{cases}
\end{equation}
Set
\[
v_{j,k}(x)=\frac{u_{T_j,k}(R_jx)}{R_ja_j}, \qquad k=1,2.
\]
For each fixed \(s>0\), Lemma~\ref{lem:uniform-harmonic-approximation} and \eqref{eq:record-ratios} bound \(v_{j,k}\) in \(L^2(A_s)\); on each \(sD^{(\ell)}\), its \(L^2\) distance to harmonic functions is \(O_s(a_j^{2\theta-1})\to0\). After passing to a subsequence, \(v_{j,k}\to U_k\) strongly in \(L^2_{\rm loc}(P_k\setminus\{0\})\) and smoothly off the axis, with \(U_k\) harmonic and
\begin{equation}\label{eq:harmonic-growth}
\int_{A_s}|U_k|^2\le \begin{cases} Cs^{4-\gamma},&0<s\le1,\\ Cs^{4-\gamma_0},&s\ge1. \end{cases}
\end{equation}
By \eqref{eq:harmonic-growth} and interior estimates, \(U_k\) extends harmonically across \(0\), with \(U_k(0)=0\), and has sublinear growth. Hence
\begin{equation}\label{eq:harmonic-zero}
U_1=U_2=0.
\end{equation}

Put \(\lambda_j=a_j^{-2}N_z^2\mu_j\). Equations~\eqref{eq:physical-budget} and~\eqref{eq:record-ratios} imply \(\sup_j\lambda_j(K)<\infty\) for each compact \(K\subset\mathbb R^3\setminus\{0\}\) and, for \(0<s\le1\),
\begin{equation}\label{eq:normal-annulus-bound}
\limsup_{j\to\infty}\lambda_j(A(s,2s))\le Cs^{2-\gamma}.
\end{equation}
After passing to a further subsequence, \(\lambda_j\rightharpoonup\nu\) on \(\mathbb R^3\setminus\{0\}\). Since \(v_{j,k}\to0\) smoothly off the axis, \(\operatorname{spt}\nu\subset\mathbb Re_3\setminus\{0\}\).

By \eqref{eq:horizontal-height}, \eqref{eq:uf-budget}, and \eqref{eq:record-ratios}, the convergence \(v_{j,k}\to0\) in \(L^2_{\mathrm{loc}}\) rules out concentration of \(a_j^{-2}(X_h\cdot N)^2\mu_j\) on \(\mathbb Re_3\setminus\{0\}\). Together with smooth convergence off the axis, this gives for every compact \(K\subset\mathbb R^3\setminus\{0\}\)
\begin{equation}\label{eq:horizontal-normal-small}
a_j^{-2}\int_{M_j\cap K}(X_h\cdot N)^2 \mathrm{d}\mu_j\longrightarrow0.
\end{equation}

For smooth \(\psi\), the field \(Y_\psi=-\psi'(z)X_h+\psi(z)e_3\) satisfies
\begin{equation}\label{eq:axial-variation}
\operatorname{div}_M Y_\psi =-2\psi'(z)N_z^2+\psi''(z)(X_h\cdot N)N_z.
\end{equation}
For \(\psi\in C_c^\infty(\mathbb R\setminus\{0\})\), choose \(\chi\in C_c^\infty(\mathbb R)\) with \(0\le\chi\le1\) and \(\chi=1\) on \([-1,1]\). Put \(r_h=|X_h|\). For the field \(\chi(r_h)Y_\psi\), the first variation formula and \eqref{eq:horizontal-normal-small} give
\[
2\int\chi(r_h)\psi'(z) \mathrm d\lambda_j =-a_j^{-2}\int_{M_j}r_h\chi'(r_h)\psi'(z) \mathrm d\mu_j+o(1).
\]
On the support of \(\chi'(r_h)\psi'(z)\), \(v_{j,k}\to0\) in \(C^1\), so \(r_h=|r|+o(a_j^2)\) and \(\mathrm d\mu_j=(1+o(a_j^2)) \mathrm d r \mathrm d z\). Since \(\int_{\mathbb R}\psi'=0\), \(\int_{M_j}r_h\chi'(r_h)\psi'(z) \mathrm d\mu_j=o(a_j^2)\). Thus \(\int\psi'(z) \mathrm d\nu=0\), so \(\nu=c_\pm dz\) on \(\{\pm z>0\}\), with \(c_\pm\ge0\).

Equation~\eqref{eq:normal-annulus-bound} gives \((c_++c_-)s\le Cs^{2-\gamma}\); since \(\gamma<1\), \(c_+=c_-=0\). Hence \(\nu=0\) and \(\lambda_j(\{1\le|X|\le2\})\to0\). Since \(X\cdot N=X_h\cdot N+zN_z\), the monotonicity formula and \eqref{eq:baseline-drop} imply
\[
\frac1K\le\frac{e_{M_j}(1)-e_{M_j}(2)}{a_j^2} \le\frac C{a_j^2}\int_{M_j\cap\{1\le|X|\le2\}} \left((X_h\cdot N)^2+N_z^2\right) \mathrm d\mu_j \longrightarrow0.\qedhere
\]
\end{proof}

\section{Height bounds and reflection symmetry}\label{sec:wings-symmetry}

\subsection{Growth of the heights}\label{sec:wing-growth}

As in Section~\ref{sec:uniform-transverse-decay}, fix \(\mathscr F=\{S,S_1,S_2,\ldots\}\subset\mathcal C(C)\) with \(|S_j|\to|S|\) as varifolds.

\begin{proposition}\label{prop:uf-tail}
For some universal \(\beta\in(0,1/2)\), there are constants \(C,R_*>0\) such that
\begin{equation}\label{eq:uf-tail}
|u_{T,k}|+|r||\nabla u_{T,k}|+r^2|\nabla^2 u_{T,k}| \le C|r|^{\beta},
\end{equation}
for \(T\in\mathscr F\), \(k=1,2\), \(|z|\le|r|\), and \(|r|\ge R_*\).
\end{proposition}

\begin{proof}
It suffices to prove
\begin{equation}\label{eq:wing-ltwo}
R^{-1}\|u_{T,k}\|_{L^2(A_R)}\le CR^{\beta},
\end{equation}
since interior estimates on \(B_{c|r|}(r,z)\) imply \eqref{eq:uf-tail}. Let \(\theta>1/2\) and the disks \(D^{(\ell)}\) be as in Lemma~\ref{lem:uniform-harmonic-approximation}. Choose \(0<\gamma<1\) so that \(\beta:=1-\gamma\theta\in(0,1/2)\), and put \(q=1-\gamma/2\). By \eqref{eq:uf-decay}, for \(R\ge R_1\),
\[
R^{-1}\|u_{T,k}\|_{L^2(A_R)}\le CR^q,
\]
and on each \(RD^{(\ell)}\) there is a harmonic function \(h\) satisfying
\[
R^{-1}\|u_{T,k}-h\|_{L^2(RD^{(\ell)})}\le CR^\beta.
\]
If \eqref{eq:wing-ltwo} fails, choose \(T_j\in\mathscr F\) and a fixed \(k\), and put \(m_j(R)=R^{-1}\|u_{T_j,k}\|_{L^2(A_R)}\). Since \(m_j(R)\le CR^q\), we can choose \(\varepsilon_j\downarrow0\) and \(R_j\to\infty\) such that
\[
\frac{m_j(R_j)}{R_j^\beta+\varepsilon_jR_j^q} \ge\frac12\sup_{R\ge R_1}\frac{m_j(R)}{R^\beta+\varepsilon_jR^q} \longrightarrow\infty.
\]
Set \(v_j(x)=u_{T_j,k}(R_jx)/m_j(R_j)\), so that \(\|v_j\|_{L^2(A_1)}=1\). For fixed \(s>0\) and large \(j\),
\[
s^{-1}\|v_j\|_{L^2(A_s)} =\frac{m_j(sR_j)}{m_j(R_j)} \le2 \begin{cases} s^\beta,&0<s\le1,\\ s^q,&s\ge1. \end{cases}
\]
On each \(sD^{(\ell)}\), there are harmonic functions \(h_j\) with
\[
\|v_j-h_j\|_{L^2(sD^{(\ell)})} \le Cs^{\beta+1}\frac{R_j^\beta}{m_j(R_j)}\longrightarrow0.
\]
After passing to a subsequence, \(v_j\to v\) strongly in \(L^2_{\rm loc}(\mathbb R^2\setminus\{0\})\), where \(v\) is harmonic. By interior estimates, \(v(x)=O(|x|^\beta)\) as \(x\to0\) and \(v(x)=O(|x|^q)\) as \(|x|\to\infty\). Since \(\beta>0\), \(v\) extends harmonically across \(0\), with \(v(0)=0\). Since \(q<1\), \(v\equiv0\), contrary to \(\|v\|_{L^2(A_1)}=1\).
\end{proof}

\subsection{Curvature bounds and extension of the graphs}\label{sec:curvature-compactness}

\begin{proposition}\label{prop:uf-flux}
For \(T\in\mathcal C(C)\), the measures \(\nu_{T,H}:=H^{-1}(z/H)_\#(N_z^2\mu_T)\) converge locally weakly on \(\mathbb R\), as \(H\to\infty\), to \(c(T)\mathcal L^1\), where
\begin{equation}\label{eq:cm-mean}
c(T):=\lim_{H\to\infty}\frac1{2H} \int_{T\cap\{|z|\le H\}}N_z^2 \mathrm d\mu_T
\end{equation}
is finite and depends continuously on \(T\) under varifold convergence. For \(\lambda>0\) and \(b\in\mathbb R^3\),
\begin{equation}\label{eq:uf-fluxlaws}
c(\lambda(T-b))=\lambda c(T).
\end{equation}
\end{proposition}
\begin{proof}
By \eqref{eq:uf-tail}, for \(T\in\mathscr F\) and an interval \(J\subset\mathbb R\),
\begin{equation}\label{eq:wing-integral-tail}
\int_{T\cap\{z\in J, r_h>R\}}N_z^2 \mathrm d\mu_T \le C|J|R^{2\beta-1}, \qquad R\ge R_*+\sup_{z\in J}|z|.
\end{equation}

For some \(\varepsilon>0\), we have
\begin{equation}\label{eq:cm-mixed}
\int_{T\cap B_H}|X_h\cdot N| |N_z| \mathrm d\mu_T \le CH^{2-\varepsilon}.
\end{equation}
Indeed, choose \(\max\{2\beta,6/7\}<\gamma<1\). Put \(Q=Q_{T,R}\) and \(u_i=u_{T,R,i}\). Equations~\eqref{eq:uf-decay}, \eqref{eq:uf-budget} and \eqref{eq:wing-ltwo} give
\[
Q\le CR^{-\gamma}, \qquad \sum_i\int_{A_{\sigma_R}(3/4,3)}u_i^2\mathrm d A\le CR^{2\beta-2}.
\]
By \eqref{eq:finite-tilt}--\eqref{eq:horizontal-height} and Caccioppoli's inequality away from the axis,
\[
\begin{aligned}
R^{-4}\int_{T\cap A(R,2R)}(X_h\cdot N)^2\,\mathrm d\mu_T
&\le C\sum_i\int_{A_{\sigma_R}(3/4,3)}u_i^2\,\mathrm dA+CQ^{4/3}, \\
R^{-2}\int_{T\cap A(R,2R)}N_z^2\,\mathrm d\mu_T&\le CQ.
\end{aligned}
\]
By the Cauchy--Schwarz inequality,
\[
\int_{T\cap A(R,2R)}|X_h\cdot N| |N_z| \mathrm d\mu_T \le C\left(R^{2+\beta-\gamma/2}+R^{3-7\gamma/6}\right) \le CR^{2-\varepsilon}.
\]
Summing over dyadic annuli proves \eqref{eq:cm-mixed}.

For \(\phi\in C_c^\infty(\mathbb R)\) with \(\int\phi=0\), let \(\psi\in C_c^\infty(\mathbb R)\) satisfy \(\psi'=\phi\). Let \(\chi_R(\rho)=\chi(\rho/R)\), where \(\chi\in C_c^\infty((-2,2))\) and \(\chi=1\) on \([-1,1]\). For \(R=AH\), where \(A>2\) and \(\operatorname{spt}\psi\subset(-A/2,A/2)\), apply the first variation formula to
\[
\chi_R(r_h)Y_{\psi(\cdot/H)}
=\chi_R(r_h)\left[-\frac{\phi(z/H)}{H}X_h+\psi(z/H)e_3\right].
\]
Since
\[
\frac1H\int_C r_h\chi_R'(r_h)\phi(z/H) \mathrm d\mu_C=0,
\]
\eqref{eq:axial-variation} and \eqref{eq:uf-tail} imply
\[
\begin{aligned}
&\frac2H\int_T\chi_R(r_h)\phi(z/H)N_z^2\,\mathrm d\mu_T \\
={}&\frac1{H^2}\int_T\chi_R(r_h)\phi'(z/H)(X_h\cdot N)N_z\,\mathrm d\mu_T+O_\phi(H^{2\beta-1}).
\end{aligned}
\]
By \eqref{eq:wing-integral-tail} and \eqref{eq:cm-mixed},
\begin{equation}\label{eq:cm-zero-mean}
\left|\int\phi \mathrm d\nu_{T,H}\right| \le\frac{C_\phi}{H^2}\int_{T\cap B_{C_\phi H}} |X_h\cdot N| |N_z| \mathrm d\mu_T +C_\phi H^{2\beta-1} \le C_\phi H^{-\varepsilon}.
\end{equation}

Fix \(0\le\eta\in C_c^\infty(\mathbb R)\) with \(\int\eta=1\), and set
\[
a_T(H)=\int\eta \mathrm d\nu_{T,H}.
\]
For \(1\le s\le2\), taking \(\phi=\eta-s^{-1}\eta(\cdot/s)\) in \eqref{eq:cm-zero-mean} gives
\[
|a_T(H)-a_T(sH)|\le CH^{-\varepsilon}.
\]
Thus \(a_T(H)\to c(T)\ge0\), with \(|a_T(H)-c(T)|\le CH^{-\varepsilon}\). For arbitrary \(\phi\in C_c^\infty(\mathbb R)\), apply \eqref{eq:cm-zero-mean} to \(\phi-(\int\phi)\eta\). Then
\begin{equation}\label{eq:cm-smooth-mean}
\int\phi \mathrm d\nu_{T,H} =c(T)\int\phi+O_\phi(H^{-\varepsilon}).
\end{equation}
Approximating \(\mathbf1_{[-1,1]}\) from above and below in \eqref{eq:cm-smooth-mean} gives \eqref{eq:cm-mean}.

For fixed \(H\), \(a_{S_j}(H)\to a_S(H)\) by varifold convergence and \eqref{eq:wing-integral-tail}. Hence \(\limsup_j|c(S_j)-c(S)|\le CH^{-\varepsilon}\to0\) as \(H\to\infty\). Changing variables in \eqref{eq:cm-mean} proves \eqref{eq:uf-fluxlaws}.
\end{proof}

For \(T\in\mathcal C(C)\), let \(A_T\) be its second fundamental form.

\begin{proposition}\label{prop:cm-curvature}
There is \(L(C)<\infty\) such that every \(T\in\mathcal C(C)\) satisfies
\[
c(T)>0, \qquad \sup_T |A_T|\le\frac{L(C)}{c(T)}.
\]
\end{proposition}
\begin{proof}
Suppose that \(T\) has bounded curvature and \(c(T)=0\). Put \(\nu=z_\#(N_z^2\mu_T)\) and fix \(\kappa_j\downarrow0\). For fixed \(j\) and large \(H\), \(\nu((z-s,z+s))\le\nu([-2s,2s])\le2\kappa_js\) for \(|z|\le H\) and \(s\ge H\), by \eqref{eq:cm-mean}. The weak \((1,1)\) estimate for the maximal function of \(\nu\llcorner[-2H,2H]\) gives
\[
\left|\left\{z\in[-H,H]: \sup_{0<s<H}\frac{\nu((z-s,z+s))}{2s}>\kappa_j\right\}\right| \le\frac{C\nu([-2H,2H])}{\kappa_j}=o(H).
\]
Thus, for large \(H\), there is \(z_j\in[-H,H]\) such that
\[
\nu((z_j-s,z_j+s))\le2\kappa_js, \qquad s>0.
\]

Let \(\alpha=\angle(P_1,P_2)\). By Lemma~\ref{lem:tilt-iteration} with \(t=z_j\), \(\kappa=\kappa_j\), and \(h\) sufficiently large, there are \(q_j\in\mathbb{R}^2\) and \(s_j\in[A_\alpha \kappa_j,2A_\alpha \kappa_j)\) such that \(T_j=s_j^{-1}(T-(q_j,z_j))\) contains two points in \(B_C\) whose tangent planes make an angle bounded below by a positive constant independent of \(j\). Since \(T\) has bounded curvature and \(s_j\to0\),
\[
\sup_{T_j}|A_{T_j}|=s_j\sup_T|A_T|\longrightarrow0.
\]
After passing to a subsequence, \(T_j\) converges locally smoothly, possibly with multiplicity, to a union of parallel planes, a contradiction.

Suppose \(S_j\in\mathcal C(C)\) satisfy \(c(S_j)=c_0\ge0\) and \(\sup_{S_j}|A_{S_j}|\to\infty\). Choose \(p_j\in S_j\) and \(\lambda_j\to\infty\) such that, for \(T_j=\lambda_j(S_j-p_j)\),
\[
|A_{T_j}(0)|=1, \qquad \sup_{T_j\cap B_L}|A_{T_j}|\le2
\]
for every fixed \(L>0\) and all sufficiently large \(j\). After passing to a subsequence, \(T_j\to T_\infty\) locally smoothly, with \(|A_{T_\infty}(0)|=1\) and \(\sup_{T_\infty}|A_{T_\infty}|\le2\). Lemma~\ref{lem:cones} and Corollary~\ref{cor:cone-limit} imply \(|T_j|\to|T_\infty|\) as varifolds and \(T_\infty\in\mathcal C(C)\). Thus \(c(T_\infty)>0\), and Proposition~\ref{prop:uf-flux} gives
\[
\lambda_jc_0=c(T_j)\longrightarrow c(T_\infty)\in(0,\infty),
\]
which is impossible as \(\lambda_j\to\infty\).

Thus the curvature is uniformly bounded for each fixed value of \(c(T)\). In particular, every \(T\in\mathcal C(C)\) has bounded curvature and hence \(c(T)>0\). Since \(c(c(T)^{-1}T)=1\),
\[
c(T)\sup_T|A_T|\le L(C).\qedhere
\]
\end{proof}

Put \(E=\{\pm\eta_1,\pm\eta_2\}\), where \(\eta_i\) is the chosen unit normal to \(P_i\).

\begin{proposition}\label{prop:cm-marked}
For \(c_0,\varepsilon>0\), the set
\[
\{T\in\mathcal C(C): 0\in T, c(T)=c_0, \operatorname{dist}(N^{(T)}(0),E)\ge\varepsilon\}
\]
is compact for local smooth convergence. Here \(N^{(T)}(0)\) is a unit normal to \(T\) at the origin.
\end{proposition}
\begin{proof}
Let \(T_j\) be a sequence in this set. Choose \(R_j>0\) such that
\[
\frac{\mu_{T_j}(B_{R_j})}{\pi R_j^2}=\frac32, \qquad S_j=R_j^{-1}T_j.
\]
After passing to a subsequence, \(|S_j|\to V\) with \(\|V\|(B_1)=3\pi/2\). By Lemma~\ref{lem:cones}, \(\Theta(V,0)=1\) and \(V\) has density two at infinity. The convergence is smooth near \(0\), so \(\operatorname{dist}(N^{(V)}(0),E)\ge\varepsilon\). By Corollary~\ref{cor:cone-limit}, \((\eta_{0,R})_\#V\to C\). Since \(N^{(V)}(0)\notin E\), we have \(V=|S_\infty|\) with \(S_\infty\in\mathcal C(C)\) by Lemma~\ref{lem:cones}.

Propositions~\ref{prop:uf-flux} and~\ref{prop:cm-curvature} give
\[
\frac{c_0}{R_j}=c(S_j)\longrightarrow c(S_\infty)>0.
\]
Thus \(T_j=R_jS_j\to (c_0/c(S_\infty))S_\infty\) locally smoothly.
\end{proof}

\begin{proposition}\label{prop:cm-moving}
There are a smooth map \(p:\mathbb{R}\to e_3^\perp\) and constants \(C,R_0>0\) such that
\begin{equation}\label{eq:cm-center}
|p'|\le C, \qquad |p(z)-p(s)|\le C(1+|z-s|)^{\beta}.
\end{equation}
For \(i=1,2\), put \(a_i=\tau_i\cdot p\), \(b_i=\eta_i\cdot p\). The surface \(M\) contains four disjoint graphs
\begin{equation}\label{eq:cm-movingwings}
X_i(r,z)=r\tau_i+u_i(r,z)\eta_i+z e_3, \qquad d=|r-a_i(z)|>R_0,
\end{equation}
with
\begin{equation}\label{eq:cm-movingbounds}
|u_i-b_i(z)|+d|\nabla u_i|+d^2|\nabla^2 u_i|\le Cd^{\beta}.
\end{equation}
Their complement in \(M\) is contained in \(\{|X_h-p(z)|\le C\}\) and has area at most \(C\) in every slab \(\{s<z<s+1\}\).
\end{proposition}
\begin{proof}
Choose \(\varepsilon_0>0\) so that the \(2\varepsilon_0\)-neighborhoods of the points of \(E\) are disjoint. Each level \(z=n\), \(n\in\mathbb Z\), contains \(q_n=p_n+ne_3\in M\) with \(\operatorname{dist}(N(q_n),E)\ge\varepsilon_0\). At a critical point of \(z|_M\), \(N=\pm e_3\). At a regular level \(z=n\), a sufficiently large circle meets the graphs of \(u_{M,1}(\cdot,n)\) and \(u_{M,2}(\cdot,n)\) twice each, in alternating order. Disjoint arcs cannot join both opposite pairs, so an arc joins graphs over different planes and contains a point \(q_n\) with \(\operatorname{dist}(N(q_n),E)\ge\varepsilon_0\).

By Proposition~\ref{prop:cm-marked}, \eqref{eq:uf-tail} holds for \(T=M-q_n\), with \(C>0,R_*>2\) independent of \(n\). The four graphs of \(M-q_n\), translated by \(q_n\), cover
\[
M\cap\{|X_h-p_n|>R_*, |z-n|\le|X_h-p_n|/2\}.
\]
On each graph over \(q_n+P_i\), \(\operatorname{dist}(N,\{\pm\eta_i\})<\varepsilon_0\). Since no \(q_m\) lies on these graphs, \(|p_n-p_m|\le R_*+2|n-m|\). Take \(r=A(R_*+|n-m|)\) with \(A\) sufficiently large. For \(i=1,2\), the point \(X\) on the graph over \(q_n+P_i\) with \(z=n\) and \((X_h-p_n)\cdot\tau_i=r\) also lies on the graph over \(q_m+P_i\). By \eqref{eq:uf-tail},
\[
|\eta_i\cdot(p_n-p_m)| \le |\eta_i\cdot(X_h-p_n)|+|\eta_i\cdot(X_h-p_m)| \le Cr^\beta.
\]
Since \(\eta_1,\eta_2\) span \(e_3^\perp\), \(|p_n-p_m|\le C(1+|n-m|)^\beta\).

Choose a smooth \(p:\mathbb R\to e_3^\perp\) satisfying \eqref{eq:cm-center}, with \(|p(z)-p_n|\le C\) for \(|z-n|<1\). For fixed \(i\) and large \(R_0\), the graphs over \(q_n+P_i\) restricted to \(|z-n|<1\) and \(d>R_0\) agree on overlaps and define \(u_i\) in \eqref{eq:cm-movingwings}. By \eqref{eq:uf-tail}, \(u_i\) satisfies \eqref{eq:cm-movingbounds}. The part of \(M\cap\{s<z<s+1\}\) outside the graphs of \(u_1\) and \(u_2\) lies in \(B_C(p(s)+se_3)\) by \eqref{eq:cm-center}, so its area is at most \(C\).
\end{proof}

\subsection{Bounded heights}\label{sec:bounded-wings}

Set \(K(x)=x/(2\pi|x|^2)\).

\begin{lemma}\label{lem:kernel-comparison}
Let \(0<\beta<1/2\) and \(C_0>0\). Suppose \(a\in C(\mathbb R)\) and a measurable \(\mathcal F:\mathbb R^2\to\mathbb R^2\) satisfy
\[
|a(z)-a(s)|\le C_0(1+|z-s|)^\beta, \qquad |\mathcal F(r,s)|\le C_0(1+|r-a(s)|)^{2\beta-2}.
\]
There are \(C,R_*>0\), depending only on \(\beta,C_0\), such that for \(R\ge R_*\) and \(z\in\mathbb R\),
\begin{equation}\label{eq:fa-collapse-error}
\int_{\mathbb R^2} |K(a(z)\pm R-r,z-s)-K(\pm R,z-s)| |\mathcal F(r,s)| \mathrm dr \mathrm ds \le CR^{2\beta-1}.
\end{equation}
\end{lemma}
\begin{proof}
Fix \(z\), and write \(h=z-s\), \(\delta(h)=a(z)-a(z-h)\), \(\xi=r-a(s)\), and \(L=R+|h|\). For \(R_*\ge1\) large enough, \(|\delta(h)|\le L/4\) whenever \(R\ge R_*\).

Using \(|DK(x)|\le C|x|^{-2}\) and the bound on \(\mathcal F\),
\[
\begin{aligned}
&\int_{|\xi|\le L/4}|K(\pm R+\delta-\xi,h)-K(\pm R,h)|\,|\mathcal F(a(s)+\xi,s)|\,\mathrm d\xi \\
\le{}& CL^{-2}(|\delta|+L^{2\beta})\le CL^{2\beta-2}.
\end{aligned}
\]
Since \(|K(\pm R,h)|\le C/L\),
\[
\int_{|\xi|>L/4}|K(\pm R,h)|\,|\mathcal F(a(s)+\xi,s)|\,\mathrm d\xi
\le \frac C L\int_{|\xi|>L/4}(1+|\xi|)^{2\beta-2}\mathrm d\xi\le CL^{2\beta-2}.
\]
We also have \(|K(\pm R+\delta-\xi,h)|\le C/L\) when \(|h|\ge R\) or \(|\xi|>4R\). Since \(2\beta<1\),
\[
\int_{\mathbb R}L^{2\beta-2}\,\mathrm dh\le CR^{2\beta-1}.
\]

On \(\{|h|<R, L/4<|\xi|\le4R\}\), \(|\mathcal F(a(s)+\xi,s)|\le CR^{2\beta-2}\). With \(q=\pm R+\delta(h)-\xi\),
\[
\begin{aligned}[b]
&\int_{|h|<R}\int_{L/4<|\xi|\le4R}|K(\pm R+\delta(h)-\xi,h)|\,|\mathcal F(a(s)+\xi,s)|\,\mathrm d\xi\,\mathrm dh \\
\le{}& CR^{2\beta-2}\int_{B_{6R}}|K(q,h)|\,\mathrm dq\,\mathrm dh\le CR^{2\beta-1}.
\end{aligned}\qedhere
\]
\end{proof}

\begin{proposition}\label{prop:fa-moving-rigidity}
For some \(C>R_1>0\), the functions \(u_i\) in Proposition~\ref{prop:cm-moving} are defined for \(|r|>R_1\) and \(z\in\mathbb R\), and satisfy
\begin{equation}\label{eq:fa-boundedwings}
|u_i|+|r||\nabla u_i|+r^2|\nabla^2 u_i|\le C.
\end{equation}
Their graphs cover \(M\cap\{|X_h|>2C\}\). There are constants \(m_1,m_2\) such that, for every \(R>R_1\),
\begin{equation}\label{eq:fa-wingmean}
\sup_{z\in\mathbb{R}}\left| \frac{u_i(R,z)+u_i(-R,z)}2-m_i\right| \le CR^{2\beta-1}, \qquad i=1,2.
\end{equation}
\end{proposition}

\begin{proof}
\emph{Step 1. \(\Delta u_i=\operatorname{div}\mathcal F_i\) on \(\mathbb R^2\).} Let \(M_0\) be the complement in \(M\) of the four graphs on \(d=|r-a_i(z)|>R_0\). Extend \(u_i\) linearly in \(r\) across \(d\le R_0\), so that
\begin{equation}\label{eq:fa-centertrace}
\|\nabla u_i\|_{L^\infty(\mathbb R^2)}\le C.
\end{equation}
For \(p\in\mathbb R^2\), set
\begin{equation}\label{eq:graph-fluxes}
W(p)=\sqrt{1+|p|^2}, \qquad G(p)=\frac{p}{W(p)}, \qquad B(p)=W(p)I-\frac{p\otimes p}{W(p)}.
\end{equation}
For \(d>R_0\), set
\begin{equation}\label{eq:fa-fluxes}
\mathcal F_i=\nabla u_i-G(\nabla u_i), \qquad B_i=B(\nabla u_i)-I.
\end{equation}
By the minimal surface equation and \eqref{eq:cm-movingbounds},
\[
\operatorname{div}B_i=0, \qquad |B_i|\le Cd^{2\beta-2}, \qquad |\mathcal F_i|\le Cd^{3\beta-3}.
\]
Put
\begin{equation}\label{eq:fa-movingflux}
g_i^\pm=\left[(1,-a_i')\cdot G(\nabla u_i)\right]_{r=a_i(z)\pm R_0}.
\end{equation}
Since \(\partial_zB_i^{rz}=-\partial_rB_i^{rr}\) for \(d>R_0\),
\[
\frac{\mathrm d}{\mathrm dz} \int_{|r-a_i(z)|>R_0}B_i^{rz}(r,z) \mathrm dr =\left[B_i^{rr}-a_i'B_i^{rz}\right]_{r=a_i(z)-R_0}^{r=a_i(z)+R_0}.
\]
For \(\varphi\in C_c^\infty(\mathbb R)\), multiply \(\Delta_MX_h=0\) by \(\varphi(z)\) and integrate by parts on \(M_0\). Then
\[
\begin{aligned}
&-\int_{M_0}\varphi'N_z(N-N_ze_3) \mathrm d\mu \\
={}& \sum_i\tau_i\int\varphi \left[B_i^{rr}-a_i'B_i^{rz}\right]_{r=a_i(z)-R_0}^{r=a_i(z)+R_0} \mathrm dz +\sum_i\eta_i\int\varphi(g_i^+-g_i^-) \mathrm dz.
\end{aligned}
\]
Define signed measures \(Q_i\) by
\begin{align}
&\sum_i\eta_iQ_i \label{eq:fa-Q-definition}\\
={}& -z_\#(N_z(N-N_ze_3)\mu\llcorner M_0) +\sum_i\tau_i\left(\int_{|r-a_i(z)|>R_0} B_i^{rz}(r,z) \mathrm dr\right)\mathrm dz.\notag
\end{align}
Then
\begin{equation}\label{eq:fa-Q}
Q_i'=-(g_i^+-g_i^-)
\end{equation}
in distributions. Since \(g_i^\pm\) are bounded, \(Q_i\) have Lipschitz densities, which we also denote by \(Q_i\). By Proposition~\ref{prop:cm-moving} and \(2\beta<1\), \(\int_s^{s+1}|Q_i(z)|\,\mathrm dz\le C\) for all \(s\in\mathbb R\), so \(|Q_i(z)|\le C\).

For \(d\le R_0\), set
\[
\mathcal F_i^z =\partial_z u_i-\frac{Q_i}{2R_0},\quad \mathcal F_i^r =\partial_r u_i-g_i^- -\frac{r-a_i+R_0}{2R_0}(g_i^+-g_i^-) -\frac{a_i'Q_i}{2R_0}.
\]
For \(d<R_0\), \eqref{eq:fa-Q} gives
\[
\operatorname{div}(\nabla u_i-\mathcal F_i) =\frac{g_i^+-g_i^-+Q_i'}{2R_0}=0.
\]
On \(r=a_i(z)\pm R_0\),
\[
(1,-a_i')\cdot(\nabla u_i-\mathcal F_i)=g_i^\pm.
\]
By \eqref{eq:fa-movingflux}, the normal component of \(\nabla u_i-\mathcal F_i\) agrees across \(r=a_i(z)\pm R_0\). Hence
\begin{equation}\label{eq:fa-PDE}
\Delta u_i=\operatorname{div}\mathcal F_i\qquad \text{in distributions on }\mathbb R^2.
\end{equation}
Moreover,
\begin{equation}\label{eq:fa-mass}
\|\mathcal F_i\|_{L^\infty(\mathbb R^2)}+ \sup_z\int_{\mathbb R}|\mathcal F_i(r,z)| \mathrm dr\le C, \qquad \int_{B_\rho(q)}|\mathcal F_i|\le C\min(\rho^2,\rho).
\end{equation}

\emph{Step 2. Estimate~\eqref{eq:fa-cancellation}.} Put \(f_i(z)=\int_{\mathbb R}\mathcal F_i^z(r,z) \mathrm dr\), so \(\|f_i\|_{L^\infty(\mathbb R)}\le C\) by \eqref{eq:fa-mass}. For every interval \(I=[a,b]\),
\begin{equation}\label{eq:fa-cancellation}
\left|\int_I f_i(z) \mathrm dz\right|\le C(1+|I|)^{2\beta}.
\end{equation}
Choose \(\psi\in C_c^\infty((a-1,b+1))\) with \(\psi|_I=1\), \(0\le\psi\le1\), and \(|\psi'|\le C\). Then \(\int\psi f_i=\int_I f_i+O(1)\).

We may assume that \(p((a+b)/2)=0\). With \(R=A(1+|I|)\) for a fixed sufficiently large \(A\), we have \(\lvert p(z)\rvert\le CR^\beta\) and \(|a_i(z)|+R_0<R/2\) on \((a-1,b+1)\). Let \(\chi_R(r)=\chi(|r|/R)\), where \(\chi\in C_c^\infty((-2,2))\) and \(\chi=1\) on \([-1,1]\). Set \(\zeta_R=1\) on \(M_0\) and \(\zeta_R(X_i(r,z))=\chi_R(r)\) for \(d>R_0\). For \(z\in(a-1,b+1)\) and \(R\le|r|\le2R\),
\begin{equation}\label{eq:fa-cutoff}
|u_i|\le CR^\beta, \qquad |\nabla u_i|\le CR^{\beta-1}.
\end{equation}

Since \(\Delta_Mz=0\) and \(\chi_R\) is independent of \(z\),
\[
0=\int_M\nabla_Mz\cdot\nabla_M\left(\zeta_R\psi(z)X_h\right) \mathrm d\mu -\sum_i\iint\chi_R\partial_z[\psi(r\tau_i+u_i\eta_i)] \mathrm dr \mathrm dz.
\]
For \(d>R_0\), \(B(\nabla u_i)\nabla u_i=\nabla u_i-\mathcal F_i\), so
\[
\begin{aligned}
&\nabla_Mz\cdot\nabla_M(\zeta_R\psi X_h) \mathrm d\mu \\
={}& \left[ \chi_R\psi\left(B_i^{rz}\tau_i+ (\partial_z u_i-\mathcal F_i^z)\eta_i\right) +\chi_R\psi'(1+B_i^{zz})X_h +\chi_R'\psi B_i^{rz}X_h\right] \mathrm dr \mathrm dz.
\end{aligned}
\]
On \(M_0\),
\[
\nabla_Mz\cdot\nabla_M(\psi X_h) =\psi'X_h(1-N_z^2)-\psi N_z(N-N_ze_3).
\]
By \eqref{eq:fa-cutoff}, \(\left|\iint_{d>R_0}\chi_R'\psi B_i^{rz}X_h\,\mathrm dr\mathrm dz\right|\le CR^{2\beta}\). Replacing \(\chi_R\) by \(1\) in the terms containing \(B_i^{rz}\) and \(\mathcal F_i^z\) gives an error bounded by
\[
CR\int_R^\infty\left(s^{2\beta-2}+s^{3\beta-3}\right) \mathrm ds\le CR^{2\beta}.
\]
By \eqref{eq:fa-Q-definition} and \(\int_{d\le R_0}(\partial_z u_i-\mathcal F_i^z) \mathrm dr=Q_i\),
\[
\begin{aligned}
\sum_i\eta_i\int\psi f_i \mathrm dz &=\int_{M_0}\psi'X_h(1-N_z^2) \mathrm d\mu \\
&-\sum_i\iint_{d\le R_0}\psi'(r\tau_i+u_i\eta_i) \mathrm dr \mathrm dz \\
&+\sum_i\iint_{d>R_0}\chi_R\psi' (r\tau_i+u_i\eta_i)B_i^{zz} \mathrm dr \mathrm dz +O(R^{2\beta}).
\end{aligned}
\]
By Proposition~\ref{prop:cm-moving}, the first two integrals are \(O(R^\beta)\), since \(\operatorname{spt}\psi'\subset[a-1,a]\cup[b,b+1]\) and \(|p|\le CR^\beta\) there.
For \(z\in\operatorname{spt}\psi'\) and \(d>R_0\), \(|X_h|\le C(d+R^\beta)\) and \(|B_i|\le Cd^{2\beta-2}\). Since \(\eta_1,\eta_2\) are linearly independent,
\[
\left|\int_I f_i(z) \mathrm dz\right| \le C\left(R^\beta+\int_{R_0}^{CR}(s+R^\beta)s^{2\beta-2} \mathrm ds\right) +CR^{2\beta} \le CR^{2\beta}.
\]

\emph{Step 3. Estimates~\eqref{eq:fa-boundedwings} and~\eqref{eq:fa-wingmean}.} Fix \(i\) and omit the subscript \(i\). Put \(x=(r,z)\) and recall that \(K(x)=x/(2\pi|x|^2)\). For \(x,x_*\in\mathbb{R}^2\),
\begin{equation}\label{eq:fa-potential-normalization}
u(x)-u(x_*)= \int_{\mathbb{R}^2}[K(x-y)-K(x_*-y)]\cdot\mathcal F(y)\mathrm d A(y).
\end{equation}
By \eqref{eq:fa-mass},
\[
\int_{\mathbb R^2}|K(x-y)-K(x_*-y)|\,|\mathcal F(y)|\,\mathrm d A(y)\le C\log(2+|x-x_*|).
\]
By \eqref{eq:cm-center} and \eqref{eq:cm-movingbounds}, \(\lvert u(r,z)\rvert\le C(1+|r|+|z|)^{\beta}\), so the difference between the two sides of \eqref{eq:fa-potential-normalization} has sublinear growth. It is harmonic on \(\mathbb R^2\) by \eqref{eq:fa-PDE} and zero at \(x=x_*\), hence identically zero.

For \(h\ge0\), \eqref{eq:fa-cancellation} gives \(\left|\int_0^h f(z\pm t) \mathrm dt\right|\le C(1+h)^{2\beta}\). Since \(2\beta<1\), integration by parts gives
\[
\int_0^\infty\frac{hf(z\pm h)}{R^2+h^2} \mathrm dh =-\int_0^\infty\frac{R^2-h^2}{(R^2+h^2)^2} \left(\int_0^h f(z\pm t) \mathrm dt\right)\mathrm dh.
\]
Thus for \(R\ge1\),
\begin{equation}\label{eq:fa-fpotentialdecay}
\sup_z\left|\int_\mathbb{R} \frac{(z-s)f(s)}{R^2+(z-s)^2} \mathrm{d} s\right| \le C\int_0^\infty\frac{(1+h)^{2\beta}}{R^2+h^2} \mathrm{d} h \le CR^{2\beta-1}.
\end{equation}
Put \(g(s)=\int_{\mathbb R}\mathcal F^r(r,s) \mathrm dr\), so \(\|g\|_{L^\infty(\mathbb R)}\le C\) by \eqref{eq:fa-mass}.

By \eqref{eq:fa-mass} and \(\mathcal F=O(d^{3\beta-3})\) for \(d>R_0\), \(|\mathcal F(r,s)|\le C(1+|r-a(s)|)^{2\beta-2}\). By Lemma~\ref{lem:kernel-comparison}, replacing \(K(a(z)\pm R-r,z-s)\) by \(K(\pm R,z-s)\) in \eqref{eq:fa-potential-normalization} introduces an error \(O(R^{2\beta-1})\) for \(R\ge R_*\), where \(R_*>R_0+1\). With \(x_*\) fixed, \eqref{eq:fa-fpotentialdecay} gives
\begin{equation}\label{eq:fa-collapsedpotential}
u(a(z)\pm R,z)=m\pm\frac1{2\pi}\int_{\mathbb R} \frac{Rg(s)}{R^2+(z-s)^2} \mathrm ds +O(R^{2\beta-1}),
\end{equation}
where \(m\in\mathbb R\) and the error is uniform in \(z\). Thus \(|u_i|\le C\) on \(d\ge R_*\) for \(i=1,2\). At \(d=R_*\), \eqref{eq:cm-movingbounds} gives \(|b_i|=|\eta_i\cdot p|\le C\). Since \(\eta_1,\eta_2\) span \(e_3^\perp\), \(|p|\le C\). Equation~\eqref{eq:cm-movingbounds} then gives \(|u_i|\le C\) for all \(d>R_0\).

Choose \(R_1>4(R_*+\sup_{i,z}|a_i(z)|)\). For \(|r|>R_1/2\), interior estimates on \(B_{|r|/2}(r,z)\) give \eqref{eq:fa-boundedwings}. By Proposition~\ref{prop:cm-moving}, the graphs on \(|r|>R_1\) cover \(M\cap\{|X_h|>2C\}\) for some \(C>R_1\). Since \(|a_i|\le C\), \eqref{eq:fa-boundedwings} and \eqref{eq:fa-collapsedpotential} give for \(R>R_1\)
\[
\begin{aligned}[b]
\frac{u_i(R,z)+u_i(-R,z)}2
&=\frac{u_i(a_i(z)+R,z)+u_i(a_i(z)-R,z)}2+O(R^{-1}) \\
&=m_i+O(R^{2\beta-1}).
\end{aligned}\qedhere
\]
\end{proof}

Let \(\mathscr H(M)\) be the closure of \(\{M-ze_3:z\in\mathbb R\}\) with respect to local smooth convergence.

\begin{proposition}\label{prop:fa-hull}\label{prop:fa-wingmean}
The set \(\mathscr H(M)\subset\mathcal C(C)\) is compact. Every \(S\in\mathscr H(M)\) satisfies \(c(S)=c(M)\) and has four graphs over \(P_1,P_2\) on \(|r|>R_1\). These graphs cover \(S\cap\{|X_h|>2C\}\) and satisfy \eqref{eq:fa-boundedwings} and \eqref{eq:fa-wingmean} with the same constants \(R_1,C,m_1,m_2\).
\end{proposition}
\begin{proof}
Given \((z_j)\subset\mathbb R\), choose \(q_j\in M\cap\{z=z_j\}\) with \(\operatorname{dist}(N(q_j),E)\ge\varepsilon_0\), as in the proof of Proposition~\ref{prop:cm-moving}. By \eqref{eq:fa-boundedwings}, \(|q_j-z_je_3|\le C\).
By Proposition~\ref{prop:cm-marked} applied to \(M-q_j\), we may assume that \(M-z_je_3\to S\in\mathcal C(C)\) locally smoothly. By continuity of \(c\),
\[
c(S)=\lim_{j\to\infty}c(M-z_je_3)=c(M).
\]

The graph functions \(u_i(r,z+z_j)\) of \(M-z_je_3\) converge smoothly on compact subsets of \(\{|r|>R_1\}\) by \eqref{eq:fa-boundedwings}. Every point of \(M-z_je_3\) with \(|X_h|>2C\) lies on a graph where \(r^2=|X_h|^2-u_i(r,z+z_j)^2>C^2>R_1^2\). The limiting graphs therefore cover \(S\cap\{|X_h|>2C\}\) and satisfy \eqref{eq:fa-boundedwings} and \eqref{eq:fa-wingmean} with the same \(C,m_1,m_2\).
\end{proof}

\subsection{Reflection symmetry}\label{sec:alexandrov-reflection}

Translate so that \(m_1=m_2=0\). Write \(X=(x,y,z)\), with the \(x\)- and \(y\)-axes along the angle bisectors of the lines \(\mathbb R\tau_1\) and \(\mathbb R\tau_2\).

\begin{theorem}\label{thm:alex-hull}\label{thm:original-alex}
Every \(T\in\mathscr H(M)\) is invariant under reflection in \(x=0\) and \(y=0\). On each side of either plane, orthogonal projection of \(T\) onto that plane is a diffeomorphism onto its image.
\end{theorem}
\begin{proof}
Fix \(T\in\mathscr H(M)\). By \eqref{eq:fa-boundedwings}, on \(|y|>Y\) the surface consists of two graphs \(x=g_\pm(y,z)\), with
\begin{equation}\label{eq:ax-far}
|g_+-a|y||+|g_-+a|y||\le B, \qquad |\nabla g_+|+|\nabla g_-|\le B,
\end{equation}
where \(a>0,Y,B\) are independent of \(T\). On \(T\cap\{|y|\le Y\}\), \(|x|\) is uniformly bounded.

For sufficiently large \(L\), the circle \(\{x^2+y^2=L^2, z=0\}\) meets every \(S\in\mathscr H(M)\) transversely once in each quadrant. Since \(S\) is connected and properly embedded, \(\mathbb R^3\setminus S\) has two components, which alternate along the circle. Let \(\Omega_X(S),\Omega_Y(S)\) be those containing \((\pm L,0,0)\) and \((0,\pm L,0)\), respectively. Direct \(N\) toward \(\Omega_Y(S)\).
For \(|y|>Y\), the region \(g_-<x<g_+\) belongs to \(\Omega_Y(T)\).

Let \(\mathscr R_\lambda\) be reflection in \(x=\lambda\), and write \(T_\lambda^+=T\cap\{x>\lambda\}\). For large \(\lambda\), \(T_\lambda^+\) lies on the graph of \(g_+\). By \eqref{eq:ax-far}, \(g_-<2\lambda-g_+<g_+\) whenever \(g_+>\lambda\), so \(\mathscr R_\lambda(T_\lambda^+)\subset\Omega_Y(T)\). Define
\begin{equation}\label{eq:ax-stop}
\lambda_0=\inf\{\lambda\ge0: \mathscr R_\mu(T_\mu^+)\subseteq\overline{\Omega_Y(T)} \text{ for every }\mu>\lambda\}.
\end{equation}

Suppose \(\lambda_0>0\). Choose \(0<\lambda_j\uparrow\lambda_0\) and \(q_j\in T_{\lambda_j}^+\) such that \(\mathscr R_{\lambda_j}(q_j)\in\Omega_X(T)\). By \eqref{eq:ax-far}, \(g_-<0<\lambda_0<g_+\) on \(|y|\ge Y\) for sufficiently large \(Y\). On \(|y|>Y\), the function \(w_j=2\lambda_j-g_+-g_-\) is bounded. For \(R>2\), choose \(\eta_R\in C_c^\infty(\mathbb R^2)\) with \(\eta_R=1\) on \(B_R\) and \(\int_{\mathbb R^2}|\nabla\eta_R|^2\le C/\log R\). If \(w_j\ge0\) on \(|y|=Y\), then by the minimal surface equations for \(2\lambda_j-g_+\) and \(g_-\) and the bound \(|\nabla g_+|+|\nabla g_-|\le B\),
\[
\int_{\{|y|>Y\}}\eta_R^2|\nabla (w_j)_-|^2 \le C\int_{\mathbb R^2}|\nabla \eta_R|^2\le\frac C{\log R}
\]
and hence \(w_j\ge0\) on \(|y|>Y\) after letting \(R\to\infty\). If \(|y(q_j)|>Y\), then \(w_j(y(q_j),z(q_j))<0\), so \(w_j<0\) at a point with \(|y|=Y\). Hence \(q_j\) may be chosen with \(|y(q_j)|\le Y\) and \(|x(q_j)|\le C\).

For \(\alpha=X,Y\), \(\Omega_\alpha(S-ze_3)=\Omega_\alpha(S)-ze_3\). If \(S_j\to S\) in \(\mathscr H(M)\), then \(K\subset\Omega_\alpha(S_j)\) for every \(K\Subset\Omega_\alpha(S)\) and all sufficiently large \(j\). Hence if \(p_j\to p\) with \(p_j\in\overline{\Omega_\alpha(S_j)}\), then
\begin{equation}\label{eq:ax-colors}
p\in\overline{\Omega_\alpha(S)}.
\end{equation}

By Proposition~\ref{prop:fa-hull}, after passing to a subsequence,
\[
T_j:=T-z(q_j)e_3 \longrightarrow S\in\mathscr H(M),\quad p_j:=q_j-z(q_j)e_3 \longrightarrow p\in S\cap\{x\ge\lambda_0\}.
\]
For \(\mu\ge\lambda_0\), \eqref{eq:ax-colors} gives
\begin{equation}\label{eq:ax-cap}
\mathscr R_\mu(S\cap\{x>\mu\})\subseteq\overline{\Omega_Y(S)}.
\end{equation}
If \(x(p)>\lambda_0\), then by \eqref{eq:ax-colors},
\[
\mathscr R_{\lambda_0}(p)\in \overline{\Omega_X(S)}\cap\overline{\Omega_Y(S)}=S.
\]
By the maximum principle, \(S=\mathscr R_{\lambda_0}(S)\) near \(\mathscr R_{\lambda_0}(p)\).

If \(x(p)=\lambda_0\), letting \(\mu\uparrow x(q)\) in \eqref{eq:ax-cap} gives \(N_x(q)\le0\) for \(q\in S\cap\{x>\lambda_0\}\). Since \(x|_S\) has no local maximum, \(N_x(p)\le0\). On the other hand,
\[
\mathscr R_{\lambda_j}(p_j)=p_j-2(x(p_j)-\lambda_j)e_x\in\Omega_X(T_j).
\]
Since \(x(p_j)-\lambda_j\to0\) and \(T_j\to S\) smoothly, \(N_x(p)\ge0\). Thus \(N_x(p)=0\), and by Hopf's lemma, \(S=\mathscr R_{\lambda_0}(S)\) near \(p\).

By properness and unique continuation, \(\mathscr R_{\lambda_0}(S)=S\).
The linear part of \(\mathscr R_{\lambda_0}\) permutes \(P_1,P_2\). By \eqref{eq:fa-boundedwings}, fixed shifts of \(r\) change the averages in \eqref{eq:fa-wingmean} by \(O(R^{-1})\). Since \(m_1=m_2=0\), the averages in \eqref{eq:fa-wingmean} for \(\mathscr R_{\lambda_0}(S)\) tend to \(2\lambda_0\eta_i\cdot e_x\). Since \(\mathscr R_{\lambda_0}(S)=S\),
\[
0=2\lambda_0\eta_i\cdot e_x, \qquad i=1,2.
\]
Thus \(\lambda_0=0\). Similarly, \(\mathscr R_0(T\cap\{x<0\})\subseteq\overline{\Omega_Y(T)}\). For sufficiently large \(|y|\), \(-g_+\ge g_-\) and \(-g_-\le g_+\), so \(g_++g_-=0\). By unique continuation, \(\mathscr R_0(T)=T\).

For \(p\in T\cap\{x>0\}\) and \(0<\mu<x(p)\), \(\mathscr R_\mu(p)\in\overline{\Omega_Y(T)}\). Letting \(\mu\uparrow x(p)\) gives \(N_x(p)\le0\). Since the Gauss map of \(T\) is open, \(N_x<0\) on \(T\cap\{x>0\}\). At every intersection with \(T\cap\{x>0\}\), a line directed by \(e_x\) passes from \(\Omega_Y(T)\) to \(\Omega_X(T)\), and hence has at most one such intersection. Projection onto \(x=0\) is therefore a diffeomorphism onto its image. The case \(x<0\) follows by reflection in \(x=0\). Interchanging \(x,y\) and \(\Omega_X,\Omega_Y\) proves both claims for \(y=0\).
\end{proof}

\section{Periodicity and Scherk rigidity}\label{sec:scherk-rigidity}

\subsection{Half-plane coordinates}\label{sec:half-plane-coordinates}

Fix \(S\in\mathscr H(M)\), with \(N\) directed toward \(\Omega_Y(S)\). Both \(\Omega_X(S)\) and \(\Omega_Y(S)\) are invariant under reflection in \(x=0\) and \(y=0\). By the proof of Theorem~\ref{thm:alex-hull} and reflection,
\begin{equation}\label{eq:hp-signs}
\operatorname{sgn}N_x=-\operatorname{sgn}x, \qquad \operatorname{sgn}N_y=\operatorname{sgn}y.
\end{equation}
Set
\[
Q=S\cap\{x>0, y>0\}.
\]
Closures and boundaries of \(Q\) are taken in \(S\).

\begin{lemma}\label{lem:hp-saddles}
The critical set of \(z|_S\) is \(S\cap\{x=y=0\}\) and is locally finite. Near each \(p=(0,0,z_0)\) in this set, \(S\) is the graph
\begin{equation}\label{eq:hp-saddle}
z-z_0=h(x,y)=A(x^2-y^2)+O((x^2+y^2)^2), \qquad A\ne0.
\end{equation}
Moreover, \(Q\cap\{z=z_0\}\) contains an arc ending at \(p\).
\end{lemma}
\begin{proof}
By \eqref{eq:hp-signs}, \(dz|_S=0\) exactly when \(x=y=0\). Thus \(S\) meets the \(z\)-axis transversely in a locally finite set. Near \(p=(0,0,z_0)\in S\), \(S\) is a graph \(z-z_0=h(x,y)\), with \(h\) even in each variable and \(\Delta h(0)=0\). \(h_x\) solves the linearized minimal surface equation, with \(h_x=0\) on \(x=0\) and \(h_x\ne0\) for \(x>0\). By the Hopf lemma, \(h_{xx}(0)\ne0\), so \eqref{eq:hp-saddle} holds. Thus \(Q\cap\{z=z_0\}\) contains an arc \(y=x+O(x^3)\) ending at \(p\).
\end{proof}

Let \(\tau=(\cos\theta,\sin\theta,0)\), \(0<\theta<\pi/2\), point into the half-plane of \(C\) in \(\{x>0,y>0\}\), and put \(\eta=(-\sin\theta,\cos\theta,0)\). By Propositions~\ref{prop:fa-moving-rigidity} and~\ref{prop:fa-hull}, a graph \(W\subset Q\) is parametrized for \(r>R_*\) and \(z\in\mathbb R\) by
\begin{equation}\label{eq:hp-wing}
X(r,z)=r\tau+u(r,z)\eta+z e_3, \qquad |u|+r|\nabla u|\le C.
\end{equation}
On all of \(S\),
\begin{equation}\label{eq:hp-cross-bounds}
|x|\le(\cot\theta)|y|+K, \qquad |y|\le(\tan\theta)|x|+K,
\end{equation}
with \(R_*,C,K\) independent of \(S\in\mathscr H(M)\).

\begin{lemma}\label{lem:hp-slices}\label{lem:hp-topology}
\(\overline Q=S\cap\{x\ge0, y\ge0\}\). For every \(t\in\mathbb R\), \(Q\cap\{z=t\}\) is a single smooth arc with endpoint \(\gamma(t)=(a(t),b(t),t)\), along which \(x\) and \(y\) strictly increase to \(+\infty\). The functions \(a,b:\mathbb R\to[0,\infty)\) are continuous and bounded and satisfy \(ab=0\). The map
\[
\overline Q\longrightarrow[0,\infty)\times\mathbb R, \qquad (x,y,z)\longmapsto(y-b(z),z),
\]
is a homeomorphism taking \(Q\) onto \((0,\infty)\times\mathbb R\).
\end{lemma}
\begin{proof}
Each level \(z=t\) meets \(W\subset Q\). By \eqref{eq:hp-signs}, each component of \(Q\cap\{z=t\}\) is a smooth curve that can be parametrized with \(x'=N_y>0\) and \(y'=-N_x>0\). By properness and \eqref{eq:hp-cross-bounds}, both \(x\) and \(y\) are unbounded above on each component. If there were two components, a line parallel to the \(x\)-axis would meet \(Q\) twice, contrary to Theorem~\ref{thm:alex-hull}. Since \(x,y\) are bounded below, the arc has an endpoint \(\gamma(t)=(a(t),b(t),t)\), which lies in \(\partial Q\) by properness.

By \eqref{eq:hp-signs}, \(S\) has local coordinates \((x,z)\) at points with \(x=0<y\), and \((y,z)\) at points with \(y=0<x\). These charts and Lemma~\ref{lem:hp-saddles} at the axis show that every \(p\in S\cap\{x\ge0, y\ge0\}\) belongs to \(\overline{Q\cap\{z=z(p)\}}\). Thus \(\partial Q=\gamma(\mathbb R)\).

By \eqref{eq:hp-cross-bounds}, \(a+b\le K\), so \(z:\partial Q\to\mathbb R\) is a proper bijection. Its inverse \(\gamma\) is therefore continuous. The map \((x,y,z)\mapsto(y-b(z),z)\) is a proper continuous bijection, hence a homeomorphism.
\end{proof}

\begin{theorem}\label{thm:halfplane}\label{rem:hp-scope}
There is a unique conformal bijection
\begin{equation}\label{eq:hp-coordinate}
w:Q\longrightarrow\mathbb H_+ :=\{\zeta\in\mathbb C:\operatorname{Im}\zeta>0\}, \qquad \operatorname{Re}w=z.
\end{equation}
It extends to a homeomorphism \(\overline Q\to\overline{\mathbb H_+}\), with \(w(\gamma(t))=t\).
\end{theorem}
\begin{proof}
Since \(\Delta_S z=0\) and \(Q\) is simply connected, \(z\) has a harmonic conjugate \(s\) with \(ds=*dz\), where \(dz+i*dz\) is of type \((1,0)\) (cf. \cite[Section~3.1]{MW} for periodic surfaces). Since \(*dz\) is closed on \(S\), \(s\) extends smoothly across \(\partial Q\).

Along each smooth arc of \(\partial Q\), the conormal of \(Q\) is parallel to \(e_x\) or \(e_y\). Hence \(s\) is constant on \(\partial Q=\gamma(\mathbb R)\), and we take \(s|_{\partial Q}=0\). Since \(dz|_Q\ne0\), \(w=z+is\) is locally conformal.

Since \(ds(N_y,-N_x,0)=|\nabla_S z|^2>0\), \(s\) is strictly increasing along each curve \(Q\cap\{z=t\}\). By \eqref{eq:hp-wing},
\[
\partial_r(s\circ X)=\frac{1+u_r^2}{\sqrt{1+|\nabla u|^2}}\ge C^{-1}.
\]
Thus \(s(Q\cap\{z=t\})=(0,\infty)\) for every \(t\in\mathbb R\), and \(w:\overline Q\to\overline{\mathbb H_+}\) is a homeomorphism.
\end{proof}

\begin{proposition}\label{prop:chain-gaps}
There are constants \(0<\delta\le D<\infty\) such that every \(S\in\mathscr H(M)\) satisfies
\[
S\cap\{x=y=0\}=\{(0,0,z_j):j\in\mathbb Z\}, \qquad \delta\le z_{j+1}-z_j\le D .
\]
The arcs \(\gamma((z_j,z_{j+1}))\) lie alternately in \(y=0\) and \(x=0\).
\end{proposition}
\begin{proof}
If \(S\) did not meet the axis, the connected curve \(\gamma(\mathbb R)\) would lie entirely in one reflection plane, say \(x=0\). Then \(S\cap\{y=0\}=\varnothing\) by reflection in \(x=0\), contrary to connectedness and symmetry in \(y=0\).

By compactness of \(\mathscr H(M)\) (Proposition~\ref{prop:fa-hull}), there is \(D>0\) such that every \(S\in\mathscr H(M)\) meets \(\{x=y=0, |z|<D/2\}\). By Lemma~\ref{lem:hp-saddles} and compactness, there is also \(\delta>0\) such that \(S\cap\{x=y=0, |z|<\delta\}=\{0\}\) whenever \(0\in S\). Since \(S-te_3\in\mathscr H(M)\) for every \(t\in\mathbb R\), \(\delta\le z_{j+1}-z_j\le D\). By \eqref{eq:hp-saddle}, the arcs \(\gamma((z_j,z_{j+1}))\) lie alternately in \(y=0\) and \(x=0\).
\end{proof}

\subsection{Weierstrass data}\label{sec:holomorphic-estimates}

For \(S\in\mathscr H(M)\), the inverse of \eqref{eq:hp-coordinate} has the form
\[
X_S(w)=(u(w),v(w),z), \qquad w=z+is.
\]
Since \(X_S\) is conformal,
\[
4(-u_w+iv_w)(u_w+iv_w)=1.
\]
By \eqref{eq:hp-signs}, \(N_x<0<N_y\) on \(Q\). Thus there is a holomorphic function \(f=f_S\) on \(\mathbb H_+\) satisfying
\[
e^f=2(-u_w+iv_w)=\frac{-N_x+iN_y}{1+N_z}, \qquad 0<\operatorname{Im}f<\pi/2.
\]
Then
\begin{equation}\label{eq:chain-weierstrass}
dX_S=\operatorname{Re}\left[(-\sinh f,-i\cosh f,1) dw\right],
\end{equation}
and
\begin{equation}\label{eq:chain-derivatives}
(u_s,v_s)=\cosh(\operatorname{Re}f) (\sin(\operatorname{Im}f),\cos(\operatorname{Im}f)).
\end{equation}

Set \(I_j=(z_j,z_{j+1})\) and \(\varepsilon_j=(-1)^j\). We may assume that \(\gamma(I_{2m})\subset\{y=0\}\). On \(I_j\), \(\operatorname{Im}f=\pi(1-\varepsilon_j)/4\). The boundary values \(u_\partial=a\) and \(v_\partial=b\) are continuous and uniformly bounded.

\begin{lemma}\label{lem:chain-coordinates}
For every \(S\in\mathscr H(M)\),
\begin{equation}\label{eq:chain-poisson-position}
u=s\cos\theta+P_s*u_\partial, \qquad v=s\sin\theta+P_s*v_\partial,
\end{equation}
where \(P_s(z)=s/[\pi(s^2+z^2)]\). If \(S_n\to S\) locally smoothly, then \(X_{S_n}\to X_S\) and \(f_{S_n}\to f_S\) locally smoothly on \(\mathbb H_+\). After reflection in \(x=0\) or \(y=0\), \(X_{S_n}\to X_S\) locally smoothly across \(\mathbb R\setminus\{z_j\}\). For \((0,0,z_{j,n})\in S_n\) with \(z_{j,n}\to z_j\), \(X_{S_n}(z_{j,n}+\zeta^2)\to X_S(z_j+\zeta^2)\) smoothly near \(\zeta=0\).
\end{lemma}
\begin{proof}
By \cite[Theorem~7.26]{ABR},
\[
u=c_1s+P_s*u_\partial, \qquad v=c_2s+P_s*v_\partial, \qquad c_1,c_2\ge0.
\]
Since \(u=c_1s+O(1)\) and \(v=c_2s+O(1)\), \eqref{eq:hp-cross-bounds} implies \(c_1\sin\theta=c_2\cos\theta\). Differentiating the Poisson integrals gives \(u_w=-ic_1/2+O(s^{-1})\) and \(v_w=-ic_2/2+O(s^{-1})\). Letting \(s\to\infty\) in \(4(u_w^2+v_w^2)=-1\), we obtain \(c_1^2+c_2^2=1\).

Under local smooth convergence \(S_n\to S\), the coordinates \(w_n\) in \eqref{eq:hp-coordinate} converge locally smoothly to \(w\), since \(dw_n=dz+i*dz\) and \(w_n=z\) on \(xy=0\). Since \(dw\ne0\) away from the axis, \(X_{S_n}\to X_S\) locally smoothly in \(\mathbb H_+\) and across \(\mathbb R\setminus\{z_j\}\). The identity \(e^{f_S}=2(-u_w+iv_w)\) also gives \(f_{S_n}\to f_S\) locally smoothly in \(\mathbb H_+\).

By Lemma~\ref{lem:hp-saddles}, \(w=z_j+\zeta^2\) in a holomorphic coordinate \(\zeta\) near \((0,0,z_j)\). In these coordinates, \(X_{S_n}(z_{j,n}+\zeta^2)\to X_S(z_j+\zeta^2)\) smoothly near \(\zeta=0\).
\end{proof}

Put
\[
\mathscr X=\{f_S:S\in\mathscr H(M)\} =\overline{\{f_M(\cdot+h):h\in\mathbb R\}}, \qquad (\tau_h f)(w)=f(w+h),
\]
with locally uniform convergence on \(\mathbb H_+\). Equation~\eqref{eq:chain-weierstrass} and \(\inf_{\mathbb H_+}u=\inf_{\mathbb H_+}v=0\) determine \(X_S\) from \(f_S\). Thus \(S\mapsto f_S\) is a homeomorphism from the compact family \(\mathscr H(M)\) onto \(\mathscr X\) (Proposition~\ref{prop:fa-hull} and Lemma~\ref{lem:chain-coordinates}). For \(h\in\mathbb R\), \(f_{S-he_3}=\tau_h f_S\).

For \(f\in\mathscr X\), put
\[
f_\infty=i(\pi/2-\theta), \qquad g=e^{2f}, \qquad g_\infty=e^{2f_\infty}, \qquad \mathcal B=\frac{g-g_\infty}{g-\bar g_\infty}.
\]
Then \(g:\mathbb H_+\to\mathbb H_+\) and \(\mathcal B:\mathbb H_+\to\mathbb D\).

For a singly periodic Scherk surface, we may assume that
\[
g(w)=-\frac{\sin(\pi w/P)}{\sin(\pi(w-a)/P)}
\]
for some \(P>0\) and \(a=2P\theta/\pi\in(0,P)\). Then \(\mathcal B(w)=e^{2\pi i(w-a)/P}\). Thus
\[
f^{-1}(f(w))=w+P\mathbb Z \qquad\text{for }w\in\mathbb H_+.
\]

\begin{lemma}\label{lem:chain-inner-data}
There is \(\rho>0\) such that for every \(f\in\mathscr X\) there are \(\lambda_j>0\) and holomorphic functions \(a_j:B_\rho(0)\to\mathbb C\setminus\{0\}\) satisfying
\begin{equation}\label{eq:chain-endpoint-expansion}
e^{f(z_j+\zeta^2)}=\zeta^{\varepsilon_j}a_j(\zeta^2) \qquad\text{for }0<|\zeta|^2<\rho.
\end{equation}
With \(r=|w-z_j|\), for \(0<r<\rho\) and \(0\le m\le3\),
\begin{equation}\label{eq:chain-vertex-norm}
|\nabla^m(u^2+v^2-\lambda_jr)|\le C_m r^{2-m}.
\end{equation}
On \(B_\rho(0)\), \(|(a_j^{\pm1})^{(k)}|\le C_k\) for \(k\ge0\). For \(s\ge1\) and \(k\ge1\),
\begin{equation}\label{eq:chain-inner-bounds}
|f-f_\infty|+|\mathcal B|\le C/s, \qquad |f^{(k)}|\le C_k/s^{k+1}.
\end{equation}
\(\mathcal B\) extends holomorphically across \(\mathbb R\), with \(|\mathcal B(z)|=1\). There are \(s_0,c,C>0\) such that, for \(0\le s\le s_0\),
\[
\mathcal B\ne0, \qquad c\le\partial_s(-\log|\mathcal B|), \qquad \left|\frac{\mathcal B'}{\mathcal B}\right|\le C.
\]
\end{lemma}
\begin{proof}
The map \(\zeta\mapsto X_S(z_j+\zeta^2)\) is a conformal immersion near \(0\), with \(u,v\) odd in \(\zeta\). The functions
\[
\zeta e^{\pm f(z_j+\zeta^2)}=\mp u_\zeta+iv_\zeta
\]
are even and holomorphic near \(0\). Since \((\zeta e^f)(\zeta e^{-f})=\zeta^2\) and \(\operatorname{Im}f=\pi(1-\varepsilon_j)/4\) on \(I_j\), \eqref{eq:chain-endpoint-expansion} follows. By conformality at \(\zeta=0\), \(u^2+v^2=\lambda_j|\zeta|^2+O(|\zeta|^4)\) with \(\lambda_j>0\), and \eqref{eq:chain-vertex-norm} follows by differentiation. Since \(S-z_je_3\in\mathscr H(M)\), the constants \(\rho,C_k,C_m\) can be chosen independently of \(f,j\) by Lemma~\ref{lem:chain-coordinates} and compactness.

By \eqref{eq:chain-poisson-position}, \(f-f_\infty=O(s^{-1})\), uniformly in \(z\). The definition of \(\mathcal B\) and Cauchy estimates imply \eqref{eq:chain-inner-bounds}.

Extend \(g\) by \(g(\bar w)=\overline{g(w)}\). By \eqref{eq:chain-endpoint-expansion}, \(\mathcal B\) extends holomorphically across \(\mathbb R\), with \(|\mathcal B(z)|=1\). By the Hopf lemma,
\[
\partial_s(-\log|\mathcal B|)(z,0) =-i\frac{\mathcal B'(z)}{\mathcal B(z)}>0.
\]
By Lemma~\ref{lem:chain-coordinates} and the compactness and translation invariance of \(\mathscr X\), there are \(s_0,c,C>0\) independent of \(f\) such that, for \(0\le s\le s_0\),
\[
\mathcal B\ne0, \qquad \partial_s(-\log|\mathcal B|)\ge c, \qquad \left|\frac{\mathcal B'}{\mathcal B}\right|\le C.\qedhere
\]
\end{proof}

By \eqref{eq:chain-endpoint-expansion}, \(f''/f'=-1/(w-z_j)+O(1)\) near \(z_j\). Define
\begin{equation}\label{eq:complete-normalized-Q}
\mathcal Q(w)=\sum_j \left(\frac1{w-z_j}+\frac{z_j}{1+z_j^2}\right), \qquad \mathcal V=\mathcal V_f=\frac{f''}{f'}+\mathcal Q .
\end{equation}
The summands are \(O(|z_j|^{-2})\) locally uniformly in \(w\) as \(|z_j|\to\infty\). Since \(z_{j+1}-z_j\ge\delta\), the series converges locally uniformly on \(\mathbb C\setminus\{z_j\}\) and depends continuously on \(f\in\mathscr X\). \(\mathcal V\) extends holomorphically across each \(z_j\). For \(s\ge1\) and \(k\ge1\),
\begin{equation}\label{eq:complete-Q-bounds}
|\operatorname{Re}\mathcal Q(w)|\le C\log(2+|z|+s), \qquad |\operatorname{Im}\mathcal Q(w)|\le C,\quad |\mathcal Q^{(k)}(w)|\le C_k s^{-k}.
\end{equation}
Here \(\operatorname{Im}\mathcal Q=-s\sum_j|w-z_j|^{-2}\) and \(\sum_j|w-z_j|^{-k-1}\le C_k s^{-k}\) for \(k\ge1\). The bound on \(\operatorname{Re}\mathcal Q\) follows by integrating \(\mathcal Q'\) along \(i\to i(s+|z|)\to w\), using \(\operatorname{Re}\mathcal Q(i)=0\).

\subsection{Inverse images and integral estimates}\label{sec:height-integral-identity}

Put \(S_0=\{\alpha+i\phi:0<\phi<\pi/2\}\). For \(\zeta\in S_0\setminus\{f_\infty\}\), put
\[
\omega=e^{2\zeta}, \qquad \lambda=\frac{\omega-g_\infty}{\omega-\bar g_\infty}.
\]
Let \(Z_j=x_j+iy_j\) be the zeros of \(f-\zeta\), counted with multiplicity. Then \(\mathcal B^{-1}(\lambda)=\{Z_j\}\). By Lemma~\ref{lem:equal-height-fibers}, \(S\) is a singly periodic Scherk surface if \(y_j=y_k\) for all \(j,k\) and almost every \(\zeta\).

\begin{lemma}\label{lem:infinite-fiber}
For \(f\in\mathscr X\) and \(\zeta\in S_0\setminus\{f_\infty\}\), \(f^{-1}(\zeta)\) is infinite. For every compact \(K\Subset S_0\setminus\{f_\infty\}\),
\[
0<c_K\le y_j\le C_K, \qquad \#\{j:x_j\in[z,z+1]\}\le C_K,
\]
for \(\zeta\in K\) and \(z\in\mathbb R\), with constants independent of \(f\). For every regular value \(\zeta\),
\begin{equation}\label{eq:infinite-height-identity}
\operatorname{Im}\mathcal V(Z_j) =-\sum_{k\ne j}a_{jk}(y_j-y_k), \qquad a_{jk}=\frac{2}{|Z_j-Z_k|^2}.
\end{equation}
For each \(j\), \(\sum_{k\ne j}a_{jk}<\infty\).
\end{lemma}
\begin{proof}
By Lemma~\ref{lem:chain-inner-data}, \(c_K\le y_j\le C_K\). By the Schwarz lemma applied to \(\Phi_\omega=(g-\omega)/(g-\bar\omega)\),
\[
\sum_j\log\left|\frac{w-\bar Z_j}{w-Z_j}\right| \le-\log|\Phi_\omega(w)|.
\]
For fixed \(Y>2C_K\),
\[
\sum_j\frac1{1+x_j^2}\le-C\log|\Phi_\omega(iY)|<\infty.
\]
Hence \(\sum_j\log|(w-\bar Z_j)/(w-Z_j)|\) converges locally uniformly on \(\overline{\mathbb H_+}\setminus\{Z_j\}\). Thus
\[
-\log|\Phi_\omega(w)|-\sum_j\log\left|\frac{w-\bar Z_j}{w-Z_j}\right|
\]
is nonnegative and harmonic on \(\mathbb H_+\), with zero boundary values, hence equals \(cs\), \(c\ge0\), by \cite[Theorem~7.22]{ABR}. Since \(\Phi_\omega(iY)\to(g_\infty-\omega)/(g_\infty-\bar\omega)\ne0\), we have \(c=0\). Differentiating in \(s\) at \(s=0\),
\begin{equation}\label{eq:infinite-poisson-count}
\begin{aligned}
q_\lambda(z)&:= -i\frac{\mathcal B'(z)}{\mathcal B(z)} \frac{1-|\lambda|^2}{|\mathcal B(z)-\lambda|^2} \\
&=-\left.\partial_s\log|\Phi_\omega(z+is)|\right|_{s=0} =\sum_j\frac{2y_j}{(z-x_j)^2+y_j^2}.
\end{aligned}
\end{equation}
By Lemma~\ref{lem:chain-inner-data}, \(c_K\le q_\lambda\le C_K\). If \(f^{-1}(\zeta)\) were finite, \eqref{eq:infinite-poisson-count} would imply \(q_\lambda(z)\to0\) as \(|z|\to\infty\), contrary to \(q_\lambda\ge c_K>0\). Moreover, by \eqref{eq:infinite-poisson-count}, \(\#\{j:x_j\in[z,z+1]\}\le C_Kq_\lambda(z+1/2)\le C_K\).

Thus \(\sum_{k\ne j}a_{jk}<\infty\) for every regular value \(\zeta\). Set
\[
H_\omega=\frac{\mathcal Q}{2}-f'+\frac{g'}{g-\omega}.
\]
Since \(g=e^{2f}\),
\[
H_\omega(w)=\frac1{w-Z_j}+\frac{\mathcal V(Z_j)}2+O(w-Z_j), \qquad w\to Z_j.
\]
By \eqref{eq:chain-endpoint-expansion} and compactness of \(\mathscr X\), \(H_\omega\) is holomorphic across \(\mathbb R\), with \mbox{\(|\operatorname{Im}H_\omega|\le C\)} for \(-1\le\operatorname{Im}w\le0\). Since \(g(\bar w)=\overline{g(w)}\), \eqref{eq:chain-inner-bounds} and \eqref{eq:complete-Q-bounds} imply \mbox{\(|\operatorname{Im}H_\omega|\le C\)} for \(\operatorname{Im}w\le-1\). On \(\mathbb R\), \(\operatorname{Im}H_\omega=q_\lambda/2\). For \(\operatorname{Im}w<0\), the Poisson formula gives
\[
\operatorname{Im}H_\omega(w) =\sum_j\operatorname{Im}\frac1{w-Z_j}.
\]
The identity holds throughout \(\mathbb C\setminus\{Z_j\}\), since both sides are harmonic there. Comparing the constant terms at \(Z_j\) gives \eqref{eq:infinite-height-identity}.
\end{proof}

For a regular value \(\zeta\) and bounded \(h=(h_j)\) with \(h_j\in\mathbb R\), define the weighted discrete Laplacian \(L_\zeta\) by
\[
(L_\zeta h)_j=2\sum_{k\ne j}\frac{h_j-h_k}{|Z_j-Z_k|^2}.
\]
For finitely supported \(h\),
\[
\sum_j h_j(L_\zeta h)_j
=\frac12\sum_{j\ne k}a_{jk}(h_j-h_k)^2\ge0.
\]
For \(y=(y_j)\), \eqref{eq:infinite-height-identity} reads \mbox{\(-\operatorname{Im}\mathcal V(Z_j)=(L_\zeta y)_j\)}.
Recall that
\[
E(f)=\int_{S_0}\cosh(2\alpha) \sum_{j\ne k}a_{jk}(y_j-y_k)^2\mathrm d A,
\]
and, for \(\psi\in C_c^\infty(\mathbb R;\mathbb R)\),
\begin{equation}\label{eq:infinite-local-work}
W_\psi=\int_{S_0}\cosh(2\alpha)\sum_j\psi(x_j)y_j(L_\zeta y)_j\,\mathrm dA.
\end{equation}
Put \(X_h(w)=(u(w),v(w),0)\), \(\mathcal P=|X_h|^2-s^2\), and \(J=\mathcal V\mathcal P_{ww}\).

\begin{lemma}\label{lem:inverse-image-integral-identity}
For \(f\in\mathscr X\),
\begin{equation}\label{eq:infinite-current}
W_\psi=-\operatorname{Re}\int_{\mathbb H_+}\psi'J\mathrm d z\mathrm d s -\frac12\int_0^\infty s\left(\int_{\mathbb R} \psi''\operatorname{Im}J\mathrm d z\right)\mathrm d s.
\end{equation}
\end{lemma}
\begin{proof}
Since \(z,s\) are harmonic conjugates, \(|\nabla_Ss|^2=|\nabla_Sz|^2=1-N_z^2\). Thus
\[
\Delta_S\mathcal P=2(1+N_z^2)-2(1-N_z^2)=4N_z^2.
\]
Also, \(N_z=-\tanh(\operatorname{Re}f)\) and \(\mathrm d\mu_S=\cosh^2(\operatorname{Re}f)\,\mathrm dz\,\mathrm ds\). Since \(\operatorname{Re}f\) is harmonic,
\[
\Delta_w\mathcal P=4\sinh^2(\operatorname{Re}f),\qquad
\Delta_w^2\mathcal P=8\cosh(2\operatorname{Re}f)|f'|^2,
\]
where \(\Delta_w=\partial_z^2+\partial_s^2\). By conformality, \(2(X_h)_w\cdot(X_h)_w=(s^2)_{ww}=-1/2\), so
\[
\mathcal P_{ww}=2X_h\cdot(X_h)_{ww}=-(u\cosh f+iv\sinh f)f'.
\]
Consequently,
\begin{equation}\label{eq:infinite-current-derivative}
J =-(u\cosh f+iv\sinh f)(f''+\mathcal Qf'),\quad J_{\bar w\bar w} =\tfrac12\cosh(2\operatorname{Re}f)|f'|^2\mathcal V.
\end{equation}

By \eqref{eq:chain-endpoint-expansion} and \(\sum_{k\ne j}|z_k-z_j|^{-2}\le C\), we have \(|\mathcal V'|\le C\) on \(|w-z_j|<\rho\), where \(\rho,C>0\) are independent of \(f\) and \(j\). For some \(s_1>0\), Lemma~\ref{lem:chain-inner-data} implies \(|f'|\ge c\) when \(0\le s\le s_1\) and \(|w-z_j|\ge\rho\) for all \(j\). The functions \(f'',f''',\mathcal Q'\) are bounded there. Hence \(\mathcal V\) is holomorphic for \(|s|<s_1\), with \(\operatorname{Im}\mathcal V=0\) on \(\mathbb R\), and satisfies for \(0\le s\le s_1\)
\begin{equation}\label{eq:infinite-collar-velocity}
|\mathcal V'|\le C, \qquad |\operatorname{Im}\mathcal V(z+is)|\le Cs.
\end{equation}
By \eqref{eq:chain-endpoint-expansion} and \eqref{eq:infinite-collar-velocity}, \(s|\operatorname{Im}J_{\bar w\bar w}|\le C|w-z_j|^{-1}\) for \(|w-z_j|<\rho\). For \(s\ge1\), \eqref{eq:chain-inner-bounds} and \eqref{eq:complete-Q-bounds} give \(s|\operatorname{Im}J_{\bar w\bar w}|\le Cs^{-3}\). By \eqref{eq:infinite-height-identity} and the area formula,
\[
W_\psi=-\int_{\mathbb H_+}s\psi(z)\cosh(2\operatorname{Re}f)|f'|^2\operatorname{Im}\mathcal V\,\mathrm dz\,\mathrm ds =-2\int_{\mathbb H_+}s\psi(z)\operatorname{Im}J_{\bar w\bar w}\,\mathrm dz\,\mathrm ds.
\]

By \eqref{eq:chain-endpoint-expansion}--\eqref{eq:chain-vertex-norm}, \(|u\cosh f+iv\sinh f|\le C\) for \(|w-z_j|<\rho\). By Lemma~\ref{lem:chain-coordinates} and compactness of \(\mathscr X\), this bound also holds for \(0<s\le1\) with \(|w-z_j|\ge\rho\) for every \(j\). For \(s\ge1\), \eqref{eq:chain-poisson-position} and \eqref{eq:chain-inner-bounds} give
\[
u\cosh f+iv\sinh f =is\sinh(f-f_\infty)+O(1)=O(1).
\]
Consequently,
\begin{equation}\label{eq:infinite-two-phase-bound}
|\mathcal P_{ww}|=|(u\cosh f+iv\sinh f)f'|\le C|f'|.
\end{equation}
By Lemma~\ref{lem:chain-inner-data} and \eqref{eq:complete-normalized-Q}--\eqref{eq:complete-Q-bounds},
\begin{align}
|J|+s|J_z|&\le Cs^{-2}\log(2+|z|+s) &&\text{for }s\ge1,\notag \\
|w-z_j| |J|&\le C\log(2+|z_j|) &&\text{for }0<|w-z_j|<\rho. \label{eq:infinite-product-bounds}
\end{align}

On \(\Omega=\{0<s<T, |w-z_j|>\varepsilon\text{ for all }j\}\), integration by parts gives
\[
\begin{aligned}
&-2\operatorname{Im}\int_\Omega s\psi J_{\bar w\bar w} \\
={}& -2\operatorname{Im}\int_\Omega \left(\tfrac14s\psi''+\tfrac i2\psi'\right)J +\operatorname{Re}\int_{\partial\Omega} \left[s\psi J_{\bar w} -\tfrac12(s\psi'+i\psi)J\right]dw.
\end{aligned}
\]
Since \(J\in\mathbb R\) on \(\mathbb R\), the integral over \(\partial\Omega\cap\mathbb R\) has real part \(0\). By \eqref{eq:chain-vertex-norm}, the integral over \(\partial\Omega\cap\partial B_\varepsilon(z_j)\) has real part
\[
-\frac{\lambda_j\mathcal V(z_j)\psi(z_j)}8\int_\pi^0\cos\phi\,\mathrm d\phi+O_\psi(\varepsilon)=O_\psi(\varepsilon),
\]
where \(w=z_j+\varepsilon e^{i\phi}\).
The integral over \(\partial\Omega\cap\{s=T\}\) is \(O_\psi(\log(2+T)/T)\).

For \(s\ge1\), the identity \(\int\psi''J\mathrm dz=-\int\psi'J_z\mathrm dz\) and \eqref{eq:infinite-product-bounds} give
\[
\left|\int_{\mathbb R}\left(\frac{s}{4}\psi''+\frac{i}{2}\psi'\right)J\mathrm dz\right|\le C_\psi\frac{\log(2+s)}{s^2}.
\]
Letting \(\varepsilon\downarrow0\) and \(T\to\infty\) proves \eqref{eq:infinite-current}.
\end{proof}

Choose a nonincreasing \(\eta\in C^\infty(\mathbb R)\) with \(\eta=1\) on \(({-\infty},1/2]\) and \(\eta=0\) on \([1,\infty)\). For \(R\ge2\), set
\begin{equation}\label{eq:complete-cutoff}
\chi_R(z)=\eta \left(\frac{\log|z|}{\log R}\right), \qquad \chi_R(0)=1, \qquad \psi_R=\chi_R^2.
\end{equation}
\(\chi_R\) will be used in the proof of Lemma~\ref{lem:complete-target-error}.

\begin{proposition}\label{prop:logarithmic-height-estimate}
For every \(f\in\mathscr X\),
\begin{equation}\label{eq:complete-local-average}
|W_{\psi_R}|\le C.
\end{equation}
If \(\mathcal V_{f_*}'=0\) for every limit \(f_*=\lim_{n\to\infty}\tau_{h_n}f=\lim_{n\to\infty}f(\cdot+h_n)\) with \(|h_n|\to\infty\), then \(W_{\psi_R}\to0\) as \(R\to\infty\).
\end{proposition}
\begin{proof}
Choose \(0<\varepsilon<1\) so that \eqref{eq:infinite-collar-velocity} holds on \(0\le s\le\varepsilon\). By compactness and translation invariance of \(\mathscr X\), \(c\le|f'(z+i\varepsilon)|\le C\) and \(|f''/f'(z+i\varepsilon)|\le C\). Set
\[
c_f(z)=\operatorname{Re}\mathcal V(z+i\varepsilon), \qquad r_f(z)=\operatorname{Re}\int_0^\infty (\mathcal V(z+is)-c_f(z))\mathcal P_{ww}(z+is) ds.
\]
The difference \(\mathcal V_{\tau_hf}-\mathcal V_f(\cdot+h)\) is a real constant. Hence
\[
\mathcal V_{\tau_hf}'=\mathcal V_f'(\cdot+h), \qquad r_{\tau_hf}=r_f(\cdot+h).
\]
For \(s\ge1\),
\[
|\mathcal Q(z+is)-c_f(z)| \le C+\int_\varepsilon^s|\mathcal Q'(z+i\rho)| \mathrm d\rho \le C\log(2+s).
\]
The integrand defining \(r_f\) is therefore uniformly \(O(\log(2+s)/s^2)\) as \(s\to\infty\) and \(O(|w-z_j|^{-1})\) near \(z_j\). Consequently,
\begin{equation}\label{eq:complete-remainder-bound}
\sup_{t\in\mathbb R}\int_t^{t+1}|r_f(z)|\,\mathrm dz\le C.
\end{equation}
By \eqref{eq:infinite-product-bounds} and \(|\psi_R^{(k)}(z)|\le C/(|z|^k\log R)\) for \(k=1,2\),
\[
\begin{aligned}
&\left|\int_0^\infty s\left(\int_{\mathbb R}\psi_R''\operatorname{Im}J\,\mathrm dz\right)\mathrm ds\right|\\
\le{}&\int_0^{\sqrt R}\int_{\mathbb R}s|\psi_R''J|\,\mathrm dz\mathrm ds
+\int_{\sqrt R}^{\infty}\int_{\mathbb R}s|\psi_R'J_z|\,\mathrm dz\mathrm ds
\le\frac{C\log R}{\sqrt R}.
\end{aligned}
\]
With \(A(z)=\operatorname{Re}\int_0^\infty\mathcal P_{ww}(z+is)\,\mathrm ds\), \eqref{eq:infinite-current} becomes
\[
W_{\psi_R}=-\int_{\mathbb R}\psi_R'(c_fA+r_f)\,\mathrm dz+O(\log R/\sqrt R).
\]

For \(e=(\cos\theta,\sin\theta,0)\), \(2e\cdot(X_h)_w=-i\cosh(f-f_\infty)=-i+O(s^{-2})\). Thus \(e\cdot X_h=s+\kappa+O(s^{-1})\) for a constant \(\kappa\). By \eqref{eq:chain-poisson-position},
\[
\mathcal P_s=2\kappa+O(s^{-1}),\qquad \mathcal P_z=O(s^{-1})\qquad\text{for }s\ge1.
\]
Also, \(\mathcal P_s(z,0)=0\) away from the \(z_j\), and \(\mathcal P_z\) is bounded by \eqref{eq:chain-vertex-norm}. For \(T\ge1\), put
\[
b_T(z)=\frac14\int_0^T \mathcal P_z(z,s)\,\mathrm ds,\qquad \|b_T\|_{L^\infty(\mathbb R)}\le C\log(2+T).
\]
For \(\psi\in C_c^\infty(\mathbb R)\), the identity \(4\operatorname{Re}\mathcal P_{ww}=\mathcal P_{zz}-\mathcal P_{ss}\) and the bound \(|\mathcal P_{ww}|\le Cs^{-2}\) for \(s\ge1\) give
\begin{equation}\label{eq:complete-general-mean}
\left|\int_{\mathbb R}\left[(A+\kappa/2)\psi+b_T\psi'\right] dz\right| \le\frac CT\int_{\mathbb R}|\psi| dz.
\end{equation}

Write \(c_f=a_f+d_f\), where
\[
a_f(z)=\operatorname{Re}\mathcal Q(i(|z|+\varepsilon)).
\]
Then \(a_f\) is even and Lipschitz, with \(|a_f(z)|\le C\log(2+|z|)\) and \(|a_f'(z)|\le C/(1+|z|)\). For \(\ell\ge2\), \(F(w)=\mathcal Q(w+i\varepsilon)-a_f(\ell)\) satisfies
\[
|\operatorname{Im}F|\le C, \qquad \operatorname{Re}F(i\ell)=0, \qquad F^2=O(\log^2(2+|w|)).
\]
By the Poisson formula for \(\operatorname{Re}(F^2)\) at \(i\ell\),
\[
\int_{\mathbb R}\frac{\ell\,(\operatorname{Re}F(z))^2}{\pi(\ell^2+z^2)}\,\mathrm dz\le C.
\]
Since \(|c_f(z)-\operatorname{Re}\mathcal Q(z+i\varepsilon)|+|a_f(z)-a_f(\ell)|\le C\) for \(\ell\le|z|\le2\ell\),
\[
\int_{\ell\le|z|\le2\ell}|d_f(z)|^2 dz\le C\ell.
\]
\(\psi_R'\) is supported in \(\{\sqrt R\le|z|\le R\}\). By the Cauchy--Schwarz inequality on \(\{2^k\le|z|\le2^{k+1}\}\) and the derivative bounds for \(\psi_R\),
\begin{equation}\label{eq:complete-cutoff-mean-bounds}
\begin{aligned}
&\|\psi_R'd_f\|_{L^1(\mathbb R)}+\sqrt R \|\psi_R''d_f\|_{L^1(\mathbb R)}\\
\le{}&\frac C{\log R}\sum_{\sqrt R/2\le2^k\le R}\left(1+\frac{\sqrt R}{2^k}\right)\le C.
\end{aligned}
\end{equation}

By \eqref{eq:complete-general-mean} and \(\int_{\mathbb R}\psi_R'a_f\,\mathrm dz=0\),
\[
\begin{aligned}
&\int_{\mathbb R}\psi_R'a_f A\,\mathrm dz \\
={}& -\int_{\mathbb R}(\psi_R''a_f+\psi_R'a_f')b_{\sqrt R}\,\mathrm dz +O\left(R^{-1/2}\int_{\mathbb R}|\psi_R'a_f|\,\mathrm dz\right) =O(\log R/\sqrt R).
\end{aligned}
\]
Take \(\psi=\psi_R'd_f\) in \eqref{eq:complete-general-mean}. Since \(\int_{\mathbb R}\psi_R'd_f\,\mathrm dz=\int_{\mathbb R}\psi_R'c_f\,\mathrm dz\), \eqref{eq:complete-cutoff-mean-bounds} gives
\begin{equation}\label{eq:complete-cutoff-work}
\begin{aligned}
W_{\psi_R}={}&\int_{\mathbb R}(\psi_R''d_f+\psi_R'd_f')b_T\,\mathrm dz-\int_{\mathbb R}\psi_R'r_f\,\mathrm dz\\
&+\frac\kappa2\int_{\mathbb R}\psi_R'c_f\,\mathrm dz+O(T^{-1})+o_R(1).
\end{aligned}
\end{equation}
Put \(\Psi_R(z,s)=\psi_R(\sqrt{z^2+s^2})\). Since
\[
c_f(z)=\operatorname{Re}\mathcal Q(z+i\varepsilon)+\partial_z\log|f'(z+i\varepsilon)|,
\]
Green's identity on \(\mathbb H_+\), applied to \(\Psi_R\) and \(\operatorname{Im}\mathcal Q(\cdot+i\varepsilon)\), gives
\[
\int_{\mathbb R}\psi_R'c_f\,\mathrm dz =-\int_{\mathbb H_+}\operatorname{Im}\mathcal Q(z+i(s+\varepsilon))\Delta\Psi_R\mathrm d z\mathrm d s -\int_{\mathbb R}\psi_R''\log|f'(z+i\varepsilon)|\mathrm d z.
\]
Since \(\|\Delta\Psi_R\|_{L^1(\mathbb H_+)}=\pi\int_0^\infty|(z\psi_R'(z))'|\mathrm d z\), the bounds on \(\operatorname{Im}\mathcal Q\) and \(\log|f'(\cdot+i\varepsilon)|\) imply
\[
\left|\int_{\mathbb R}\psi_R'c_f\,\mathrm dz\right| \le C\int_0^\infty|(z\psi_R'(z))'| dz+C\|\psi_R''\|_{L^1(\mathbb R)} \le C/\log R.
\]
Taking \(T=1\) in \eqref{eq:complete-cutoff-work} and using \eqref{eq:complete-remainder-bound}, \eqref{eq:complete-cutoff-mean-bounds}, and \(|d_f'|\le C\) proves \eqref{eq:complete-local-average}.

Suppose that \(\mathcal V_{f_*}'=0\) whenever \(f(\cdot+h_n)\to f_*\) and \(|h_n|\to\infty\). Then \(r_{f_*}=0\).

The map \(f\mapsto r_f\) is continuous from \(\mathscr X\) to \(L^1([0,1])\). For \(f_n=f_{S_n}\to f_*\) in \(\mathscr X\), let \((0,0,z_{j,n})\) be the intersections of \(S_n\) with the axis, indexed so that \(z_{j,n}\to z_{j,*}\). Fix \(0<r<\rho\) and \(T\ge1\). Lemma~\ref{lem:chain-coordinates} and \eqref{eq:infinite-current-derivative} imply uniform convergence of the integrands defining \(r_{f_n}\) on
\[
\{z+is:0\le z\le1,\ 0\le s\le T\}\setminus\bigcup_jB_r(z_{j,*}).
\]
The integrands are uniformly \(O(|w-z_{j,n}|^{-1})\) near \(z_{j,n}\) and \(O(\log(2+s)/s^2)\) as \(s\to\infty\). Hence
\[
\limsup_{n\to\infty}\|r_{f_n}-r_{f_*}\|_{L^1([0,1])}\le Cr+\frac{C\log(2+T)}{T},
\]
with \(C\) independent of \(r,T\). Letting \(r\downarrow0\) and \(T\to\infty\) proves \(r_{f_n}\to r_{f_*}\) in \(L^1([0,1])\).

By compactness of \(\mathscr X\), the continuity of \(f\mapsto r_f\), and \(a_f'(z)\to0\),
\[
d_f'(z)=\operatorname{Re}\mathcal V_f'(z+i\varepsilon)-a_f'(z)\longrightarrow0, \qquad \int_z^{z+1}|r_f(u)| du\longrightarrow0
\]
as \(|z|\to\infty\).
Thus \eqref{eq:complete-cutoff-mean-bounds} gives
\[
\|\psi_R''d_f\|_{L^1(\mathbb R)}+\|\psi_R'd_f'\|_{L^1(\mathbb R)}+\|\psi_R'r_f\|_{L^1(\mathbb R)}\longrightarrow0
\]
as \(R\to\infty\).
Since \(b_T\) is bounded for fixed \(T\ge1\) and \(\int_{\mathbb R}\psi_R'c_f\,\mathrm dz=O(1/\log R)\), \eqref{eq:complete-cutoff-work} gives
\[
\limsup_{R\to\infty}|W_{\psi_R}|\le\frac CT.
\]
Letting \(T\to\infty\), we obtain \(W_{\psi_R}\to0\).
\end{proof}

\subsection{Periodicity}\label{sec:axis-periodicity}

For a regular value \(\zeta=\alpha+i\phi\in S_0\setminus\{f_\infty\}\), recall that \(Z_j=x_j+iy_j\) are the points of \(f^{-1}(\zeta)\) and \(a_{jk}=2/|Z_j-Z_k|^2\).

\begin{lemma}\label{lem:complete-target-error}
For every \(f\in\mathscr X\),
\begin{equation}\label{eq:complete-height-energy}
E(f)\le2\liminf_{R\to\infty}W_{\psi_R}\le C.
\end{equation}
\end{lemma}
\begin{proof}
Fix a regular value \(\zeta\ne f_\infty\). Since \(\psi_R=\chi_R^2\) and \(a_{jk}=a_{kj}\),
\begin{gather}
\frac12\sum_{j\ne k}a_{jk} \left(\chi_R(x_j)y_j-\chi_R(x_k)y_k\right)^2 =\sum_{j\ne k}a_{jk}\psi_R(x_j)y_j(y_j-y_k)+\mathcal E_R(\zeta),\notag \\
\mathcal E_R(\zeta)=\frac12\sum_{j\ne k} a_{jk}y_jy_k\left(\chi_R(x_j)-\chi_R(x_k)\right)^2. \label{eq:complete-ims}
\end{gather}
Integrating \eqref{eq:complete-ims} against \(\cosh(2\alpha)\mathrm d A\) and applying Fatou's lemma, we obtain \(E(f)\le2\liminf_{R\to\infty}W_{\psi_R}\) provided
\[
\int_{S_0}\cosh(2\alpha)\mathcal E_R\,\mathrm d A\longrightarrow0
\]
as \(R\to\infty\).

For \(w,\xi\in\mathbb H_+\), write \(w=z+is\) and \(\xi=x+iy\). By symmetry in \(j,k\) and the area formula,
\begin{align}
&\int_{S_0}\cosh(2\alpha)\mathcal E_R\mathrm d A \label{eq:complete-source-estimate}\\
\le{}& 2\int_{\mathbb H_+}\cosh(2\operatorname{Re}f)|f'|^2 \sum_{\substack{f(\xi)=f(w), \xi\ne w\\0<y\le s}} \frac{sy}{|w-\xi|^2} |\chi_R(z)-\chi_R(x)|^2\mathrm d z\mathrm d s .\notag
\end{align}
By \eqref{eq:infinite-poisson-count},
\[
q(z):=\sum_j\frac{2y_j}{(z-x_j)^2+y_j^2}
=\frac{2\operatorname{Im}(e^{2\zeta})g'(z)}{|g(z)-e^{2\zeta}|^2}.
\]
Since \(g\) increases from \(-\infty\) to \(+\infty\) on each \((z_{2j+1},z_{2j+3})\),
\[
\int_{z_{2j+1}}^{z_{2j+3}}q(z)\,\mathrm dz=2\pi.
\]
Since \(z_{j+1}-z_j\ge\delta\), \(\int_Iq\le C(|I|+1)\) for every interval \(I\subset\mathbb R\), uniformly in \(\zeta\ne f_\infty\). If \(|I|=r\ge1\),
\begin{equation}\label{eq:complete-scale-count}
\#\{j:x_j\in I, 0<y_j\le r\} \le C\int_{I+[-r,r]}q(z) dz\le Cr.
\end{equation}

Write \([\chi]_{\dot H^{1/2}(\mathbb R)}^2 =\iint_{\mathbb R^2}\frac{|\chi(z)-\chi(u)|^2}{|z-u|^2}\mathrm d z\mathrm d u\).
For \(r>0\), set \(I_n=[nr,(n+1)r)\) and \(c_n=r^{-1}\int_{I_n}\chi_R\). For \(z\in I_n\) and \(x\in I_m\),
\[
|\chi_R(z)-\chi_R(x)|^2 \le3\left(|c_n-c_m|^2+\operatorname{osc}_{I_n}(\chi_R)^2+\operatorname{osc}_{I_m}(\chi_R)^2\right).
\]
On \(I_n\times I_m\), \(|z-u|^2\le2r^2(1+(n-m)^2)\), and \(\sum_n\operatorname{osc}_{I_n}(\chi_R)\le2\). Hence
\begin{equation}\label{eq:complete-interval-comparison}
\sum_{n,m}\frac{|c_n-c_m|^2}{1+(n-m)^2} \le2[\chi_R]_{\dot H^{1/2}(\mathbb R)}^2,\quad \sum_n\operatorname{osc}_{I_n}(\chi_R)^2 \le2\sqrt r \|\chi_R'\|_{L^2(\mathbb R)} .
\end{equation}

By Lemma~\ref{lem:chain-inner-data}, there is \(0<\sigma\le1\) such that \(\mathcal B\) is injective on \(B_{2\sigma}(z)\) for every \(z\in\mathbb R\). If \(w\ne \xi\), \(f(w)=f(\xi)\), and \(0<y\le s\le\sigma\), then \(|z-x|\ge\sigma\) and
\[
|w-\xi|^{-2}\le\frac{C}{1+(z-x)^2}.
\]
By \eqref{eq:chain-endpoint-expansion} and compactness of \(\mathscr X\),
\[
\sup_n\int_n^{n+1}\int_0^\sigma s^2\cosh(2\operatorname{Re}f)|f'|^2\mathrm d s\mathrm d z\le C.
\]
Using \eqref{eq:complete-scale-count} and \eqref{eq:complete-interval-comparison} with \(r=1\),
\[
\begin{aligned}
&\int_0^\sigma\int_{\mathbb R}\cosh(2\operatorname{Re}f)|f'|^2
\sum_{\substack{f(\xi)=f(w),\ \xi\ne w\\0<y\le s}}
\frac{sy}{|w-\xi|^2}|\chi_R(z)-\chi_R(x)|^2\mathrm d z\mathrm d s \\
\le{}& C\left([\chi_R]_{\dot H^{1/2}(\mathbb R)}^2+\|\chi_R'\|_{L^2(\mathbb R)}\right).
\end{aligned}
\]

For \(s\ge\sigma\), put \(r=1+s\) and take \(\xi\ne w\) with \(f(\xi)=f(w)\) and \(0<y\le s\). By compactness of \(\mathscr X\) and \eqref{eq:chain-inner-bounds}, \(\cosh(2\operatorname{Re}f(w))\le C\), \(|f'(w)|\le Cr^{-2}\), and \(|f''|\le Cr^{-3}\) on \(B_{s/2}(w)\). If \(|w-\xi|<s/2\), then \(0=f(\xi)-f(w)=f'(w)(\xi-w)+O(r^{-3}|w-\xi|^2)\). If \(|w-\xi|\ge s/2\), then \(r\le C|w-\xi|\). Thus \(|f'(w)|\le C\min\{r^{-2},r^{-3}|w-\xi|\}\), and
\[
|f'(w)|^2\frac{sy}{|w-\xi|^2} \le\frac{C}{r^2(r^2+(z-x)^2)}.
\]
For \(z\in I_n\), \eqref{eq:complete-scale-count} gives
\[
|f'(w)|^2 \sum_{\substack{f(\xi)=f(w), \xi\ne w\\x\in I_m, 0<y\le s}} \frac{sy}{|w-\xi|^2} \le\frac{C}{r^3(1+(n-m)^2)}.
\]
By \eqref{eq:complete-interval-comparison},
\[
\begin{aligned}
&\int_{\mathbb R}\cosh(2\operatorname{Re}f)|f'|^2
\sum_{\substack{f(\xi)=f(w),\ \xi\ne w\\0<y\le s}}
\frac{sy}{|w-\xi|^2}|\chi_R(z)-\chi_R(x)|^2\mathrm d z \\
\le{}& C\left(r^{-2}[\chi_R]_{\dot H^{1/2}(\mathbb R)}^2+r^{-3/2}\|\chi_R'\|_{L^2(\mathbb R)}\right).
\end{aligned}
\]
Integrating over \(s\ge\sigma\) and combining with the bound for \(0<s\le\sigma\),
\begin{equation}\label{eq:complete-target-error}
\int_{S_0}\cosh(2\alpha)\mathcal E_R\mathrm d A \le C\left([\chi_R]_{\dot H^{1/2}(\mathbb R)}^2+\|\chi_R'\|_{L^2(\mathbb R)}\right).
\end{equation}

Since \(|\chi_R'|\le C/(|z|\log R)\) on \(\sqrt R<|z|<R\) and \(\chi_R'=0\) elsewhere,
\[
\begin{gathered}
\int_{\mathbb H_+}\left|\nabla\bigl[\chi_R(\sqrt{z^2+s^2})\bigr]\right|^2\mathrm d z\mathrm d s
\le\frac{C}{(\log R)^2}\int_{\sqrt R}^R\frac{\mathrm d r}{r}\le\frac C{\log R},\\
\|\chi_R'\|_{L^2(\mathbb R)}\le\frac{C}{R^{1/4}\log R}.
\end{gathered}
\]
The trace inequality \cite[Chapter~10, Proposition~17.1]{DiBenedetto} applied to \(\chi_R(\sqrt{z^2+s^2})\) gives
\begin{equation}\label{eq:complete-cutoff-bounds}
[\chi_R]_{\dot H^{1/2}(\mathbb R)}^2+\|\chi_R'\|_{L^2(\mathbb R)}\le C/\log R.\qedhere
\end{equation}
\end{proof}

\begin{lemma}\label{lem:equal-height-fibers}
If \(f\in\mathscr X\) and \(E(f)=0\), there are \(P>0\) and \(h\in\mathbb R\) such that
\begin{equation}\label{eq:chain-two-gaps}
e^{2f(w+h)}=-\frac{\sin(\pi w/P)}{\sin(\pi(w-a)/P)}, \qquad a=2P\theta/\pi\in(0,P).
\end{equation}
\end{lemma}
\begin{proof}
Since \(E(f)=0\), \(\operatorname{Im}w\) is constant on \(\mathcal B^{-1}(\lambda)\) for almost every \(\lambda\in\mathbb D\). If \(0\in\mathcal B(\mathbb H_+)\), there is a disk \(D_0\Subset\mathbb H_+\) with \(0\in\mathcal B(D_0)\). Since \(|\mathcal B(z+is)|\le C/s\) by Lemma~\ref{lem:chain-inner-data},
\[
\mathcal B(\{s>H\})\subset\mathcal B(D_0)
\]
for some \(H>\sup_{D_0}\operatorname{Im}w\). For every \(\lambda\) in the nonempty open set \(\mathcal B(\{s>H\})\), there are \(w_1\in D_0\) and \(w_2\in\mathbb H_+\) such that
\[
\mathcal B(w_1)=\mathcal B(w_2)=\lambda, \qquad \operatorname{Im}w_1<H<\operatorname{Im}w_2,
\]
contradicting \(E(f)=0\).

Since \(\mathcal B\ne0\) on \(\mathbb H_+\), \(-\log|\mathcal B|\) is positive and harmonic on \(\mathbb H_+\). By Lemma~\ref{lem:chain-inner-data}, it is continuous on \(\overline{\mathbb H_+}\) and zero on \(\mathbb R\). Hence~\(-\log|\mathcal B(z+is)|=2\pi s/P\) for some \(P>0\)~\cite[Theorem~7.22]{ABR}. For \(h=z_1-a\), \(\mathcal B(a+h)=\mathcal B(z_1)=1\), so
\[
\begin{aligned}
\mathcal B(w+h)&=e^{2\pi i(w-a)/P},\\
g(w+h)&=\frac{g_\infty-\bar g_\infty\mathcal B(w+h)}{1-\mathcal B(w+h)}=-\frac{\sin(\pi w/P)}{\sin(\pi(w-a)/P)}.
\end{aligned}
\]
\end{proof}

\begin{proof}[Proof of Theorem~\ref{thm:transverse-scherk-classification}]
Put \(f=f_M\), so \(E(f)<\infty\) by \eqref{eq:complete-height-energy}. Let \(f(\cdot+h_n)\to f_*\) with \(|h_n|\to\infty\). If \(f_*(w_1)=f_*(w_2)\) is a regular value and \(\operatorname{Im}w_1\ne\operatorname{Im}w_2\), then for some \(L,c>0\) and all large \(n\),
\[
c\le\int_{S_0}\cosh(2\alpha) \sum_{\substack{j\ne k\\|x_j-h_n|<L}} a_{jk}(y_j-y_k)^2\mathrm d A\longrightarrow0,
\]
a contradiction. Thus \(E(f_*)=0\), and \eqref{eq:infinite-height-identity} implies that \(\mathcal V_{f_*}'=0\). By Proposition~\ref{prop:logarithmic-height-estimate} and \eqref{eq:complete-height-energy},
\[
0\le E(f)\le2\liminf_{R\to\infty}W_{\psi_R}=0.
\]
By Lemma~\ref{lem:equal-height-fibers}, we may assume that \eqref{eq:chain-two-gaps} holds with \(h=0\). The Weierstrass representation \eqref{eq:chain-weierstrass} is therefore that of a singly periodic Scherk surface, with
\[
X_M(w+P)-X_M(w)=X_M(P)-X_M(0)=Pe_3.\qedhere
\]
\end{proof}

\section{Convergence near a plane of multiplicity two}\label{sec:height-limits}

\begin{theorem}\label{thm:catenoid-classification}
Suppose that every tangent cone of \(M\) at infinity is a plane of multiplicity two. Then \(M\) is a catenoid.
\end{theorem}

Until the proof of Theorem~\ref{thm:classification}, we assume that every tangent cone of \(M\) at infinity is a plane of multiplicity two.

\subsection{Local height estimates and compactness}\label{cat:sec:height-compactness}

In this section and Section~\ref{sec:nonconstant-heights}, use coordinates \(X=(x,h,z)\) and identify \(\{h=0\}\) with \(\mathbb R^2\) by \(\xi=(x,z)\). Let \(\pi(x,h,z)=(x,z)\). For a unit normal \(N=(N_x,N_h,N_z)\), we have \(\pi N=(N_x,N_z)\). Put
\[
D_r=\{\xi\in\mathbb R^2:|\xi|<r\}, \qquad Z_r=\{X:\pi X\in D_r, |h|<r\}.
\]
Let \(\Sigma\subset Z_4\) be a smooth properly embedded minimal surface without boundary, put \(\mu=\mathcal H^2\llcorner\Sigma\), and assume
\begin{equation}\label{cat:eq:local-assumptions}
\Sigma\subset\{|h|\le1\}, \qquad \mu(\Sigma)\le M_0, \qquad a^2=\int_\Sigma h^2 \mathrm d\mu\le1.
\end{equation}
Since \(|h|\le1\), \(\pi|_\Sigma:\Sigma\to D_4\) is proper.

\begin{lemma}\label{cat:lem:height-tilt}
\[
\sup_{\Sigma\cap Z_{7/2}}|h|\le Ca, \qquad \int_{\Sigma\cap Z_3}|\pi N|^2 \mathrm d\mu\le Ca^2.
\]
\end{lemma}
\begin{proof}
Since \(\Delta_\Sigma h=0\) and \(|\nabla_\Sigma h|^2=|\pi N|^2\), the estimates follow from the mean-value inequality \cite{MichaelSimon} and Caccioppoli's inequality.
\end{proof}

On \(D_4\), put
\[
\nu=\pi_\#\mu, \qquad S=\pi_\#((\pi N)\otimes(\pi N)\mu), \qquad \kappa=\operatorname{tr}S.
\]
For \(Y\in C_c^\infty(D_4;\mathbb R^2)\), \(\int_\Sigma\operatorname{div}_\Sigma(Y\circ\pi)\,\mathrm d\mu=0\), hence
\begin{subequations}\label{cat:eq:moments}
\begin{equation}\label{cat:eq:moments-divergence}
\nabla \nu=\operatorname{div}S.
\end{equation}
For \(\varphi\in C_c^\infty(D_4)\), \(\Delta_\Sigma\varphi=\Delta\varphi-\nabla^2\varphi(\pi N,\pi N)\). Multiplying \(\Delta_\Sigma h=0\) and \(\Delta_\Sigma h^2=2|\pi N|^2\) by \(\varphi\) and integrating by parts, we obtain
\begin{align}
\int_\Sigma h\Delta\varphi \mathrm d\mu &=\int_\Sigma h\nabla^2\varphi(\pi N,\pi N) \mathrm d\mu, \label{cat:eq:moments-height}\\
\int_\Sigma h^2\Delta\varphi \mathrm d\mu &=2\int\varphi \mathrm d\kappa +\int_\Sigma h^2\nabla^2\varphi(\pi N,\pi N) \mathrm d\mu. \label{cat:eq:moments-height-squared}
\end{align}
\end{subequations}

Fix a nonnegative radial \(\eta\in C_c^\infty(D_1)\) with \(\int\eta=1\), and put \(\eta_\rho(\xi)=\rho^{-2}\eta(\xi/\rho)\). For \(0<\rho<4\), define on \(D_{4-\rho}\)
\[
f_\rho=\eta_\rho*\nu, \qquad m_\rho(\xi)=\frac{\int_\Sigma h\eta_\rho(\xi-\pi X) \mathrm d\mu(X)}{f_\rho(\xi)} \quad \text{ where }f_\rho(\xi)>0.
\]

\begin{proposition}\label{cat:prop:center}
If \(0<\rho\le\sigma\le1\) and \(f_\rho,f_\sigma\ge c_0>0\) on \(D_{3/2}\), then
\begin{equation}\label{cat:eq:cubic}
\|m_\rho\|_{C^2(D_{3/2})}\le Ca, \qquad \|\Delta m_\rho\|_{L^\infty(D_{3/2})}+\|m_\rho-m_\sigma\|_{L^\infty(D_{3/2})}\le Ca^3.
\end{equation}
The constant \(C\) depends only on \(\eta,\rho,c_0\) and \(M_0\).
\end{proposition}
\begin{proof}
Since \(f_\rho\ge c_0\), \(\|m_\rho\|_{C^2(D_{3/2})}\le Ca\). By \eqref{cat:eq:moments} and Lemma~\ref{cat:lem:height-tilt}, \(\|\nabla f_\rho\|_{C^1(D_{3/2})}\le Ca^2\) and
\[
f_\rho\Delta m_\rho =\int_\Sigma(h-m_\rho)\nabla^2\eta_\rho(\xi-\pi X) (\pi N,\pi N) \mathrm d\mu -2\nabla f_\rho\cdot \nabla m_\rho=O(a^3).
\]
The radial function \(K\) defined by \(K(\sigma)=0\) and \(rK'(r)=\int_0^r t(\eta_\rho(t)-\eta_\sigma(t))\mathrm d t\) belongs to \(C_c^\infty(D_\sigma)\) and satisfies
\begin{equation}\label{cat:eq:uniform-radial-potential}
\Delta K=\eta_\rho-\eta_\sigma, \qquad |\nabla^2 K(\xi)|\le\frac{C_\eta}{(\rho+|\xi|)^2}.
\end{equation}
With \(\varphi=K(\xi-\cdot)\), \eqref{cat:eq:moments} reads
\begin{equation}\label{cat:eq:center-comparison}
f_\rho(m_\rho-m_\sigma) =\int_\Sigma(h-m_\sigma)\nabla^2K(\xi-\pi X) (\pi N,\pi N) \mathrm d\mu=O(a^3).\qedhere
\end{equation}
\end{proof}

For an integral \(2\)-varifold \(V\) in \(\pi^{-1}(D_4)\), put \(\mu=\|V\|\) and let \(N\) be a unit normal to its approximate tangent plane. Set
\[
E_0(V,D_r(\xi)) =\frac2{\pi r^2}\int_{\pi^{-1}(D_r(\xi))}|\pi N|^2\mathrm d\mu.
\]
For \(\lambda>0\), put
\[
K_\lambda=\left\{\xi\in D_1: \sup_{\xi\in D_r(y)\subset D_4}E_0(V,D_r(y))\le\lambda\right\}.
\]

\begin{proposition}[Lipschitz approximation]
\label{cat:prop:lipschitz-approximation}
There are \(\varepsilon_0,\lambda_0,C_0>0\) such that, if
\begin{equation}\label{cat:fold:eq:relative-BDD-hypotheses}
\begin{gathered}
\delta V=0\text{ in }\pi^{-1}(D_4), \qquad E_0(V,D_4)\le\varepsilon_0,\\
\frac{\mu(\pi^{-1}(D_4))}{|D_4|}\le3, \qquad \left|\frac{\mu(\pi^{-1}(D_3))}{|D_3|}-2\right|<\frac12,
\end{gathered}
\end{equation}
then for every \(0<\lambda\le\lambda_0\) there are Lipschitz functions \(u_1\le u_2\) on \(D_1\) with \(\operatorname{Lip}u_i\le C_0\lambda^{1/4}\) and \(\Gamma=\sum_{i=1}^2|\operatorname{graph}u_i|=V\) over \(K_\lambda\). For \(2\le q<\infty\),
\begin{equation}\label{cat:eq:lipschitz-exceptional-set}
|D_1\setminus K_\lambda| +\mu\left(\pi^{-1}(D_1\setminus K_\lambda)\right) \le C_0^q\lambda^{-q/2} \int_{\pi^{-1}(D_4)}|\pi N|^q\mathrm d\mu.
\end{equation}
If \(\operatorname{spt} V\cap\pi^{-1}(D_1)\subset\{|h|\le H\}\), the functions may also be chosen with \(|u_i|\le H\).

For each \(\xi\in D_1\setminus K_\lambda\), there is a disk \(\xi\in D_r(y)\subset D_4\) with \(E_0(V,D_r(y))>\lambda\) and
\begin{equation}\label{cat:fold:eq:relative-cover-mass}
\mu(\pi^{-1}(D_r(y))) +\mu(\pi^{-1}(D_1\cap D_{5r}(y)))\le C_0r^2.
\end{equation}
If \(V=|\Sigma|\) for a smooth properly embedded surface \(\Sigma\), then for every \(\xi\in K_\lambda\) the fiber \(\Sigma\cap\pi^{-1}(\xi)\) contains exactly two points, both with \(|N_h|>1/2\).
\end{proposition}
\begin{proof}
For \(m=Q=2\) and \(n=1\), \cite[Theorem 3.3, (3.11)]{BDD} gives the functions \(u_1\le u_2\) and \eqref{cat:eq:lipschitz-exceptional-set}. If \(|h|\le H\), replace \(u_i\) by \(\max\{-H,\min\{u_i,H\}\}\).

The disks chosen in Step~5 of the proof of \cite[Theorem~3.3]{BDD} satisfy \eqref{cat:fold:eq:relative-cover-mass} if \(D_{5r}(y)\subset D_{7/2}\); otherwise use \(r\ge5/12\) and \(\mu(\pi^{-1}(D_4))\le48\pi\).

For \(V=|\Sigma|\), we have \(\Theta(V,X)=1\) for \(X\in\Sigma\), so by \cite[(1.11)]{BDD}
\[
\#\bigl(\Sigma\cap\pi^{-1}(\xi)\bigr)=2
\qquad\text{for }\xi\in K_\lambda.
\]
Letting \(s\downarrow0\) in \cite[(1.9)]{BDD} gives \(|\pi N|\le C_0\lambda^{1/4}<\sqrt3/2\) for small \(\lambda_0\).
\end{proof}

For the rest of this subsection, let \(\Sigma_j\) be smooth properly embedded minimal surfaces without boundary in \(Z_{16}\), and suppose \(|\Sigma_j|\to2|P|\) there, with \(P=\{h=0\}\). Let \(a_j>0\) tend to zero and satisfy
\begin{equation}\label{cat:eq:compat-local-bounds}
\sup_{\Sigma_j\cap Z_{16}}|h|\le Ca_j, \qquad \mu_j(Z_{16})\le C, \qquad \int_{\Sigma_j\cap Z_{16}}|\pi N|^2 \mathrm d\mu_j \le Ca_j^2.
\end{equation}
Since \(\sup_{\Sigma_j}|h|\le Ca_j<16\) for large \(j\), \(|\Sigma_j|\) is stationary in \(\pi^{-1}(D_{16})\). Write \(\nu_j,S_j,\kappa_j\) for the measures \(\nu,S,\kappa\) with \(\mu=\mu_j\).

Fix \(3<r<4\) and \(0<\lambda\le\lambda_0\). By Proposition~\ref{cat:prop:lipschitz-approximation} with \(q=2\), there are functions \(u_{j,-}\le u_{j,+}\) on \(D_r\) with \(|u_{j,\pm}|\le Ca_j\) and \(\operatorname{Lip}u_{j,\pm}\le C\lambda^{1/4}\) such that \(|\Sigma_j|=\sum_\pm|\operatorname{graph}u_{j,\pm}|\) over \(D_r\setminus E_j\) and
\begin{equation}\label{cat:eq:compat-bad-two}
|E_j|+\mu_j(\pi^{-1}(E_j))\le C_\lambda a_j^2.
\end{equation}

\begin{lemma}
\label{cat:lem:compat-weak}
After passing to a subsequence, there are \(m,v\in W^{1,2}_{\mathrm{loc}}(D_4)\) with \(v\ge0\) and \(\Delta m=0\). For every \(0<r<4\) and \(0<\lambda\le\lambda_0\), every pair \(u_{j,-}\le u_{j,+}\) on \(D_r\) from Proposition~\ref{cat:prop:lipschitz-approximation} with \(|u_{j,\pm}|\le Ca_j\) satisfies
\[
a_j^{-1}u_{j,\pm}\longrightarrow m\pm v \quad \text{ strongly in }L^2(D_r) \text{ and weakly in }W^{1,2}(D_r).
\]
\end{lemma}
\begin{proof}
By the area formula and \eqref{cat:eq:compat-bad-two},
\begin{equation}\label{cat:eq:compat-moment-errors}
\int_{D_r}(|\nabla u_{j,+}|^2+|\nabla u_{j,-}|^2)\mathrm d A \le C_\lambda a_j^2.
\end{equation}
After passing to a subsequence, \(a_j^{-1}u_{j,\pm}\to m\pm v\) in \(L^2_{\mathrm{loc}}(D_4)\). Two choices \(u_{j,\pm},\widetilde u_{j,\pm}\) on \(D_r\) agree outside \(E_j\cup\widetilde E_j\), so
\[
a_j^{-1}\|u_{j,\pm}-\widetilde u_{j,\pm}\|_{L^2(D_r)} \le C|E_j\cup\widetilde E_j|^{1/2}=O(a_j).
\]
The area formula and \eqref{cat:eq:compat-bad-two}--\eqref{cat:eq:compat-moment-errors} give, locally in total variation,
\begin{equation}\label{cat:eq:compat-strong-moments}
\frac{\pi_\#(h\mu_j)}{a_j}\longrightarrow2m \mathcal L^2, \qquad \frac{\pi_\#(h^2\mu_j)}{a_j^2}\longrightarrow 2(m^2+v^2) \mathcal L^2.
\end{equation}
Dividing \eqref{cat:eq:moments-height} by \(a_j\) and letting \(j\to\infty\) proves \(\Delta m=0\). For fixed \(\rho\),
\[
f_{j,\rho}\longrightarrow2, \qquad \frac{m_{j,\rho}}{a_j}\longrightarrow\eta_\rho*m=m \quad \text{ in }C^\infty_{\mathrm{loc}}(D_{4-\rho}).\qedhere
\]
\end{proof}

Fix \(0<\rho<1\), define \(m_j=m_{j,\rho}\) on \(D_{16-\rho}\), and write \(w_j=h-m_j\circ\pi\).

\begin{proposition}
\label{cat:prop:compat-high-graphs}
Fix \(K\Subset D_{16-\rho}\) and \(\delta>0\). For all large \(j\) and every \(\xi\in K\cap\pi(\{|w_j|\ge\delta a_j\})\), \(\Sigma_j\) consists of two smooth graphs \(h_{j,-}<m_j<h_{j,+}\) over a neighborhood of \(\xi\). Uniformly for such \(\xi\),
\[
|\nabla h_{j,+}|+|\nabla h_{j,-}|+ \frac{|h_{j,+}+h_{j,-}-2m_j|}{a_j}\longrightarrow0.
\]
\end{proposition}
\begin{proof}
Fix \(r_0>0\). For \(B_{2r_0}(p)\Subset Z_{16}\) and \(0<r<r_0/2\), take \(\chi=1\) on \(B_r(p)\), \(\chi=0\) outside \(B_{r_0}(p)\), and \(\chi=\log(r_0/|X-p|)/\log(r_0/r)\) between them. Since \(\Delta_{\Sigma_j}h=0\) and \(\mu_j(B_s(p))\le Cs^2\),
\begin{equation}\label{cat:eq:compat-log-tilt}
\int_{B_r(p)}|\pi N|^2 \mathrm d\mu_j \le4\|h\|_{L^\infty(\Sigma_j\cap B_{r_0}(p))}^2\int|\nabla_{\Sigma_j}\chi|^2 \mathrm d\mu_j \le\frac{Ca_j^2}{\log(r_0/r)}.
\end{equation}
For \(r\ge a_j\), \(\Sigma_j\cap\pi^{-1}(D_r(\xi))\subset B_{Cr}(\xi)\). Thus \eqref{cat:eq:compat-log-tilt} gives \(\kappa_j(D_r(\xi))\le Ca_j^2/\log(r_0/r)\) for \(a_j\le r\le cr_0\), with a fixed \(c\in(0,1)\).

Fix \(L>0\). By \eqref{cat:eq:moments}, \eqref{cat:eq:uniform-radial-potential}, and \eqref{cat:eq:compat-log-tilt},
\[
\begin{aligned}
|f_{j,La_j}(\xi)-f_{j,\rho}(\xi)| &\le C_L\int\frac{\mathrm d\kappa_j(y)}{(|\xi-y|+a_j)^2} \\
&\le C_L\left(\frac{\kappa_j(D_{\sqrt{a_j}}(\xi))}{a_j^2} +\frac{\kappa_j(D_{16})}{a_j}\right)\longrightarrow0.
\end{aligned}
\]
Since \(f_{j,\rho}\to2\), \eqref{cat:eq:center-comparison} and \(|h-m_j|\le Ca_j\) give
\begin{equation}\label{cat:eq:compat-shrinking-center}
f_{j,La_j}=2+o(1), \qquad m_{j,La_j}-m_j=o(a_j),
\end{equation}
locally uniformly on \(D_{16-\rho}\).

For \(\xi_j\in K\), the surfaces \(a_j^{-1}(\Sigma_j-\xi_j-m_j(\xi_j)e_h)\) have uniformly bounded area on each fixed ball and satisfy \(|h|\le C\). After passing to a subsequence, these surfaces converge to an integral stationary varifold \(V\). By \eqref{cat:eq:compat-log-tilt}, \(\pi N=0\) almost everywhere on \(V\). For \(\varphi\in C_c^\infty(\mathbb R^3)\) and \(i=x,z\),
\[
0=\delta V(\varphi e_i)=\int\partial_i\varphi\,\mathrm d\|V\|.
\]
Thus \(V=\sum_\ell m_\ell|\{h=t_\ell\}|\) with positive integers \(m_\ell\). Equation~\eqref{cat:eq:compat-shrinking-center} gives \(\sum_\ell m_\ell=2\) and \(\sum_\ell m_\ell t_\ell=0\). Hence, for some \(t\ge0\),
\begin{equation}\label{cat:eq:compat-height-scale-limit}
V=|\{h=t\}|+|\{h=-t\}|.
\end{equation}

For \(X_j\in\Sigma_j\) with \(\pi(X_j)=\xi_j\) and \(|w_j(X_j)|\ge\delta a_j\), Lemma~\ref{lem:varifold-limits} gives \(|w_j(X_j)|/a_j\to t\ge\delta\) and smooth convergence to \(V\) near both planes. Since \(|h|\le Ca_j\) on \(\Sigma_j\), \(\Sigma_j\cap\pi^{-1}(D_{a_j}(\xi_j))\) consists of two smooth graphs for~large~\(j\).
\end{proof}

\begin{proposition}
\label{cat:prop:compat-pressure}
For every disk \(D\Subset D_4\),
\begin{equation}\label{cat:eq:compat-fourth-tilt}
\int_{\pi^{-1}(D)}|\pi N|^4 \mathrm d\mu_j=o(a_j^2), \qquad |\nu_j-2\mathcal L^2-\kappa_j/2|(D)=o(a_j^2).
\end{equation}
\end{proposition}
\begin{proof}
Since \(\Delta_{\Sigma_j}w_j=-\Delta m_j+ \nabla^2m_j(\pi N,\pi N)\), Proposition~\ref{cat:prop:center} gives
\[
\int_{\pi^{-1}(D_4)}|\Delta_{\Sigma_j}w_j| \mathrm d\mu_j\le Ca_j^3, \qquad \int_{\pi^{-1}(D_4)}|\nabla_{\Sigma_j}w_j|^2 \mathrm d\mu_j\le Ca_j^2.
\]
Choose \(\chi\in C_c^\infty(D_4)\) with \(\chi=1\) near \(\overline D\). Multiplying \(\Delta_{\Sigma_j}w_j\) by \(\chi^2\max\{-\delta a_j,\min\{w_j,\delta a_j\}\}\) and integrating by parts gives
\begin{equation}\label{cat:eq:compat-band-energy}
\int_{\{|w_j|<\delta a_j\}}\chi^2|\nabla_{\Sigma_j}w_j|^2 \mathrm d\mu_j\le C\delta a_j^2.
\end{equation}
By Proposition~\ref{cat:prop:compat-high-graphs}, \(|\pi N|\to0\) uniformly on \(\pi^{-1}(D)\cap\{|w_j|\ge\delta a_j\}\). On \(\{|w_j|<\delta a_j\}\), \(|\pi N|^4\le C|\nabla_{\Sigma_j}w_j|^2+Ca_j^4\). Thus \eqref{cat:eq:compat-local-bounds} and \eqref{cat:eq:compat-band-energy} give
\[
\limsup_{j\to\infty}a_j^{-2}\int_{\pi^{-1}(D)}|\pi N|^4\mathrm d\mu_j\le C_D\delta.
\]
Letting \(\delta\downarrow0\) proves the first estimate in \eqref{cat:eq:compat-fourth-tilt}.

By Proposition~\ref{cat:prop:lipschitz-approximation} with \(q=4\) and a fixed \(0<\lambda\le\lambda_0\), there is \(E_j\subset D\) such that
\[
|E_j|+\mu_j(\pi^{-1}(E_j))=o(a_j^2).
\]
For almost every \(\xi\in D\setminus E_j\), \(\#(\Sigma_j\cap\pi^{-1}(\xi))=2\). Since \(0\le1-|N_h|-\tfrac12|\pi N|^2\le\tfrac12|\pi N|^4\), the area formula gives
\[
|\nu_j-2\mathcal L^2-\kappa_j/2|(D) \le\mu_j(\pi^{-1}(E_j))+2|E_j| +\tfrac12\int_{\pi^{-1}(D)}|\pi N|^4 \mathrm d\mu_j =o(a_j^2).\qedhere
\]
\end{proof}

Pass to a further subsequence with \(S_j/a_j^2\rightharpoonup\mathcal S\) locally as measures. Put \(H=m^2+v^2\) and \(\mathcal S^0=\mathcal S-\tfrac12(\operatorname{tr}\mathcal S)I\). By \eqref{cat:eq:compat-local-bounds},
\[
a_j^{-2}\int_{\Sigma_j}h^2|\pi N|^2 \mathrm d\mu_j \le Ca_j^2\longrightarrow0.
\]
By \eqref{cat:eq:moments}, \eqref{cat:eq:compat-strong-moments}, and \eqref{cat:eq:compat-fourth-tilt},
\begin{equation}\label{cat:eq:compat-limit-moments}
\frac{\nu_j-2\mathcal L^2}{a_j^2}\rightharpoonup\tfrac12\Delta H, \qquad \operatorname{tr}\mathcal S=\Delta H\ge0, \qquad \operatorname{div}\mathcal S^0=0.
\end{equation}

\begin{lemma}
\label{cat:lem:compat-area-cap}
Suppose \(\mu_j(B_r)\le2\pi r^2\) for \(0<r<4\). Then
\begin{equation}\label{cat:eq:compat-height-monotonicity}
r\longmapsto r^{-4}\int_{D_r}H\mathrm d A \quad \text{ is nonincreasing on }(0,4).
\end{equation}
\end{lemma}
\begin{proof}
For \(0<r<4\),
\[
\int(r^2-|X|^2)_+ \mathrm d\mu_j =2\int_0^r s\mu_j(B_s) \mathrm ds\le\pi r^4.
\]
Since \(|h|\le Ca_j\), \eqref{cat:eq:compat-strong-moments} gives
\[
\begin{aligned}
0&\le\int\left[(r^2-|X|^2)_+-(r^2-|\xi|^2)_+ +h^2\mathbf1_{D_r}(\xi)\right] \mathrm d\mu_j \\
&\le\int_{\{r^2-C^2a_j^2<|\xi|^2<r^2\}}h^2 \mathrm d\mu_j=o(a_j^2).
\end{aligned}
\]
By \eqref{cat:eq:compat-limit-moments}, for almost every \(r\),
\begin{equation}\label{cat:eq:compat-area-cap-limit}
\begin{aligned}
r^5\frac{d}{dr}\left(r^{-4}\int_{D_r}H\mathrm d A\right) &=\tfrac12\int(r^2-|\xi|^2)_+ \mathrm d(\Delta H)-2\int_{D_r}H\mathrm d A \\
&=\lim_j a_j^{-2}\left[\int(r^2-|X|^2)_+ \mathrm d\mu_j-\pi r^4\right]\le0.\qedhere
\end{aligned}
\end{equation}
\end{proof}

\subsection{Approximation by planes}\label{cat:sec:height-and-planes}

Fix a smooth function \(\omega(|X|)\) that is positive on \(B_1\) and zero outside \(\overline B_1\), and set
\[
C(R)=R^{-4}\int_M X\otimes X \omega(|X|/R) \mathrm d\mu.
\]
Choose a unit eigenvector \(\mathbf{n}_R\) for the smallest eigenvalue of \(C(R)\) and put \(P_R=\mathbf{n}_R^\perp\). Let \(\Pi_R\) be the orthogonal projection onto \(P_R\), and define
\begin{equation}\label{cat:eq:aR}
a_R^2 =R^{-4}\int_{M\cap B_{8R}}\operatorname{dist}(X,P_R)^2 \mathrm d\mu,\quad \mathsf H_R =Ra_R.
\end{equation}
For large \(R\), \(a_R>0\), since \(M\) is not a plane.

\begin{lemma}
\label{cat:lem:initial-comparison}
For all sufficiently large \(R\), \(P_R\) is unique and \(\Pi_R\) depends smoothly on \(R\). For \(R\le t\le2R\),
\[
\|\Pi_t-\Pi_{2R}\|\le Ca_{2R}, \qquad a_t\le Ca_{2R}.
\]
As \(R\to\infty\), \(a_R\to0\). In coordinates with \(P_R=\{h=0\}\), \(|R^{-1}M|\to2|\{h=0\}|\) as varifolds, and \(\sup_{R^{-1}M\cap B_L}|h|\to0\) for each fixed \(L\).
\end{lemma}
\begin{proof}
By the monotonicity formula and Lemma~\ref{lem:varifold-limits}, for any \(R_j\to\infty\) a subsequence satisfies \(|R_j^{-1}M|\to2|P|\) for a plane \(P\). Then
\[
C(R_j)\longrightarrow c_\omega\Pi_P, \qquad c_\omega>0.
\]
The smallest eigenvalue of \(C(R)\) is therefore simple for all large \(R\). Since \(C(R)\) is smooth, so is \(\Pi_R\), and \(P_{R_j}\to P\).

Lemma~\ref{lem:varifold-limits} and \(P_{R_j}\to P\) give \(\sup_{R_j^{-1}M\cap B_L}\operatorname{dist}(X,P_{R_j})\to0\) for each fixed \(L\), and \(a_{R_j}\to0\) by \eqref{eq:area-bound}.

For \(R\le t\le2R\), \(C(t)\ge c\Pi_t\). By \eqref{cat:eq:aR} and \eqref{eq:area-bound},
\[
\begin{aligned}[b]
\|\Pi_t-\Pi_{2R}\|^2 &\le C \mathbf{n}_{2R}\cdot C(t)\mathbf{n}_{2R}\le Ca_{2R}^2, \\
a_t^2&\le C\left(a_{2R}^2+\|\Pi_t-\Pi_{2R}\|^2\right) \le Ca_{2R}^2.
\end{aligned}\qedhere
\]
\end{proof}

\begin{proposition}\label{cat:prop:plane-drift}
For all sufficiently large \(R\),
\begin{align}
\left\|\frac{d\Pi_R}{d\log R}\right\| &\le Ca_R^3,\label{cat:eq:q-drift} \\
\|\Pi_{2R}-\Pi_R\|&\le Ca_{2R}^3.\label{cat:eq:drift}
\end{align}
\end{proposition}
\begin{proof}
Put \(r=|X|\) and \(s=X\cdot N\). For fixed \(u,v\in\mathbb R^3\),
\[
\begin{aligned}
&\operatorname{div}_M\left((X\cdot u)(X\cdot v)X\right) \\
={}& 4(X\cdot u)(X\cdot v) -s\left((X\cdot u)(N\cdot v)+(X\cdot v)(N\cdot u)\right), \\
&X\cdot\nabla_Mr=r-s^2/r.
\end{aligned}
\]
By the first variation formula,
\begin{equation}\label{cat:eq:covariance-derivative}
\begin{aligned}
C'(R)&=-R^{-5}\int_M X\otimes X\left(4\omega(r/R)+\frac rR\omega'(r/R)\right) \mathrm d\mu \\
&=-R^{-5}\int_M\left[ s(X\otimes N+N\otimes X)\omega(r/R) +\frac{\omega'(r/R)}{Rr}s^2X\otimes X\right] \mathrm d\mu.
\end{aligned}
\end{equation}

For a fixed \(R\), take coordinates with \(\mathbf{n}_R=e_h\). For \(e\in\{e_x,e_z\}\), put \(Y_e=(X\cdot e)e_h-he\) and \(\Psi(r)=-\int_r^\infty t\omega(t/R) \mathrm dt\). Then
\[
\nabla_M\Psi(r)=\omega(r/R)(X-sN), \qquad \operatorname{div}_M Y_e=0, \qquad X\cdot Y_e=0.
\]
Applying the first variation formula to \(\Psi(r)Y_e\) gives
\begin{equation}\label{cat:eq:torque}
\int_M s\left((X\cdot e)N_h-h(N\cdot e)\right)\omega(r/R) \mathrm d\mu=0.
\end{equation}
By \eqref{cat:eq:covariance-derivative} and \eqref{cat:eq:torque},
\[
\pi C'(R)e_h=-R^{-5}\int_M h \left(2s\pi N\omega(r/R)+\xi\frac{\omega'(r/R)}{Rr}s^2\right) \mathrm d\mu.
\]
By Lemma~\ref{cat:lem:height-tilt}, \(\|h\|_{L^\infty(M\cap B_R)}\le Ca_RR\) and \(\|\pi N\|_{L^2(M\cap B_R)}\le Ca_RR\). Since \(|s|\le|h|+R|\pi N|\) and \(|\omega'(r/R)|\le Cr/R\),
\[
R|\pi C'(R)e_h|\le CR^{-5}\|h\|_{L^\infty(M\cap B_R)} \int_{M\cap B_R}(h^2+R^2|\pi N|^2) \mathrm d\mu \le Ca_R^3.
\]

Choose \(\mathbf{n}_R\) smoothly in \(R\) and write \(C(R)\mathbf{n}_R=\lambda_R\mathbf{n}_R\). Differentiating and projecting onto \(P_R\) gives
\[
(C(R)-\lambda_RI)\mathbf{n}_R'=-\Pi_R C'(R)\mathbf{n}_R.
\]
Since \(\lambda_R\to0\) and the other two eigenvalues stay bounded away from zero (Lemma~\ref{cat:lem:initial-comparison}), we have \(C(R)-\lambda_RI\ge cI\) on \(P_R\). Since \(\mathbf{n}_R'\in P_R\) and \(\Pi_R=I-\mathbf{n}_R\otimes \mathbf{n}_R\),
\[
R\|\Pi_R'\|\le CR|\mathbf{n}_R'| \le CR|\pi C'(R)e_h|\le Ca_R^3.
\]
Integrating over \([R,2R]\) and using \(a_t\le Ca_{2R}\) from Lemma~\ref{cat:lem:initial-comparison} proves \eqref{cat:eq:drift}.
\end{proof}

For a sequence \(R_j\to\infty\), put \(a_j=a_{R_j}\) and rotate \(\Sigma_j=R_j^{-1}M\) so that \(P_{R_j}=\{h=0\}\).

\begin{lemma}
\label{cat:lem:compat-controlled-scales}
Suppose that for every integer \(k\ge0\) there is \(C_k<\infty\) such that
\begin{equation}\label{cat:eq:compat-outer-bounds}
a_{2^kR_j}\le C_ka_j \quad \text{ for all sufficiently large }j.
\end{equation}
After passing to a subsequence, the convergence in Lemma~\ref{cat:lem:compat-weak} holds on \(D_r\) for every \(r>0\). The identities \eqref{cat:eq:compat-limit-moments} hold on \(\mathbb R^2\), and \eqref{cat:eq:compat-height-monotonicity} holds for every \(r>0\). Locally uniformly for \(s>0\),
\begin{equation}\label{cat:eq:scale-model-compatibility}
\frac{a_{sR_j}}{a_j}\longrightarrow \left(\frac2{s^4}\int_{D_{8s}}(m^2+v^2)\mathrm d A\right)^{1/2}.
\end{equation}
\end{lemma}
\begin{proof}
For each interval \(I\Subset(0,\infty)\), Lemma~\ref{cat:lem:initial-comparison} and \eqref{cat:eq:compat-outer-bounds} give \(a_{sR_j}\le C_Ia_j\) for \(s\in I\). Integrating \eqref{cat:eq:q-drift} gives
\[
\sup_{s\in I}\|\Pi_{sR_j}-\Pi_{R_j}\|\le C_Ia_j^3.
\]
For fixed \(L>0\), \eqref{cat:eq:aR} and \eqref{eq:area-bound} imply
\begin{equation}\label{cat:eq:compat-expanded-height}
\int_{\Sigma_j\cap B_L}h^2 \mathrm d\mu_j\le C_La_j^2.
\end{equation}
By Lemma~\ref{cat:lem:height-tilt},
\[
\sup_{\Sigma_j\cap Z_L}|h|\le C_La_j, \qquad \int_{\Sigma_j\cap Z_L}|\pi N|^2 \mathrm d\mu_j \le C_La_j^2.
\]
After passing to a subsequence, Lemmas~\ref{cat:lem:compat-weak} and \ref{cat:lem:compat-area-cap} applied on each fixed disk give limits \(m,v\) on \(\mathbb R^2\) satisfying \eqref{cat:eq:compat-limit-moments} and \eqref{cat:eq:compat-height-monotonicity}.

By \eqref{cat:eq:compat-strong-moments},
\[
a_j^{-2}\int_{\Sigma_j\cap B_r}h^2 \mathrm d\mu_j \longrightarrow2\int_{D_r}(m^2+v^2)\mathrm d A
\]
locally uniformly for \(r>0\). Since \(\|\Pi_{sR_j}-\Pi_{R_j}\|\le C_Ia_j^3\), \eqref{eq:area-bound} gives for \(s\in I\)
\[
\frac{a_{sR_j}}{a_j} =\frac1{s^2a_j} \left(\int_{\Sigma_j\cap B_{8s}}h^2 \mathrm d\mu_j\right)^{1/2} +O_I(a_j^2).\qedhere
\]
\end{proof}

\begin{proposition}
\label{cat:prop:compat-upper-doubling}
For every \(\varepsilon>0\), all sufficiently large \(R\) satisfy
\begin{equation}\label{cat:eq:compat-upper-doubling}
a_{2R}\le(1+\varepsilon)a_R.
\end{equation}
\end{proposition}
\begin{proof}
Suppose \(a_{2R_j}>Ba_{R_j}\) for some \(B>1\) and \(R_j\to\infty\). Let \(k_j\) maximize \(a_{2^kR_j}/B^k\) over integers \(k\ge0\). Then \(k_j\ge1\), and \(T_j=2^{k_j}R_j\) satisfies, for every integer \(\ell\ge0\),
\begin{equation}\label{cat:eq:compat-selected-scales}
a_{T_j/2}\le a_{T_j}/B, \qquad a_{2^\ell T_j}\le B^\ell a_{T_j}.
\end{equation}
By Lemma~\ref{cat:lem:compat-controlled-scales} applied at \(T_j\) and \eqref{cat:eq:compat-height-monotonicity}, a subsequence of \(s\mapsto a_{sT_j}/a_{T_j}\) converges locally uniformly to a nonincreasing function on \((0,\infty)\). Hence \eqref{cat:eq:compat-selected-scales} gives
\[
1\le\lim_j\frac{a_{T_j/2}}{a_{T_j}}\le B^{-1}<1.\qedhere
\]
\end{proof}

\begin{lemma}\label{cat:lem:physical-height}
\(\mathsf H_R\to\infty\).
\end{lemma}
\begin{proof}
Suppose \(\mathsf H_{R_j}\le C\) and \(\mathbf{n}_{R_j}\to\mathbf{n}\) for some \(R_j\to\infty\). By Lemma~\ref{cat:lem:height-tilt}, \(|\mathbf{n}\cdot X|\le C\) on \(M\). For \(r>0\), as in \eqref{cat:eq:compat-log-tilt},
\[
\int_{M\cap B_r}|\nabla_M(\mathbf{n}\cdot X)|^2 \mathrm d\mu \le\lim_{R\to\infty}\frac{C}{\log(R/r)}=0.
\]
By properness, \(M\) is a plane, contrary to the assumption that \(M\) has density two at infinity.
\end{proof}

Let \(R_j\to\infty\). Iterating \eqref{cat:eq:compat-upper-doubling} gives \eqref{cat:eq:compat-outer-bounds}. Fix \(0<\lambda\le\lambda_0\) and take \(u_{j,-}\le u_{j,+}\) from Lemma~\ref{cat:lem:compat-weak}. After passing to a subsequence, let \(m\pm v\) be the limits of \(u_{j,\pm}/a_j\) in Lemma~\ref{cat:lem:compat-controlled-scales}.

\begin{lemma}\label{cat:lem:compat-zero-mean}
\(m=0\) and
\begin{equation}\label{cat:eq:compat-gap-growth}
\|v\|_{L^\infty(D_R)}\le C(1+R), \qquad 2\int_{D_8}v^2\mathrm d A=1.
\end{equation}
\end{lemma}
\begin{proof}
For \(H=m^2+v^2\), \eqref{cat:eq:scale-model-compatibility} at \(s=1\) gives \(2\int_{D_8}H\mathrm d A=1\). Since \(\Delta H\ge0\) by \eqref{cat:eq:compat-limit-moments}, the mean-value inequality and \eqref{cat:eq:compat-height-monotonicity} give \(H(\xi)\le C(1+|\xi|)^2\). Thus the harmonic function \(m\) is affine.

Since \(\pi C(R_j)e_h=0\) and \(\omega(|X|)-\omega(|\xi|)=O(h^2)\), \eqref{cat:eq:compat-strong-moments} implies
\[
0=\lim_j\frac1{a_j}\int_{\Sigma_j\cap Z_2} h\xi\omega(|X|) \mathrm d\mu_j =2\nabla m\int x^2\omega(|\xi|)\mathrm d A.
\]
Thus \(m\) is constant.

Since \(R_ja_j\to\infty\) by Lemma~\ref{cat:lem:physical-height}, a subsequence satisfies \(|a_j^{-1}\Sigma_j|\to2|P|\) for a plane \(P\) through \(0\). Also \(m_j(0)/a_j\to m(0)\) by Lemma~\ref{cat:lem:compat-weak}. Equation~\eqref{cat:eq:compat-height-scale-limit} at \(\xi=0\) gives, for some \(t\ge0\),
\[
2|P|=|\{h=m(0)+t\}|+|\{h=m(0)-t\}|.
\]
Since \(0\in P\), \(t=m(0)=0\).
\end{proof}

\subsection{Classification of the limiting function \texorpdfstring{\(v\)}{v}}\label{cat:sec:height-classification}

For a \(2\times2\) matrix \(F\), write \(F^0=F-\tfrac12(\operatorname{tr}F)I\).

\begin{lemma}[{\cite{Delort}}]
\label{cat:lem:compat-delort}
Let \(\Omega\subset\mathbb R^2\) be open. Suppose \(w_j\rightharpoonup w\) in \(W^{1,2}_{\mathrm{loc}}(\Omega)\) and \(A_j\to I\) in \(L^\infty_{\mathrm{loc}}(\Omega;\mathbb R^{2\times2})\), with \(\operatorname{div}(A_j\nabla w_j)\ge0\). Then \(\Delta w\ge0\) and
\begin{equation}\label{cat:eq:compat-delort}
(\nabla w_j\otimes \nabla w_j)^0\longrightarrow(\nabla w\otimes \nabla w)^0 \quad \text{ in distributions on }\Omega.
\end{equation}
\end{lemma}
\begin{proof}
Since \(A_j\nabla w_j\rightharpoonup \nabla w\) locally in \(L^2\), we have \(\Delta w\ge0\). Fix \(U\Subset\Omega\) and \(0\le\chi\in C_c^\infty(\Omega)\) with \(\chi=1\) near \(\overline U\). The nonnegative measures \(\lambda_j=\chi\operatorname{div}(A_j\nabla w_j)\) on \(\mathbb R^2\) are supported in \(\operatorname{spt}\chi\) and have uniformly bounded mass and \(H^{-1}(\mathbb R^2)\) norm. They converge weakly to \(\lambda=\chi\Delta w\). Their potentials \(p_j=(2\pi)^{-1}\log|\cdot|*\lambda_j\) converge weakly in \(W^{1,2}_{\mathrm{loc}}(\mathbb R^2)\) to \(p=(2\pi)^{-1}\log|\cdot|*\lambda\).

For \(a\in\mathbb R^2\) and \(0<r<1\),
\[
\lambda(D_r(a))^2\log(1/r) \le\int_r^1\frac{\lambda(D_s(a))^2}{s} \mathrm ds \le2\pi\int_{D_1(a)}|\nabla p|^2.
\]
Thus \(\lambda\) has no atoms, and neither does any \(\lambda_j\).

Identify \(\xi=(x,z)\) with \(x+iz\in\mathbb C\). For \(\varphi\in C_c^\infty(\mathbb C)\), put \(C_\varphi(a)=\int_{\mathbb C}\varphi(\zeta)/(\zeta-a)\mathrm d A(\zeta)\). If \(\lambda_j\) has a smooth density,
\[
2\pi(\partial_xp_j-i\partial_zp_j)(\zeta)=\int_{\mathbb C}\frac{\mathrm d\lambda_j(a)}{\zeta-a}.
\]
Hence
\[
\begin{aligned}
4\pi^2\int_{\mathbb C}\varphi(\partial_xp_j-i\partial_zp_j)^2\mathrm d A &=\iint\left(\int_{\mathbb C}\frac{\varphi(\zeta)}{(\zeta-a)(\zeta-b)}\mathrm d A(\zeta)\right)\mathrm d\lambda_j(a)\mathrm d\lambda_j(b) \\
&=\iint\frac{C_\varphi(a)-C_\varphi(b)}{a-b}\mathrm d\lambda_j(a)\mathrm d\lambda_j(b).
\end{aligned}
\]
Since \(C_\varphi\) is smooth, \((C_\varphi(a)-C_\varphi(b))/(a-b)\) is bounded and continuous for \(a,b\in\operatorname{spt}\chi\), \(a\ne b\). The diagonal \(\{a=b\}\) has measure zero for both \(\lambda_j\otimes\lambda_j\) and \(\lambda\otimes\lambda\). Approximating \(\lambda_j\) and \(\lambda\) by measures with smooth densities gives the same identity. Since \(\lambda_j\otimes\lambda_j\rightharpoonup\lambda\otimes\lambda\), taking real and imaginary parts gives
\[
(\nabla p_j\otimes\nabla p_j)^0\longrightarrow(\nabla p\otimes\nabla p)^0 \quad \text{ in distributions}.
\]

On \(U\),
\[
\Delta(w_j-p_j)=\operatorname{div}((I-A_j)\nabla w_j) \longrightarrow0 \quad \text{in }H^{-1}_{\mathrm{loc}}(U).
\]
Since \(w_j-p_j\to w-p\) in \(L^2_{\mathrm{loc}}(U)\), interior elliptic estimates give
\[
\nabla(w_j-p_j)\longrightarrow\nabla(w-p) \quad \text{strongly in }L^2_{\mathrm{loc}}(U).
\]
Expanding \(\nabla w_j\otimes \nabla w_j\) proves \eqref{cat:eq:compat-delort} on \(U\).
\end{proof}

\begin{proposition}
\label{cat:prop:compat-local-stress}
\(v\) is subharmonic and satisfies
\begin{equation}\label{cat:eq:compat-gap-stress}
\operatorname{div}(2\nabla v\otimes \nabla v-|\nabla v|^2I)=0.
\end{equation}
\end{proposition}
\begin{proof}
Fix \(L>0\) and \(\delta>0\), and restrict \(\mu_j\) to \(Z_{2L}\). For large \(j\), define on \(D_L\)
\[
r_j^\pm(\xi)=\left(a_j^{-1}\max_{\substack{X\in\Sigma_j\cap Z_L\\ \pi(X)=\xi}}(\pm h(X))-\delta\right)_+,
\]
with \(r_j^\pm=0\) outside \(\pi(\Sigma_j\cap Z_L)\). By Proposition~\ref{cat:prop:compat-high-graphs} and \(m_j/a_j\to0\), \(r_j^\pm\) is locally Lipschitz and \(\|a_j\nabla r_j^\pm\|_{L^\infty(K)}\to0\) for every \(K\Subset D_L\).

Let \(r_j\) be either \(r_j^+\) or \(r_j^-\), and put \(b_j=(1+a_j^2|\nabla r_j|^2)^{-1/2}\), so that \(b_j\to1\) locally uniformly. Since \(r_j\ge0\) and \(\operatorname{div}(b_j\nabla r_j)=0\) on \(\{r_j>0\}\),
\begin{equation}\label{cat:eq:compat-roof-inequality}
\operatorname{div}(b_j\nabla r_j)\ge0 \quad \text{in }D_L.
\end{equation}

For \(K\Subset D_L\), the area formula and Lemma~\ref{cat:lem:height-tilt} give
\[
\int_K|\nabla r_j|^2\le\frac{2}{a_j^2} \int_{\pi^{-1}(K)}|\pi N|^2 \mathrm d\mu_j\le C_K.
\]
Outside \(E_j\) from Lemma~\ref{cat:lem:compat-weak}, \(r_j^\pm=(\pm u_{j,\pm}/a_j-\delta)_+\). Since \(|E_j|\to0\) and both sides are uniformly bounded on compact subsets of \(D_L\), \(r_j^\pm\to(v-\delta)_+\) strongly in \(L^2_{\mathrm{loc}}\) and weakly in \(W^{1,2}_{\mathrm{loc}}\). Lemma~\ref{cat:lem:compat-delort} with \(A_j=b_jI\) gives
\begin{equation}\label{cat:eq:compat-roof-stress}
(\nabla r_j^\pm\otimes \nabla r_j^\pm)^0 \longrightarrow \left(\nabla (v-\delta)_+\otimes \nabla (v-\delta)_+\right)^0
\end{equation}
in distributions, and \((v-\delta)_+\) is subharmonic.

Since \(\Delta_{\Sigma_j}h=0\), the proof of \eqref{cat:eq:compat-band-energy} applies directly to \(h\). Thus, for \(K\Subset D_L\),
\begin{equation}\label{cat:eq:compat-band-tilt}
a_j^{-2}\int_{\pi^{-1}(K)\cap\{|h|\le\delta a_j\}} |\pi N|^2 \mathrm d\mu_j\le C_K\delta.
\end{equation}
The area formula gives
\[
\frac{\pi_\#\left((\pi N)\otimes(\pi N) \mathbf1_{\{|h|>\delta a_j\}}\mu_j\right)}{a_j^2} =\sum_\pm\frac{\nabla r_j^\pm\otimes \nabla r_j^\pm} {\sqrt{1+a_j^2|\nabla r_j^\pm|^2}} \mathcal L^2.
\]
Recall that \(a_j^{-2}\pi_\#((\pi N)\otimes(\pi N)\mu_j)\rightharpoonup\mathcal S\). Taking trace-free parts and letting first \(j\to\infty\) and then \(\delta\downarrow0\), we obtain from \eqref{cat:eq:compat-roof-stress} and \eqref{cat:eq:compat-band-tilt}
\[
\mathcal S^0=2(\nabla v\otimes \nabla v)^0 \mathcal L^2, \qquad \Delta v\ge0.
\]
Since \(\operatorname{div}\mathcal S^0=0\) by \eqref{cat:eq:compat-limit-moments}, \eqref{cat:eq:compat-gap-stress} follows.
\end{proof}

\begin{theorem}
\label{cat:thm:compat-arbitrary-limits}
Every limit \(v\) has the form
\[
v(\xi)=c+|p\cdot \xi-b|, \qquad c\ge0, \quad p\in\mathbb R^2, \quad b\in\mathbb R,
\]
and satisfies \(2\int_{D_8}v^2\mathrm d A=1\).
\end{theorem}
\begin{proof}
By \eqref{cat:eq:compat-gap-stress}, \(q=(v_{x}-iv_{z})^2\) is holomorphic on \(\mathbb C\). For \(R\ge1\) and \(w\in D_{R/2}\), \eqref{cat:eq:compat-gap-growth} implies
\[
|q(w)|\le CR^{-2}\int_{D_R}|\nabla v|^2\mathrm d A \le CR^{-4}\int_{D_{2R}}v^2\mathrm d A\le C.
\]
Hence \(q\) is constant. If \(q=0\), \(v\) is constant. If \(q\ne0\), we may assume that \(q=k^2\) for some \(k>0\). Then \(v_z=0\) and \(v_x\in\{-k,k\}\) almost everywhere. Since \(v\ge0\) and \(v_{xx}=\Delta v\ge0\), we have \(v=c+k|x-x_0|\) with \(c\ge0\).
\end{proof}

If \(v\) is constant, then
\begin{equation}\label{cat:eq:constant-profile}
v=b_*=(2|D_8|)^{-1/2}.
\end{equation}

\begin{theorem}
\label{cat:thm:regimes}
As \(R\to\infty\), \(a_{2R}/a_R\to\ell\in\{1/2,1\}\). Let \(\varepsilon>0\), and let \(v\) be any limit in Theorem~\ref{cat:thm:compat-arbitrary-limits}.
\begin{enumerate}
\item If \(\ell=1/2\), then \(\mathsf H_R\le C_\varepsilon R^\varepsilon\) and \(v=b_*\).
\item If \(\ell=1\), then \(a_R\ge c_\varepsilon R^{-\varepsilon}\) and
\(v=|p\cdot \xi|\).
\end{enumerate}
\end{theorem}
\begin{proof}
For a limit \(v=c+|p\cdot \xi-b|\) and \(s>0\), put
\begin{equation}\label{cat:eq:model-scale-data}
T_s(v)=\left(\frac2{s^4}\int_{D_{8s}}v^2\mathrm d A\right)^{1/2}.
\end{equation}
Set \(\rho_R=a_{2R}/a_R\). By \eqref{cat:eq:scale-model-compatibility},
\[
\rho_{R_j}\to T_2(v), \qquad \rho_{2R_j}\to\frac{T_4(v)}{T_2(v)}
\]
for every subsequence with \(u_{j,\pm}/a_j\to\pm v\) in \(L^2_{\mathrm{loc}}\).
By the symmetry of \(D_8\),
\begin{equation}\label{cat:eq:model-parameter-norm}
\begin{aligned}[b]
&J_v(s):=s^2T_s(v)^2=2\int_{D_8}(c+|sp\cdot\xi-b|)^2\,\mathrm d A\\
={}&2\int_{D_8}\left(c^2+b^2+s^2(p\cdot \xi)^2 +2c\max\{|b|,s|p\cdot \xi|\}\right)\mathrm d A.
\end{aligned}
\end{equation}
Thus \(s^{-1}\le T_s(v)\le1\) for \(s>1\). We have \(T_s(v)=s^{-1}\) exactly when \(v=b_*\), and \(T_s(v)=1\) exactly when \(v=|p\cdot\xi|\). Hence \(\operatorname{dist}(\rho_R,[1/2,1])\to0\).

For \(A,B\ge0\),
\[
\max\{A,2B\}^2\le\max\{A,B\}\max\{A,4B\}.
\]
By the Cauchy--Schwarz inequality for the sum and integral in \eqref{cat:eq:model-parameter-norm}, \(J_v(2)^2\le J_v(1)J_v(4)=J_v(4)\). If \(p\ne0\) and \(c^2+b^2>0\), the constant and \(s^2\) terms in \eqref{cat:eq:model-parameter-norm} are both nonzero, so \(J_v(2)^2<J_v(4)\). Hence
\begin{equation}\label{cat:eq:model-doubling-increase}
T_4(v)\ge T_2(v)^2.
\end{equation}
Equality holds precisely when \(p=0\) or \(c=b=0\), equivalently \(T_2(v)\in\{1/2,1\}\).

Fix \(1/2<t<1\). Suppose that \(\rho_{R_j}\ge t\ge\rho_{2R_j}\) for a sequence \(R_j\to\infty\). After passing to a subsequence, \(u_{j,\pm}/a_j\to\pm v\) in \(L^2_{\mathrm{loc}}\). By \eqref{cat:eq:scale-model-compatibility},
\[
T_2(v)\ge t\ge\frac{T_4(v)}{T_2(v)}\ge T_2(v),
\]
contrary to the equality case in \eqref{cat:eq:model-doubling-increase}. Thus there is \(R_t\) such that
\begin{equation}\label{cat:eq:ratio-interior-increase}
\rho_{2R}>t\qquad\text{whenever }R\ge R_t\text{ and }\rho_R\ge t.
\end{equation}

Fix \(R_0\) and put \(r_n=\rho_{2^nR_0}\). For every \(1/2<t<\limsup r_n\), \eqref{cat:eq:ratio-interior-increase} gives \(r_n\ge t\) for all sufficiently large \(n\). Thus \(r_n\to\ell\). If \(u_{j,\pm}/a_j\to\pm v\) for \(R_j=2^{n_j}R_0\), then \(T_2(v)=\ell\) and \(T_4(v)=\ell^2\) by \eqref{cat:eq:scale-model-compatibility}. By the equality case in \eqref{cat:eq:model-doubling-increase}, either \(v=b_*\) and \(\ell=1/2\), or \(v=|p\cdot\xi|\) and \(\ell=1\). In either case \(T_{2s}(v)=\ell T_s(v)\) for every \(s>0\). By \eqref{cat:eq:scale-model-compatibility},
\[
\rho_{sR_j}=\frac{a_{2sR_j}/a_{R_j}}{a_{sR_j}/a_{R_j}} \longrightarrow\frac{T_{2s}(v)}{T_s(v)}=\ell
\]
uniformly for \(1\le s\le2\). By compactness,
\[
\sup_{1\le s\le2}\left|\rho_{s2^nR_0}-\ell\right|\longrightarrow0.
\]
Thus \(\rho_R\to\ell\). If \(\ell=1/2\), every limit \(v\) equals \(b_*\). Iterating \(a_{2R}\le2^{-1+\varepsilon}a_R\) for large \(R\) gives \(\mathsf H_R=Ra_R\le C_\varepsilon R^\varepsilon\).

If \(\ell=1\), every limit \(v\) has the form \(|p\cdot\xi|\). Iterating \(a_{2R}\ge2^{-\varepsilon}a_R\) for large \(R\) gives \(a_R\ge c_\varepsilon R^{-\varepsilon}\).
\end{proof}

\begin{proposition}\label{cat:prop:strong-height-limits}
The functions \(u_{j,\pm}\) can be chosen so that
\begin{equation}\label{cat:fold:eq:fold}
u_{j,\pm}/a_j\longrightarrow\pm v \quad \text{ strongly in }W^{1,2}_{\mathrm{loc}}.
\end{equation}
Moreover,
\begin{equation}\label{cat:fold:eq:moment-input}
\frac{\kappa_j}{a_j^2}\rightharpoonup2|\nabla v|^2 \mathcal L^2.
\end{equation}
Every varifold limit \(V\) of \(|a_j^{-1}(\Sigma_j-\xi_j)|\) with \((\xi_j)\) bounded has the form
\begin{equation}\label{cat:fold:eq:height-scale-planes}
V=|\{h=t\}|+|\{h=-t\}|
\end{equation}
for some \(t\ge0\).
\end{proposition}
\begin{proof}
By Theorem~\ref{cat:thm:regimes}, \(v^2\) is constant or equals \((p\cdot\xi)^2\). Thus \(\Delta(v^2)=2|\nabla v|^2\), and \eqref{cat:fold:eq:moment-input} follows from \eqref{cat:eq:compat-limit-moments}.

Fix \(L>0\) and restrict \(\Sigma_j\) to \(Z_{4L}\). By Propositions~\ref{cat:prop:compat-pressure} and \ref{cat:prop:lipschitz-approximation} with \(q=4\), choose \(u_{j,-}\le u_{j,+}\) and \(E_j\subset D_L\) such that \(|\Sigma_j|=\sum_\pm|\operatorname{graph}u_{j,\pm}|\) over \(D_L\setminus E_j\) and
\begin{equation}\label{cat:fold:eq:approximation-error}
|E_j|+\mu_j(\pi^{-1}(E_j))=o(a_j^2).
\end{equation}

For \(0\le\chi\in C_c^\infty(D_L)\), by \eqref{cat:eq:compat-fourth-tilt} and \eqref{cat:fold:eq:approximation-error},
\[
\sum_\pm\int\chi\left|\frac{\nabla u_{j,\pm}}{a_j}\right|^2\mathrm d A =a_j^{-2}\int(\chi\circ\pi)|\pi N|^2\mathrm d\mu_j+o(1) \longrightarrow2\int\chi|\nabla v|^2\mathrm d A.
\]
Since \(u_{j,\pm}/a_j\rightharpoonup\pm v\) in \(W^{1,2}_{\mathrm{loc}}\) by Lemma~\ref{cat:lem:compat-weak}, \eqref{cat:fold:eq:fold} follows.

By \eqref{cat:eq:compat-height-scale-limit}, every varifold limit of \(|a_j^{-1}(\Sigma_j-\xi_j-m_j(\xi_j)e_h)|\) has the form \eqref{cat:fold:eq:height-scale-planes}. Since \(m_j(\xi_j)/a_j\to0\), the same holds for \(|a_j^{-1}(\Sigma_j-\xi_j)|\).
\end{proof}

\section{Constancy of the limiting function \texorpdfstring{\(v\)}{v}}\label{sec:nonconstant-heights}

\begin{theorem}
\label{cat:thm:folded-exclusion}
\(a_{2R}/a_R\to1/2\) as \(R\to\infty\).
\end{theorem}

We argue by contradiction. By Theorem~\ref{cat:thm:regimes}, we may assume throughout this section that \(a_{2R}/a_R\to1\).

\subsection{Approximation by two planes}\label{cat:sec:affine-approximation}

For a pair \(C\) of planes with unit normals \(\mathbf{n}_1,\mathbf{n}_2\), put
\[
\mathcal T_C(N)=\min_{i=1,2}\{1-(N\cdot \mathbf{n}_i)^2\}.
\]
For a plane \(T\) with unit normal \(N\), also write \(\mathcal T_C(T)=\mathcal T_C(N)\). For \(0\le k\le1\), put \(P_\pm=\{h=\pm kx\}\) and
\begin{equation}\label{cat:fold:eq:relative-tilt}
\mathcal T_k(N)=\mathcal T_{P_+\cup P_-}(N) =N_{z}^2+\frac{(|N_{x}|-k|N_h|)^2}{1+k^2}.
\end{equation}

Let \(V\) be an integral \(2\)-varifold satisfying \eqref{cat:fold:eq:relative-BDD-hypotheses}, and put \(\mu=\|V\|\). Fix a sufficiently small \(\tau>0\), independent of \(k\), and take \(K_\tau\), \(u_1\le u_2\), and \(\Gamma\) from Proposition~\ref{cat:prop:lipschitz-approximation} with \(\lambda=\tau\). Write \(L=C_0\tau^{1/4}\), so that \(\operatorname{Lip}u_i\le L\). Write \(\|V-\Gamma\|_{\mathrm{TV},U} =|V-\Gamma|(U\times\mathrm{Gr}(2,3))\), where \(|V-\Gamma|\) is the total variation measure and \(\mathrm{Gr}(2,3)\) denotes the two-dimensional subspaces of \(\mathbb R^3\).

\begin{proposition}\label{cat:fold:prop:relative-bad-set}
There are \(c_*>0\) and \(C_\tau<\infty\) such that, if both planes in \(C\) make an angle at most \(\delta\) with \(\{h=0\}\) and \(\sin^2\delta\le c_*\tau\), then
\begin{equation}\label{cat:fold:eq:relative-bad-set}
|D_1\setminus K_\tau| +\mu(\pi^{-1}(D_1\setminus K_\tau)) +\|V-\Gamma\|_{\mathrm{TV},\pi^{-1}(D_1)} \le C_\tau\int_{\pi^{-1}(D_4)}\mathcal T_C^2\mathrm d\mu.
\end{equation}
For \(C=P_+\cup P_-\), also
\begin{equation}\label{cat:fold:eq:relative-reconstruction}
\sum_{i=1}^2\int_{D_1}\min_{\sigma=\pm1}|\nabla u_i-\sigma ke_x|^2\mathrm d A \le C_\tau\int_{\pi^{-1}(D_4)}\mathcal T_k\mathrm d\mu.
\end{equation}
\end{proposition}
\begin{proof}
For each \(\xi\in D_1\setminus K_\tau\), choose the disk \(D_r(y)\) from Proposition~\ref{cat:prop:lipschitz-approximation}. By \eqref{cat:fold:eq:relative-cover-mass},
\[
\tau\pi r^2 <2\int_{\pi^{-1}(D_r(y))}\mathcal T_0\mathrm d\mu \le4\int_{\pi^{-1}(D_r(y))}\mathcal T_C\mathrm d\mu+4C_0\sin^2\delta r^2.
\]
For sufficiently small \(c_*\), by the Cauchy--Schwarz inequality and \eqref{cat:fold:eq:relative-cover-mass},
\[
r^2\le C_\tau r^{-2}\left(\int_{\pi^{-1}(D_r(y))}\mathcal T_C\mathrm d\mu\right)^2 \le C_\tau\int_{\pi^{-1}(D_r(y))}\mathcal T_C^2\mathrm d\mu.
\]
By Vitali's covering lemma, choose pairwise disjoint disks \(D_{r_i}(y_i)\) such that \(D_1\setminus K_\tau\subset\bigcup_i D_{5r_i}(y_i)\). Since \(\Gamma=V\) over \(K_\tau\) and \(\|\Gamma\|(\pi^{-1}(D_1\setminus K_\tau)) \le2\sqrt{1+L^2}|D_1\setminus K_\tau|\),
\begin{multline*}
|D_1\setminus K_\tau|+\mu(\pi^{-1}(D_1\setminus K_\tau)) +\|V-\Gamma\|_{\mathrm{TV},\pi^{-1}(D_1)} \le C\sum_i r_i^2 \\
\le C_\tau\int_{\pi^{-1}(\bigcup_iD_{r_i}(y_i))}\mathcal T_C^2\mathrm d\mu.
\end{multline*}

On a graph \(h=u(\xi)\) with \(p=(p_x,p_z)=\nabla u\), \eqref{cat:fold:eq:relative-tilt} becomes
\begin{equation}\label{cat:fold:eq:relative-graph}
\mathcal T_k(N)\mathrm d\mu =\frac{(|p_x|-k)^2+(1+k^2)p_z^2}{(1+k^2)\sqrt{1+|p|^2}}\mathrm d A.
\end{equation}
For \(|p|\le L\),
\[
\frac{\min_{\sigma=\pm1}|p-\sigma ke_x|^2}{2\sqrt{1+L^2}}\mathrm d A \le\mathcal T_k(N)\mathrm d\mu \le\min_{\sigma=\pm1}|p-\sigma ke_x|^2\mathrm d A.
\]
For \(C=P_+\cup P_-\), apply \eqref{cat:fold:eq:relative-graph} on \(K_\tau\) and \(|\nabla u_i|\le L\) on \(D_1\setminus K_\tau\). Since \(\mathcal T_k^2\le \mathcal T_k\), \eqref{cat:fold:eq:relative-bad-set} gives \eqref{cat:fold:eq:relative-reconstruction}.
\end{proof}

Let \(u,v\) be \(L_0\)-Lipschitz functions on a domain \(D\subset\mathbb R^2\). For an orthogonal map \(Q\) sufficiently close to \(I\), write
\[
Q(\operatorname{graph}u)=\operatorname{graph}u', \qquad Q(\operatorname{graph}v)=\operatorname{graph}v'.
\]
Writing \(\pi(QX)=A\xi+bh\) and \(h(QX)=c\cdot\xi+dh\), we have
\[
\begin{aligned}
T_u(\xi)&:=A\xi+bu(\xi),& \qquad u'(T_u(\xi))&=c\cdot \xi+du(\xi), \\
T_v(\xi)&:=A\xi+bv(\xi),& \qquad v'(T_v(\xi))&=c\cdot \xi+dv(\xi).
\end{aligned}
\]
Then \(T_u,T_v\) are bi-Lipschitz, and \(u',v'\) are Lipschitz. The Lipschitz and bi-Lipschitz constants depend only on \(L_0\). For \(T_u(\xi)\in T_u(D)\cap T_v(D)\),
\begin{equation}\label{cat:fold:eq:graph-residual-transport}
\begin{aligned}
(u'-v')(T_u(\xi)) &=d(u-v)(\xi)+v'(T_v(\xi))-v'(T_u(\xi)) \\
&=(u-v)(\xi)+O_{L_0}\left(\|Q-I\| |u-v|(\xi)\right).
\end{aligned}
\end{equation}
For \(K\Subset T_u(D)\cap T_v(D)\),
\begin{equation}\label{cat:fold:eq:graph-height-comparison}
\int_K|u'-v'|^2\le C_{L_0} \int_{T_u^{-1}(K)}|u-v|^2.
\end{equation}

For Lipschitz functions \(u_1\le u_2\) and \(\widetilde u_1\le\widetilde u_2\) from \(D\) to \([-H,H]\), put \(\Gamma=\sum_i|\operatorname{graph}u_i|\) and \(\widetilde\Gamma=\sum_i|\operatorname{graph}\widetilde u_i|\). The area formula gives
\begin{equation}\label{cat:fold:eq:reconstruction-height-comparison}
\sum_i\int_D|u_i-\widetilde u_i|^2 \le CH^2\|\Gamma-\widetilde\Gamma\|_{\mathrm{TV},\pi^{-1}(D)}.
\end{equation}

Let \(\omega\) be as in Section~\ref{cat:sec:height-and-planes}, and write
\[
\|v\|_{L^2_\omega(D_1)}^2=\int_{D_1}\omega(|\xi|)|v(\xi)|^2\mathrm d A.
\]
Let \(\mathcal A\) be the affine functions on \(\mathbb R^2\). For \(m,\ell\in\mathcal A\), put
\[
\mathcal F(m,\ell)=(m-|\ell|,m+|\ell|).
\]
\begin{lemma}
\label{cat:fold:lem:height-fit-expansion}
For all \((m,\ell),(\widetilde m,\widetilde\ell)\in\mathcal A^2\),
\begin{equation}\label{cat:fold:eq:fit-stability}
\begin{aligned}
&\|m-\widetilde m\|_{L^2_\omega(D_1)}^2+\min_\pm\|\ell\pm\widetilde\ell\|_{L^2_\omega(D_1)}^2 \\
\le{}& C\|\mathcal F(m,\ell)-\mathcal F(\widetilde m,\widetilde\ell)\|_{L^2_\omega(D_1)}^2,
\end{aligned}
\end{equation}
where \(C\ge1\) depends only on \(\omega\).
\end{lemma}
\begin{proof}
By the definition of \(\mathcal F\),
\begin{align*}
&\|\mathcal F(m,\ell)-\mathcal F(\widetilde m,\widetilde\ell)\|_{L^2_\omega(D_1)}^2 \\
={}&2\|m-\widetilde m\|_{L^2_\omega(D_1)}^2 +2\int_{D_1}\omega(|\xi|)\min\{(\ell-\widetilde\ell)^2, (\ell+\widetilde\ell)^2\}\mathrm d A.
\end{align*}
The unit sphere of \(\mathcal A\) is compact, and \(|\{p=0\}|=0\) for \(p\in\mathcal A\setminus\{0\}\). Thus, by homogeneity, for \(p,q\in\mathcal A\),
\[
\int_{D_1}\omega(|\xi|)\min\{p^2,q^2\}\mathrm d A \ge c\min\{\|p\|_{L^2_\omega(D_1)}^2,\|q\|_{L^2_\omega(D_1)}^2\}.
\]
Taking \(p=\ell-\widetilde\ell\) and \(q=\ell+\widetilde\ell\) proves \eqref{cat:fold:eq:fit-stability}.
\end{proof}

Fix \(\tau\) as in Proposition~\ref{cat:fold:prop:relative-bad-set}. Choose coordinates with \(P_R=\{h=0\}\) and restrict \(R^{-1}M\) to \(\pi^{-1}(D_4)\cap\{|h|<1\}\). For large \(R\), \(|h|\le C a_R<1\), so this varifold is stationary in \(\pi^{-1}(D_4)\) and satisfies \eqref{cat:fold:eq:relative-BDD-hypotheses}. Proposition~\ref{cat:prop:lipschitz-approximation} with \(\lambda=\tau\) gives \(U_{R,1}\le U_{R,2}\) on \(D_1\), with \(|U_{R,i}|\le C a_R\). Write \(U_R=(U_{R,1},U_{R,2})\) and define
\begin{equation}\label{cat:fold:eq:buffered-height-fit-definition}
E_h(R)=\min_{m,\ell\in\mathcal A}\|U_R-\mathcal F(m,\ell)\|_{L^2_\omega(D_1)}^2.
\end{equation}
Choose a minimizer \((m_R,\ell_R)\), put \(f_R=\mathcal F(m_R,\ell_R)\), and let \(C_R\) be the two planes \(h=R(m_R\pm\ell_R)(\xi/R)\). Set
\begin{equation}\label{cat:fold:eq:geometric-buffer-error}
E_C(R)=R^{-2}\int_{M\cap B_{8R}}\mathcal T_{C_R}\mathrm d\mu, \qquad F(R)=e(R)+E_h(R)+E_C(R).
\end{equation}

Put \(q(R)=R^{-4}\int_{M\cap B_R}(X\cdot N)^2\mathrm d\mu\). Fix \(\Lambda=256\).

\begin{proposition}
\label{cat:fold:prop:geometric-buffer-improvement}
For every \(\varepsilon>0\), there are \(C_\varepsilon,R_\varepsilon<\infty\) such that, for \(R\ge R_\varepsilon\),
\begin{equation}\label{cat:fold:eq:geometric-buffer-improvement}
F(R)\le\varepsilon F(\Lambda R)+C_\varepsilon q(\Lambda R).
\end{equation}
\end{proposition}

We first compare \(F(R)\) with \(a_R^2\).

\begin{lemma}\label{cat:fold:lem:first-scale-relative-small}
\begin{equation}\label{cat:fold:eq:first-scale-relative-small}
F(R)=o(a_R^2).
\end{equation}
\end{lemma}
\begin{proof}
Given \(R_j\to\infty\), pass to a subsequence as in Proposition~\ref{cat:prop:strong-height-limits} and write \(R=R_j\). By \eqref{cat:fold:eq:relative-bad-set} and \eqref{cat:fold:eq:approximation-error}, \(U_R=(u_{j,-},u_{j,+})\) on \(D_1\) outside a set of area \(O(a_R^2)\), so
\[
U_R/a_R\longrightarrow(-v,v)\quad \text{ in }L^2_\omega(D_1), \qquad v=|p\cdot \xi|.
\]
By \eqref{cat:fold:eq:buffered-height-fit-definition},
\[
E_h(R)\le\|U_R-\mathcal F(0,a_Rp\cdot \xi)\|_{L^2_\omega(D_1)}^2=o(a_R^2).
\]
By \eqref{cat:fold:eq:fit-stability}, replacing \(\ell_R\) by \(-\ell_R\) if necessary, we have
\begin{equation}\label{cat:fold:eq:buffered-fit-first-scale}
m_R/a_R\to0, \qquad \ell_R/a_R\to p\cdot \xi \quad \text{ in }L^2_\omega(D_1).
\end{equation}

By \eqref{cat:fold:eq:fold}, \eqref{cat:fold:eq:approximation-error}, and \eqref{cat:fold:eq:buffered-fit-first-scale},
\[
E_C(R)\le C\sum_{\sigma=\pm1}\int_{D_8} \min_\pm|\nabla u_{j,\sigma}-\nabla(m_R\pm\ell_R)|^2\mathrm d A +o(a_R^2)=o(a_R^2).
\]
Since \(r^{-4}\int_{D_r}(p\cdot \xi)^2\mathrm d A\) is constant, \eqref{cat:eq:compat-area-cap-limit} gives for almost every \(r\in(0,1)\)
\[
0\le\frac{\pi r^4}{2}e(R) \le2\pi\int_0^r s^3e(Rs) \mathrm ds =\pi r^4-\int_{\Sigma_j}(r^2-|X|^2)_+\mathrm d\mu_j =o(a_R^2).\qedhere
\]
\end{proof}

\subsection{Affine limits}\label{cat:sec:affine-limits}

Throughout Sections~\ref{cat:sec:affine-limits} and \ref{cat:sec:height-tilt-convergence}, let \(s_j\to\infty\) and \(\beta_j>0\) satisfy
\begin{equation}\label{cat:fold:eq:secondary-sequence-assumptions}
\beta_j=o(a_{s_j}), \qquad F(s_j)\le C\beta_j^2, \qquad q(s_j)=o(\beta_j^2).
\end{equation}
Pass to a subsequence as in Proposition~\ref{cat:prop:strong-height-limits}.

Take orthonormal coordinates with \(P_{s_j}=\{h=0\}\) and \(\nabla\ell_{s_j}\) in the positive \(x\)-direction. By \eqref{cat:fold:eq:buffered-fit-first-scale}, a rotation \(Q_j\) fixing \(0\) with \(\|Q_j-I\|=o(a_{s_j})\) makes a bisector plane of \(s_j^{-1}C_{s_j}\) parallel to \(\{h=0\}\) and the intersection of the two planes parallel to the \(z\)-axis. Then \(Q_j(s_j^{-1}C_{s_j})\) is the union of the graphs of
\[
f_{j,1}=m_j-|\ell_j|, \qquad f_{j,2}=m_j+|\ell_j|, \qquad \ell_j(\xi)=k_j(x-d_j),
\]
where \(m_j\in\mathbb R\), \(C^{-1}a_{s_j}\le k_j\le Ca_{s_j}\), \(d_j\to0\), and \(|m_j|\le Ck_j\). Put
\[
\Sigma_j=Q_j\left(s_j^{-1}M\cap\pi^{-1}(D_8)\cap\{|h|<2\}\right), \qquad V_j=|\Sigma_j|, \qquad \mu_j=\|V_j\|.
\]
Then \(V_j\) satisfies \eqref{cat:fold:eq:relative-BDD-hypotheses}, \(|h|\le Ck_j\) on \(\Sigma_j\), and
\begin{equation}\label{cat:fold:eq:relative-secondary-scale}
k_j\longrightarrow0, \qquad \beta_j/k_j\longrightarrow0, \qquad \int_{\pi^{-1}(D_4)}\mathcal T_{k_j}\mathrm d\mu_j \le CE_C(s_j)\le C\beta_j^2.
\end{equation}
Since \(Q_j\) fixes \(0\), for every \(\Omega\Subset D_1\) and all large \(j\),
\begin{equation}\label{cat:fold:eq:relative-radial-small}
\int_{\pi^{-1}(\Omega)}(X\cdot N)^2\mathrm d\mu_j \le\int_{\Sigma_j\cap B_1}(X\cdot N)^2\mathrm d\mu_j =q(s_j)=o(\beta_j^2).
\end{equation}

Proposition~\ref{cat:prop:lipschitz-approximation} applied to \(V_j\) with \(\lambda=\tau\) gives \(K_j\) and \(u_{j,1}\le u_{j,2}\), with \(|u_{j,i}|\le Ck_j\) on \(D_1\). Put \(\Gamma_j=\sum_{i=1}^2|\operatorname{graph}u_{j,i}|\). Write \(Q_j(\operatorname{graph}U_{s_j,i})=\operatorname{graph}\bar u_{j,i}\), with \(\bar u_{j,1}\le\bar u_{j,2}\). For each \(0<r<1\), \(|\bar u_{j,i}|\le Ck_j\) on \(D_r\). Equations \eqref{cat:fold:eq:relative-bad-set} and \eqref{cat:fold:eq:reconstruction-height-comparison} give
\begin{equation}\label{cat:fold:eq:rotated-graph-comparison}
\sum_i\int_{D_r}|u_{j,i}-\bar u_{j,i}|^2\mathrm d A \le C_ra_{s_j}^2E_C(s_j)=o(\beta_j^2).
\end{equation}
Since \(\|U_{s_j}-f_{s_j}\|_{L^2_\omega(D_1)}^2=E_h(s_j)\le C\beta_j^2\) and \(\omega(|\xi|)>0\) on \(D_1\), \eqref{cat:fold:eq:graph-height-comparison} gives, for \(0<r<1\),
\begin{equation}\label{cat:fold:eq:relative-height-fit}
\sum_{i=1}^2\|u_{j,i}-f_{j,i}\|_{L^2(D_r)}^2\le C_r\beta_j^2.
\end{equation}

Put \(\sigma_{j,i}=1\) where \(\partial_xu_{j,i}\ge0\) and \(\sigma_{j,i}=-1\) elsewhere, and set
\begin{equation}\label{cat:fold:eq:relative-A}
A_j=(A_{j,x},A_{j,z})=\frac1{\beta_j}\sum_{i=1}^2 (\sigma_{j,i}\nabla u_{j,i}-k_je_x).
\end{equation}
By \eqref{cat:fold:eq:relative-reconstruction}, \(A_j\) is bounded in \(L^2(D_1;\mathbb R^2)\).

Put \(b_j=(u_{j,1}+u_{j,2})/2-m_j\).
\begin{lemma}\label{cat:fold:lem:linear-limits}
If \(A_j\rightharpoonup A=(A_x,A_z)\) in \(L^2(D_1)\), then \(A_x-iA_z\) is holomorphic. The sequence \(b_j/\beta_j\) is bounded in \(W^{1,2}_{\mathrm{loc}}(D_1)\), and every weak limit is harmonic.
\end{lemma}
\begin{proof}
Recall from \eqref{eq:graph-fluxes} that
\[
G(p)=\frac{p}{\sqrt{1+|p|^2}}, \qquad B(p)=\sqrt{1+|p|^2}I-\frac{p\otimes p}{\sqrt{1+|p|^2}}.
\]
For \(D\Subset D_1\), \(Y\in C_c^1(D;\mathbb R^2)\), and \(\varphi\in C_c^1(D)\), applying \(\delta\Gamma_j=\delta(\Gamma_j-V_j)\) to \(Y\circ\pi\) and \((\varphi\circ\pi)e_h\) gives
\begin{equation}\label{cat:fold:eq:relative-stress-error}
\left|\int_D\sum_i\operatorname{tr}\left(B(\nabla u_{j,i}) dY\right)\mathrm d A\right| \le C\|V_j-\Gamma_j\|_{\mathrm{TV},\pi^{-1}(D)}\|dY\|_{L^\infty(D)}
\end{equation}
and
\begin{equation}\label{cat:fold:eq:relative-vertical-error}
\left|\int_D\{G(\nabla u_{j,1})+G(\nabla u_{j,2})\}\cdot \nabla \varphi\mathrm d A\right| \le\|V_j-\Gamma_j\|_{\mathrm{TV},\pi^{-1}(D)}\|\nabla \varphi\|_{L^\infty(D)}.
\end{equation}

Since \(B(-p)=B(p)\), \eqref{cat:fold:eq:relative-reconstruction} gives
\[
\left\|\frac{\sum_iB(\nabla u_{j,i})-2B(k_je_x)}{k_j\beta_j} +\begin{pmatrix}A_{j,x}&A_{j,z}\\A_{j,z}&-A_{j,x}\end{pmatrix} \right\|_{L^1(D_1)}\le C\left(k_j^2+\frac{\beta_j}{k_j}\right)\longrightarrow0.
\]
By \eqref{cat:fold:eq:relative-bad-set} and \eqref{cat:fold:eq:relative-secondary-scale},
\[
\frac{\|V_j-\Gamma_j\|_{\mathrm{TV},\pi^{-1}(D_1)}}{k_j\beta_j}\le C\frac{\beta_j}{k_j}\longrightarrow0.
\]
In distributions on \(D_1\), \eqref{cat:fold:eq:relative-stress-error} implies
\[
\operatorname{div}\begin{pmatrix}A_x&A_z\\A_z&-A_x\end{pmatrix}=0.
\]
For \(\Omega\Subset D_1\), choose \(r<1\) and \(\chi\in C_c^1(D_r)\) with \(\chi=1\) on \(\Omega\). In \eqref{cat:fold:eq:relative-vertical-error}, take \(\varphi=\chi^2b_j\), for which \(\|\nabla \varphi\|_{L^\infty(D_r)}\le C_\chi\). Since \(G\) is odd and \((1+L^2)^{-3/2}I\le dG\le I\) on \(\{|p|\le L\}\),
\begin{equation}\label{cat:fold:eq:relative-mean-gradient}
\|\nabla b_j\|_{L^2(\Omega)}^2 \le C_\Omega\left(\|b_j\|_{L^2(D_r)}^2 +\|V_j-\Gamma_j\|_{\mathrm{TV},\pi^{-1}(D_r)}\right) \le C_\Omega\beta_j^2.
\end{equation}
For \(|p|,|q|\le L\),
\[
|G(p)+G(q)-p-q|\le C_L(|p|+|q|)|p+q|.
\]
Thus \eqref{cat:fold:eq:relative-reconstruction}, \eqref{cat:fold:eq:relative-mean-gradient}, and \eqref{cat:fold:eq:relative-vertical-error} give for \(\varphi\in C_c^\infty(D_1)\)
\[
\left|\int \nabla (b_j/\beta_j)\cdot \nabla \varphi\mathrm d A\right| \le C_\varphi(k_j+\beta_j)\longrightarrow0.\qedhere
\]
\end{proof}

\begin{lemma}
\label{cat:fold:lem:relative-offfold-gradient}
For every \(\Omega\Subset D_1\setminus\{x=0\}\),
\[
\sum_{i=1}^2\int_\Omega|\nabla u_{j,i}-\nabla f_{j,i}|^2\mathrm d A \le C_\Omega\beta_j^2.
\]
\end{lemma}
\begin{proof}
Choose \(\chi\in C_c^1(D_1\setminus\{x=0\})\) with \(\chi=1\) on \(\Omega\). For some \(\delta>0\), \(|\ell_j|\ge2\delta k_j\) on \(\operatorname{spt}\chi\) for all large \(j\). In this proof we suppress \(j\). Put \(p_0=\nabla f_2=\pm ke_x\), \(\zeta=h-f_2(\xi)\) and \(a=\nabla \zeta=e_h-p_0\). Choose a nondecreasing smooth \(\eta\) with \(\eta=0\) on \((-\infty,\delta/2]\) and \(\eta=1\) on \([\delta,\infty)\). By \eqref{cat:fold:eq:relative-bad-set} and \eqref{cat:fold:eq:relative-height-fit},
\begin{equation}\label{cat:fold:eq:relative-transition-mass}
\int_{\pi^{-1}(\operatorname{spt}\chi)} \min_{i=1,2}|h-f_i(\xi)|^2\mathrm d\mu\le C_\chi\beta^2.
\end{equation}

For \(h\ge m\), \(\zeta^2=\min_i|h-f_i(\xi)|^2\). Where \(\eta'((h-m)/k)\ne0\), \(|\zeta|\le Ck\) and \(k^2\le\delta^{-2}\min_i|h-f_i(\xi)|^2\). Since \(\int\operatorname{div}_\Sigma(\chi^2\zeta\eta((h-m)/k)a)\mathrm d\mu=0\) and \(|\nabla_\Sigma\zeta\cdot e_h|\le C(\mathcal T_k+k^2)\), \eqref{cat:fold:eq:relative-transition-mass} gives
\[
\int\chi^2\eta|\nabla_\Sigma\zeta|^2\mathrm d\mu \le C_{\delta,\chi}\int_{\pi^{-1}(\operatorname{spt}\chi)} \left(\mathcal T_k+\min_i|h-f_i(\xi)|^2\right)\mathrm d\mu \le C_{\delta,\chi}\beta^2.
\]
On \(\Omega\cap K\cap\{u_2-m\ge\delta k\}\), \(\chi=\eta=1\) and \(|\nabla u_2-p_0|^2\le C|\nabla_\Sigma\zeta|^2\). Where \(u_2-m<\delta k\), \(|u_2-f_2|\ge\delta k\), so
\[
|\nabla u_2-p_0|^2 \le C\min_{\sigma=\pm1}|\nabla u_2-\sigma ke_x|^2+C_\delta|u_2-f_2|^2.
\]
Equations \eqref{cat:fold:eq:relative-bad-set}, \eqref{cat:fold:eq:relative-reconstruction}, and \eqref{cat:fold:eq:relative-height-fit} imply \(\int_\Omega|\nabla u_2-p_0|^2\mathrm d A\le C_\Omega\beta^2\). Replacing \((h,m,u_2,f_2)\) by \((-h,-m,-u_1,-f_1)\) gives the same bound for \(i=1\).
\end{proof}

Define
\[
\eta_j=\operatorname{sgn}(x-d_j) \frac{(u_{j,2}-u_{j,1})/2-|\ell_j|}{\beta_j}.
\]
\begin{proposition}
\label{cat:fold:prop:relative-radial-affine}
After passing to a subsequence,
\[
b_j/\beta_j \longrightarrow0 \text{ in }L^2_{\mathrm{loc}}(D_1),\quad \eta_j \longrightarrow0 \text{ in }L^2_{\mathrm{loc}}(D_1\setminus\{x=0\}),
\]
and
\begin{equation}\label{cat:fold:eq:relative-radial-affine}
m_j=o(\beta_j), \qquad k_jd_j=o(\beta_j).
\end{equation}
\end{proposition}
\begin{proof}
By Lemmas~\ref{cat:fold:lem:linear-limits} and \ref{cat:fold:lem:relative-offfold-gradient}, pass to a subsequence with \(b_j/\beta_j\to b\) in \(L^2_{\mathrm{loc}}(D_1)\), \(\eta_j\to\eta\) weakly in \(L^2_{\mathrm{loc}}(D_1)\) and strongly in \(L^2_{\mathrm{loc}}(D_1\setminus\{x=0\})\), and \(A_j\rightharpoonup A\) in \(L^2(D_1)\). Then \(b\) is harmonic and \(A_x-iA_z\) is holomorphic.

Where \(\sigma_{j,i}k_je_x\ne \nabla f_{j,i}\), \(2k_j|\nabla u_{j,i}|\le|\nabla u_{j,i}-\nabla f_{j,i}|^2\). For \(\Omega\Subset D_1\setminus\{x=0\}\),
\[
\|A_j-2\nabla \eta_j\|_{L^1(\Omega)} \le\frac1{k_j\beta_j}\sum_i\|\nabla u_{j,i}-\nabla f_{j,i}\|_{L^2(\Omega)}^2 \le C_\Omega\beta_j/k_j\longrightarrow0.
\]
Hence \(A=2\nabla \eta\) in \(D_1\setminus\{x=0\}\).

Put \(\mathcal E=\xi\cdot \nabla\). On a graph \(h=u(\xi)\), \(|u-\mathcal Eu|\mathrm d A=|X\cdot N|\mathrm d\mu\). By \eqref{cat:fold:eq:relative-bad-set} and \eqref{cat:fold:eq:relative-radial-small}, for \(\Omega\Subset D_1\),
\[
\sum_i\int_\Omega|u_{j,i}-\mathcal Eu_{j,i}|\mathrm d A \le C_\Omega\left(\int_{\pi^{-1}(\Omega)} (X\cdot N)^2\mathrm d\mu_j\right)^{1/2} +C_\Omega\beta_j^2=o(\beta_j).
\]
For \(\Omega\Subset D_1\setminus\{x=0\}\),
\[
\left\|\frac{(\mathcal E-1)b_j-m_j}{\beta_j}\right\|_{L^1(\Omega)} +\left\|(\mathcal E-1)\eta_j+\frac{k_jd_j}{\beta_j}\right\|_{L^1(\Omega)} \longrightarrow0.
\]
Passing to the limit, \(m_j/\beta_j\to-b(0)\) and \(k_jd_j/\beta_j\to c\) for some \(c\in\mathbb R\), with
\begin{equation}\label{cat:fold:eq:relative-mean-Euler}
(\mathcal E-1)b=-b(0), \qquad (\mathcal E-1)\eta=-c\quad \text{ in }D_1\setminus\{x=0\}.
\end{equation}
By \eqref{cat:fold:eq:relative-mean-Euler}, \(\nabla b\) and \(A\) are homogeneous of degree zero. They are continuous at \(0\), hence constant. Thus \(b\) is affine and \(\eta(\xi)=c+\tfrac12A(0)\cdot \xi\).

By \eqref{cat:fold:eq:rotated-graph-comparison} and \eqref{cat:fold:eq:graph-residual-transport},
\[
\frac{U_{s_j}-f_{s_j}}{\beta_j} \rightharpoonup \left(b-\operatorname{sgn}(x)\eta,b+\operatorname{sgn}(x)\eta\right) \quad \text{ in }L^2_\omega(D_1).
\]
Since \(b,\eta\) are affine and \((m_{s_j},\ell_{s_j})\) attains the minimum in \eqref{cat:fold:eq:buffered-height-fit-definition},
\[
\sum_{i=1}^2\int_{D_1}\omega(|\xi|)\frac{U_{s_j,i}-f_{s_j,i}}{\beta_j}
\bigl(b+(-1)^i\operatorname{sgn}(\ell_{s_j})\eta\bigr)\,\mathrm d A=0.
\]
Since \(\ell_{s_j}/k_j\to x\) in \(L^2_\omega(D_1)\), letting \(j\to\infty\) gives
\[
\int_{D_1}\omega(|\xi|)(b^2+\eta^2)\,\mathrm d A=0.
\]
Thus \(b=\eta=0\) and \(c=\eta(0)=0\).
\end{proof}

\subsection{Strong convergence of heights and tilt}\label{cat:sec:height-tilt-convergence}

Since \(C^{-1}a_{s_j}\le k_j\le Ca_{s_j}\) and \(\|Q_j-I\|=o(k_j)\), Proposition~\ref{cat:prop:strong-height-limits} implies that every varifold limit of \(|k_j^{-1}(\Sigma_j-\xi_j)|\), with \(\xi_j\) in a fixed compact subset of \(D_1\), has the form
\[
|\{h=t\}|+|\{h=-t\}|
\]
for some \(t\ge0\).

\begin{proposition}
\label{cat:fold:prop:strong-affine-approximation}
\label{cat:fold:prop:relative-strong-height}
For every \(K\Subset D_1\),
\begin{equation}\label{cat:fold:eq:strong-affine-approximation}
\sum_{i=1}^2\int_K|u_{j,i}-f_{j,i}|^2\mathrm d A +\int_{\pi^{-1}(K)}\mathcal T_{k_j}\mathrm d\mu_j=o(\beta_j^2).
\end{equation}
\end{proposition}

We first establish strong convergence away from \(\{x=0\}\). Fix a connected \(U\Subset D_1\setminus\{x=0\}\) and \(0<4\delta<\inf_U|x|\). For \(\xi\in\pi(\Sigma_j)\), define
\begin{equation}\label{cat:fold:eq:actual-roofs}
\begin{aligned}
h_{j,1}(\xi)&=\min\bigl\{-\delta k_j,\min_{X\in\Sigma_j,\,\pi(X)=\xi}h(X)\bigr\},\\
h_{j,2}(\xi)&=\max\bigl\{\delta k_j,\max_{X\in\Sigma_j,\,\pi(X)=\xi}h(X)\bigr\}.
\end{aligned}
\end{equation}
For \(\xi\notin\pi(\Sigma_j)\), set \(h_{j,1}(\xi)=-\delta k_j\) and \(h_{j,2}(\xi)=\delta k_j\). Let \(E_j\) be the set of \(\xi\) for which \(\Sigma_j\cap\pi^{-1}(\xi)\) does not consist of exactly two points with \(|N_h|>1/2\) at both.

\begin{lemma}
\label{cat:fold:prop:relative-offfold-compactness}
For every \(K\Subset U\),
\begin{equation}\label{cat:fold:eq:relative-roof-strong}
\sum_{i=1}^2\|h_{j,i}-f_{j,i}\|_{W^{1,2}(K)}^2 +\int_{\pi^{-1}(K)}\mathcal T_{k_j}\mathrm d\mu_j=o(\beta_j^2),
\end{equation}
and
\begin{equation}\label{cat:fold:eq:relative-offfold-band-little}
\int_{\pi^{-1}(K)\cap\{|h|<\delta k_j\}} (|\pi N|^2+k_j^2)\mathrm d\mu_j +|E_j\cap K| +\mu_j\left(\pi^{-1}(E_j\cap K)\right)=o(\beta_j^2).
\end{equation}
\end{lemma}
\begin{proof}
The functions \(h_{j,i}\) are locally Lipschitz, with \(\|\nabla h_{j,i}\|_{L^\infty(K)}\to0\) for \(K\Subset U\), by \eqref{cat:fold:eq:height-scale-planes} and Lemma~\ref{lem:varifold-limits}. For large \(j\), \(x-d_j\) has constant sign on \(U\), so \(f_{j,2}\) is affine there, with \(f_{j,2}-\delta k_j\ge ck_j>0\) and \(|\nabla f_{j,2}|=k_j\). Put \(w_j=(h_{j,2}-f_{j,2})/\beta_j\). On \(K_j\), \(h_{j,2}=\max\{u_{j,2},\delta k_j\}\). By \eqref{cat:fold:eq:relative-bad-set}, Lemma~\ref{cat:fold:lem:relative-offfold-gradient}, and Proposition~\ref{cat:fold:prop:relative-radial-affine}, for \(K\Subset U\),
\begin{equation}\label{cat:fold:eq:relative-to-roof-offfold}
\begin{aligned}
\int_Kw_j^2 &\le\frac1{\beta_j^2}\left(\int_K|u_{j,2}-f_{j,2}|^2 +Ck_j^2|K\setminus K_j|\right)=o(1), \\
\int_K|\nabla w_j|^2 &\le\frac C{\beta_j^2}\left(\|u_{j,2}-f_{j,2}\|_{W^{1,2}(K)}^2 +|K\setminus K_j|\right)=O(1).
\end{aligned}
\end{equation}
By \eqref{cat:fold:eq:relative-radial-affine},
\[
\xi\cdot\nabla w_j-w_j
=\frac{\xi\cdot\nabla h_{j,2}-h_{j,2}}{\beta_j}+o(1).
\]
On \(\{h_{j,2}=\delta k_j\}\), \(\nabla h_{j,2}=0\) and \(|\xi\cdot \nabla w_j|\le C|w_j|\) almost everywhere. On \(\{h_{j,2}>\delta k_j\}\), \(\xi+h_{j,2}(\xi)e_h\in\Sigma_j\), so the area formula and \eqref{cat:fold:eq:relative-radial-small} give
\[
\int_K|\xi\cdot \nabla w_j|^2 \le C_K\int_Kw_j^2 +\frac{C}{\beta_j^2}\int_{\pi^{-1}(K)}(X\cdot N)^2\mathrm d\mu_j +o(1)\longrightarrow0.
\]

As in \eqref{cat:eq:compat-roof-inequality}, \(\operatorname{div}G(\nabla h_{j,2})\ge0\). Since \(f_{j,2}\) is affine, \(w_j\) satisfies \(\operatorname{div}(B_j\nabla w_j)\ge0\) on \(U\), where
\[
B_j=\int_0^1dG(\nabla f_{j,2}+t\beta_j\nabla w_j)\mathrm d t\longrightarrow I \quad\text{locally uniformly on }U.
\]
By Lemma~\ref{cat:lem:compat-delort}, for every \(0\le\chi\in C_c^\infty(U)\),
\[
\int\chi|\xi|^2|\nabla w_j|^2 =2\int\chi(\xi\cdot \nabla w_j)^2 -2\int\chi \xi\cdot(\nabla w_j\otimes \nabla w_j)^0\xi\longrightarrow0.
\]
Similarly, \((h_{j,1}-f_{j,1})/\beta_j\to0\) in \(W^{1,2}_{\mathrm{loc}}(U)\).

To prove \eqref{cat:fold:eq:relative-offfold-band-little}, it suffices to take \(K=\overline{D_r(\xi_0)}\) and \(D=D_{4r}(\xi_0)\Subset U\). Since \(h_{j,1}/k_j\to-|x|\) and \(h_{j,2}/k_j\to|x|\) strongly in \(W^{1,2}(D)\), we can choose disks \(D_j=D_{r_j}(\xi_0)\) with \(2r<r_j<3r\) such that
\[
h_{j,1}<-\delta k_j<\delta k_j<h_{j,2}
\quad\text{on }\partial D_j.
\]

The measure \(\lambda_j=\operatorname{div}G(\nabla h_{j,2})\ge0\) is supported on \(\{h_{j,2}=\delta k_j\}\), where \(\beta_jw_j\le-ck_j\). Take \(\chi\in C_c^\infty(D)\) with \(\chi=1\) on \(D_{3r}(\xi_0)\). Since \(\|dG\|\le1\),
\[
\begin{aligned}
ck_j\lambda_j(D_j)&\le-\beta_j\int\chi^2w_j\mathrm d\lambda_j \\
&=\beta_j\int\{G(\nabla h_{j,2})-G(\nabla f_{j,2})\}\cdot \nabla (\chi^2w_j) \\
&\le C\beta_j^2\|w_j\|_{W^{1,2}(D)}^2=o(\beta_j^2).
\end{aligned}
\]
The surface \(\Sigma_j\cap\pi^{-1}(\overline D_j)\) is compact, and by \eqref{cat:fold:eq:height-scale-planes} its boundary consists of the graphs of \(h_{j,1}\) and \(h_{j,2}\) over \(\partial D_j\). Let \(\nu\) be the outward unit normal to \(\partial D_j\). Integrating \(\Delta_{\Sigma_j}h=0\) gives
\[
\int_{\partial D_j}G(\nabla h_{j,1})\cdot\nu\,\mathrm d s
=-\int_{\partial D_j}G(\nabla h_{j,2})\cdot\nu\,\mathrm d s.
\]
Multiplying \(\Delta_{\Sigma_j}h=0\) by \(\max\{-\delta k_j,\min\{h,\delta k_j\}\}\) and integrating by parts on \(\Sigma_j\cap\pi^{-1}(D_j)\) gives
\begin{equation}\label{cat:fold:eq:band-flux-identity}
\begin{aligned}
\int_{\pi^{-1}(D_j)\cap\{|h|<\delta k_j\}}|\pi N|^2\mathrm d\mu_j &=2\delta k_j\int_{\partial D_j}G(\nabla h_{j,2})\cdot\nu\mathrm d s \\
&=2\delta k_j\lambda_j(D_j)=o(\beta_j^2).
\end{aligned}
\end{equation}

By \eqref{cat:fold:eq:height-scale-planes} and Lemma~\ref{lem:varifold-limits},
\[
|k_j^{-1}(\Sigma_j-\xi_j)|\longrightarrow2|\{h=0\}| \qquad\text{for }\xi_j\in E_j\cap K.
\]
Indeed, if a limit were \(|\{h=t\}|+|\{h=-t\}|\) with \(t>0\), then \(\Sigma_j\cap\pi^{-1}(\xi_j)\) would consist of exactly two points, both with \(|N_h|>1/2\). This contradicts \(\xi_j\in E_j\).
By Lemma~\ref{lem:varifold-limits} and \(|h|\le Ck_j\), for all sufficiently large \(j\) and every \(\xi\in E_j\cap K\),
\[
\Sigma_j\cap\pi^{-1}(D_{4k_j}(\xi))\subset\{|h|<\delta k_j\}.
\]
Apply \eqref{cat:fold:eq:relative-bad-set} and Proposition~\ref{cat:prop:lipschitz-approximation} to \(k_j^{-1}(\Sigma_j-\xi)\) with \(P_+=P_-=\{h=0\}\). Since \(\mathcal T_0=|\pi N|^2\),
\[
|E_j\cap D_{k_j}(\xi)| +\mu_j\left(\pi^{-1}(E_j\cap D_{k_j}(\xi))\right) \le C \int_{\pi^{-1}(D_{4k_j}(\xi))}|\pi N|^2\mathrm d\mu_j.
\]
Summing over a maximal \(k_j\)-separated subset of \(E_j\cap K\) and using \eqref{cat:fold:eq:band-flux-identity} and the bounded overlap of the disks \(D_{4k_j}(\xi)\), we obtain
\[
|E_j\cap K|+\mu_j\left(\pi^{-1}(E_j\cap K)\right)=o(\beta_j^2).
\]
If \(\xi\in K\setminus E_j\) and \(\Sigma_j\cap\pi^{-1}(\xi)\) meets \(\{|h|<\delta k_j\}\), then \(h_{j,1}(\xi)=-\delta k_j\) or \(h_{j,2}(\xi)=\delta k_j\), so \(\sum_i|h_{j,i}(\xi)-f_{j,i}(\xi)|^2\ge ck_j^2\). By the area formula,
\[
k_j^2\mu_j\left(\pi^{-1}(K\setminus E_j) \cap\{|h|<\delta k_j\}\right) \le C\sum_{i=1}^2\int_K|h_{j,i}-f_{j,i}|^2 =o(\beta_j^2).
\]

Since \(\mathcal T_{k_j}\le C(|\pi N|^2+k_j^2)\), \eqref{cat:fold:eq:relative-graph} and \eqref{cat:fold:eq:relative-offfold-band-little} imply
\[
\begin{aligned}[b]
&\int_{\pi^{-1}(K)}\mathcal T_{k_j}\mathrm d\mu_j \\
\le{}& \sum_i\int_K|\nabla h_{j,i}-\nabla f_{j,i}|^2\mathrm d A +C\int_{\pi^{-1}(K)\cap\{|h|<\delta k_j\}} (|\pi N|^2+k_j^2)\mathrm d\mu_j \\
={}& o(\beta_j^2).
\end{aligned}\qedhere
\]
\end{proof}

\begin{lemma}
\label{cat:fold:prop:axial-defect-exclusion}
Put \(S=xN_x+hN_h\) and \(C_j=|\{h=k_jx\}|+|\{h=-k_jx\}|\). For every \(K\Subset D_1\),
\begin{equation}\label{cat:fold:eq:axial-and-transverse-support-small}
\sum_{i=1}^2\int_K|u_{j,i}-f_{j,i}|^2\mathrm d A +\int_{\pi^{-1}(K)}(N_z^2+S^2)\mathrm d\mu_j=o(\beta_j^2).
\end{equation}
For every \(0\le\varphi\in C_c^\infty(B_1)\), radial, constant near \(0\), and nonincreasing in \(|X|\),
\begin{equation}\label{cat:fold:eq:strong-radial-mass}
\int\varphi\mathrm d\mu_j-\int\varphi\mathrm d\|C_j\|=o(\beta_j^2).
\end{equation}
\end{lemma}
\begin{proof}
Write \(r=\sqrt{x^2+h^2}\). For \(\delta>0\), let \(\psi\in C_c^\infty(\pi^{-1}(D_1)\cap\{r>8\delta\})\) depend only on \(r\) and \(z\). Since \(|h|\le Ck_j\), \(\operatorname{spt}\psi\cap\Sigma_j\subset\{|x|>4\delta\}\) for large \(j\). On \(\{|h|>\delta k_j\}\setminus\pi^{-1}(E_j)\), \(\Sigma_j\) agrees with the graphs of \(h_{j,i}\). For \(|h|\le\delta k_j\), \(\psi(x,h,z)-\psi(x,\pm\delta k_j,z)=O_\psi(k_j^2)\). The area formula, \(1-|N_h|\le|\pi N|^2\), and \eqref{cat:fold:eq:relative-offfold-band-little} give
\[
\int \psi\mathrm d\mu_j =\sum_{i=1}^2\int \psi(x,h_{j,i}(x,z),z) \sqrt{1+|\nabla h_{j,i}|^2}\mathrm d x\mathrm d z+o(\beta_j^2).
\]

For \(i=1,2\), put \(g=h_{j,i}\) and \(p=(-1)^ik_j\). By \eqref{cat:fold:eq:relative-radial-affine} and \eqref{cat:fold:eq:relative-roof-strong}, \(g-px=o(\beta_j)\) in \(W^{1,2}_{\mathrm{loc}}(D_1\cap\{x>4\delta\})\). Since \(\psi_h(x,px,z)=p\psi_x(x,px,z)\),
\[
\begin{aligned}
&\int_{D_1\cap\{x>4\delta\}}\left[\psi(x,g,z)\sqrt{1+|\nabla g|^2}-\sqrt{1+p^2}\psi(x,px,z)\right]\mathrm d x\mathrm d z \\
={}& \frac p{\sqrt{1+p^2}}\int_{D_1\cap\{x>4\delta\}}\partial_x\left(\psi(x,px,z)(g-px)\right)\mathrm d x\mathrm d z +o(\beta_j^2)=o(\beta_j^2).
\end{aligned}
\]
Repeating the calculation on \(\{x<-4\delta\}\) with \(p=(-1)^{i+1}k_j\) and summing over \(i=1,2\), we obtain
\begin{equation}\label{cat:fold:eq:radial-scalar-comparison}
\int \psi \mathrm d(\mu_j-\|C_j\|)=o(\beta_j^2).
\end{equation}

Write \(\varphi(X)=F(r,z)=f(\sqrt{r^2+z^2})\), with \(\operatorname{spt}f\subset[0,R)\) for some \(R<1\). Choose \(\varepsilon>0\) with \(R^2+\varepsilon^2<1\) and \(\zeta\in C_c^\infty((-\varepsilon,\varepsilon))\) with \(\zeta=1\) on \([-\varepsilon/2,\varepsilon/2]\). Put
\[
c=\frac23(zF+r^2F_z), \qquad b=\zeta(r)\int_0^1(F-c_z)(tr,z)\mathrm d t, \qquad Y=(bx,bh,c).
\]
Then \(Y\) is compactly supported in \(\pi^{-1}(D_1)\), with \(b+rb_r+c_z=F\) near the axis, and
\begin{equation}\label{eq:axis-test}
\begin{aligned}
\operatorname{div}_{\Sigma_j}Y
&=\operatorname{div}Y-\langle DY\,N,N\rangle \\
&=b+rb_r+c_z+(b-c_z)N_z^2 -\left(b_z+\frac{c_r}{r}\right)SN_z-\frac{b_r}{r}S^2.
\end{aligned}
\end{equation}
Since \(zF_{rr}=F_z\) at \(r=0\),
\[
b=\frac13F-\frac23zF_z, \qquad \frac{c_r}{r}=2F_z, \qquad \frac{b_r}{r}=\frac13F_{rr}-\frac23F_{zz} \quad \text{ at }r=0.
\]
Hence the coefficient of \(N_z^2\) after substituting \(S=X\cdot N-zN_z\) in \eqref{eq:axis-test} satisfies
\[
b-c_z+z\left(b_z+\frac{c_r}{r}\right)-z^2\frac{b_r}{r}=-\frac13F \quad \text{ at }r=0.
\]
Since \(N_z^2\le\mathcal T_{k_j}\), the measures \(\beta_j^{-2}N_z^2\mu_j\) have uniformly bounded mass on \(\pi^{-1}(D_4)\) by \eqref{cat:fold:eq:relative-secondary-scale}. By \eqref{cat:fold:eq:relative-roof-strong} and \(|h|\le Ck_j\), they converge weakly to zero on \(\pi^{-1}(D_1)\cap\{r>0\}\).
By the Cauchy--Schwarz inequality, \eqref{cat:fold:eq:relative-secondary-scale}, and \eqref{cat:fold:eq:relative-radial-small},
\[
\int_{\Sigma_j\cap\operatorname{spt}Y} \left((X\cdot N)^2+|X\cdot N| |N_z|\right)\mathrm d\mu_j=o(\beta_j^2).
\]
Subtract the first variation identities for \(\Sigma_j\) and \(C_j\), using \eqref{eq:axis-test} and \(N_z=S=0\) on \(C_j\). Since \(b+rb_r+c_z=F\) near the axis, \eqref{cat:fold:eq:radial-scalar-comparison} gives
\begin{equation}\label{cat:fold:eq:axial-mass-identity}
\begin{aligned}
0\le\frac13\int FN_z^2\mathrm d\mu_j &=\int F \mathrm d(\mu_j-\|C_j\|)+o(\beta_j^2) \\
&=-\int_0^R f'(\rho)\left(\mu_j(B_\rho(0))-2\pi\rho^2\right)\mathrm d\rho +o(\beta_j^2)\le o(\beta_j^2),
\end{aligned}
\end{equation}
since \(f'\le0\) and \(\mu_j(B_\rho(0))\le2\pi\rho^2\).
This proves \eqref{cat:fold:eq:strong-radial-mass}. For \(K\Subset D_1\), choose \(a<1\) and \(F\) so that \(K\Subset D_a\) and \(F=1\) on \(B_a(0)\). Since \(|h|\le Ck_j\) and \(S=X\cdot N-zN_z\),
\[
\int_{\pi^{-1}(K)}(N_z^2+S^2)\mathrm d\mu_j=o(\beta_j^2).
\]

For \(i=1,2\), put \(q_j=(u_{j,i}-f_{j,i})/\beta_j\). On the graph of \(u_{j,i}\) over \(K_j\), take the unit normal \(N\) with \(N_h>0\). Then
\[
\beta_j(q_j-x\partial_xq_j) =\sqrt{1+|\nabla u_{j,i}|^2} S-(f_{j,i}-x\partial_xf_{j,i}).
\]
By \eqref{cat:fold:eq:relative-height-fit} and Proposition~\ref{cat:fold:prop:relative-radial-affine}, \(\|q_j\|_{L^2(D_r)}\le C_r\) for \(0<r<1\), and \(q_j\to0\) in \(L^2_{\mathrm{loc}}(D_1\setminus\{x=0\})\). Choose \(0\le\chi\in C_c^\infty(D_1)\) with \(\chi=1\) near \(K\) and \(\partial_x\chi=0\) near \(\{x=0\}\). By Proposition~\ref{cat:fold:prop:relative-bad-set} and \eqref{cat:fold:eq:relative-radial-affine},
\begin{equation}\label{cat:fold:eq:height-transverse-euler-coercivity}
\begin{aligned}
\frac32\int\chi q_j^2 &=\int\chi q_j(q_j-x\partial_xq_j) -\frac12\int x\partial_x\chi q_j^2 \\
&\le\frac{C_\chi}{\beta_j} \left(\int_{\pi^{-1}(\operatorname{spt}\chi)}S^2\mathrm d\mu_j\right)^{1/2} +C_\chi k_j+o(1)\longrightarrow0.
\end{aligned}\qedhere
\end{equation}
\end{proof}

By Lemma~\ref{cat:fold:prop:axial-defect-exclusion} and \eqref{cat:fold:eq:relative-tilt}, Proposition~\ref{cat:fold:prop:strong-affine-approximation} follows once we prove that, for every \(K\Subset D_1\),
\[
\int_{\pi^{-1}(K)}\frac{(|N_x|-k_j|N_h|)^2}{1+k_j^2}\mathrm d\mu_j=o(\beta_j^2).
\]

For finite nonnegative Radon measures \(E,U\) in \(\mathbb R^d\), the Besicovitch covering theorem \cite[Chapter 1, Lemma 3.12]{SimonGMT} gives
\begin{equation}\label{cat:fold:eq:relative-axial-maximal}
E\left(\left\{p\in\operatorname{spt} E: \sup_{r>0}\frac{U(B_r(p))}{E(B_r(p))}>t\right\}\right) \le\frac{C_d}{t}U(\mathbb R^d), \qquad t>0.
\end{equation}

\begin{lemma}
\label{cat:fold:prop:fixed-tilt-tail-small}
For every \(D\Subset D_1\),
\begin{equation}\label{cat:fold:eq:relative-fourth-tilt}
\int_{\pi^{-1}(D)}\mathcal T_{k_j}^{ 2}\mathrm d\mu_j=o(\beta_j^2),
\end{equation}
and
\begin{equation}\label{cat:fold:eq:relative-bad-set-little}
|D\setminus K_j| +\mu_j\left(\pi^{-1}(D\setminus K_j)\right) +\|V_j-\Gamma_j\|_{\mathrm{TV},\pi^{-1}(D)} =o(\beta_j^2).
\end{equation}
\end{lemma}
\begin{proof}
Fix \(D\Subset D_1\) and \(0<\alpha\le1\), and put \(F_j=\pi^{-1}(D)\cap\{\mathcal T_0\ge\alpha\}\). Since \(\mathcal T_{k_j}\le2\mathcal T_0+2k_j^2\), for all large \(j\),
\[
\int_{\pi^{-1}(D)}\mathcal T_{k_j}^{ 2}\mathrm d\mu_j \le C\alpha\beta_j^2+\mu_j(F_j).
\]
It suffices to prove \(\mu_j(F_j)=o(\beta_j^2)\) for every fixed \(0<\alpha\le1\).

Choose a cylinder \(\Omega\Subset\pi^{-1}(D_1)\) with \(\overline D\subset\Omega\cap\{h=0\}\). Put
\[
\nu_j=\mathcal T_{k_j}\mu_j\llcorner\Omega, \qquad \lambda_j=N_z^2\mu_j\llcorner\Omega.
\]
By \eqref{cat:fold:eq:relative-secondary-scale} and \eqref{cat:fold:eq:axial-and-transverse-support-small}, \(\nu_j(\Omega)\le C\beta_j^2\) and \(\lambda_j(\Omega)=o(\beta_j^2)\).

Choose \(0<\epsilon_j\to0\) with \(\lambda_j(\Omega)=o(\epsilon_j\beta_j^2)\). By \eqref{cat:fold:eq:relative-axial-maximal}, there are sets \(G_j\) with \(\nu_j(\mathbb R^3\setminus G_j)=o(\beta_j^2)\) such that, for \(p\in G_j\) and \(r>0\),
\begin{equation}\label{cat:fold:eq:all-radius-axial-good-centers}
\lambda_j(\overline B_r(p))\le\epsilon_j\nu_j(\overline B_r(p)).
\end{equation}

On \(F_j\), \(\mathcal T_{k_j}\ge\alpha/2\) for large \(j\), so
\[
\mu_j(F_j\setminus G_j) \le\frac2\alpha\nu_j(\mathbb R^3\setminus G_j)=o(\beta_j^2).
\]

Put \(\delta=\pi\alpha/4\). For \(p\in F_j\cap G_j\), set
\begin{equation}\label{cat:fold:eq:tilt-tail-stopping-radius}
r=r_j(p):=\sup\{t>0:\nu_j(\overline B_t(p))\ge\delta t^2\}.
\end{equation}
\(t^{-2}\nu_j(\overline B_t(p))\to\pi\mathcal T_{k_j}(p) \ge\pi\alpha/2\) as \(t\downarrow0\), so \(r>0\). The bound \(\nu_j(\Omega)\le C\beta_j^2\) and continuity from above give for \(R\ge1\)
\begin{equation}\label{cat:fold:eq:tilt-tail-stopping-growth}
r\le C_\alpha\beta_j, \qquad \nu_j(\overline B_r(p))=\delta r^2,\quad \nu_j(\overline B_{Rr}(p))\le\delta R^2r^2.
\end{equation}
For any sequence \(p_j\in F_j\cap G_j\), put \(r_j=r_j(p_j)\). After passing to a subsequence, \(|r_j^{-1}(\Sigma_j-p_j)|\to W\) as varifolds on \(\mathbb R^3\). For every fixed \(R\ge1\), \eqref{cat:fold:eq:all-radius-axial-good-centers} gives
\[
r_j^{-2}\lambda_j(B_{Rr_j}(p_j))\le\epsilon_j\delta R^2\longrightarrow0.
\]
Thus \(e_z\) is tangent to \(W\) almost everywhere. By Lemma~\ref{lem:constant-direction}, \(W=\sum_{l=1}^q|P_l|\) for affine planes \(P_l\) and \(q\le2\). By \eqref{cat:fold:eq:tilt-tail-stopping-growth},
\[
\pi\sum_l \mathcal T_0(P_l)=\lim_{R\to\infty}R^{-2}\int_{B_R}\mathcal T_0\mathrm d\|W\|\le\delta.
\]
Since \(\mathcal T_0(P_l)\le\alpha/4\),
\[
\limsup_j r_j^{-2}\mu_j\left(\overline B_{r_j}(p_j)\cap F_j\right) \le W\left(\overline B_1\times\{T:\mathcal T_0(T)\ge\alpha\}\right)=0.
\]
Thus, uniformly for \(p\in F_j\cap G_j\),
\begin{equation}\label{cat:fold:eq:tilt-tail-stopping-small}
\mu_j\left(\overline B_{r_j(p)}(p)\cap F_j\right)=o(1)\nu_j\left(\overline B_{r_j(p)}(p)\right).
\end{equation}
By the Besicovitch covering theorem, choose a subcover \(\{\overline B_{r_j(p_i)}(p_i)\}\) of \(F_j\cap G_j\) with bounded overlap. By \eqref{cat:fold:eq:tilt-tail-stopping-small},
\[
\mu_j(F_j)\le\mu_j(F_j\setminus G_j)+o(1)\sum_i\nu_j(\overline B_{r_j(p_i)}(p_i)) \le o(\beta_j^2)+o(1)\nu_j(\Omega)=o(\beta_j^2).
\]

As in the proof of Proposition~\ref{cat:fold:prop:relative-bad-set}, we can choose pairwise disjoint disks \(D_{r_i}(y_i)\) such that
\[
D\setminus K_j\subset\bigcup_iD_{5r_i}(y_i), \qquad r_i=O(\beta_j).
\]
Since \(\bigcup_iD_{r_i}(y_i)\) lies in a fixed compact subset of \(D_1\), \eqref{cat:fold:eq:relative-fourth-tilt} gives
\[
\begin{aligned}[b]
&|D\setminus K_j|+\mu_j\left(\pi^{-1}(D\setminus K_j)\right)+\|V_j-\Gamma_j\|_{\mathrm{TV},\pi^{-1}(D)} \\
\le{}& C_\tau\int_{\pi^{-1}(\bigcup_iD_{r_i}(y_i))}\mathcal T_{k_j}^2\mathrm d\mu_j=o(\beta_j^2).
\end{aligned}\qedhere
\]
\end{proof}

\begin{proof}[Proof of Proposition~\ref{cat:fold:prop:strong-affine-approximation}]
With \(b_j=(u_{j,1}+u_{j,2})/2-m_j\), put \(g_j=(u_{j,2}-u_{j,1})/2\ge0\) and choose \(D\Subset D_1\). On the graphs over \(K_j\), \(N_z^2\mathrm d\mu_j=(\partial_z u_{j,i})^2(1+|\nabla u_{j,i}|^2)^{-1/2}\mathrm d x\mathrm d z\). By \eqref{cat:fold:eq:relative-mean-gradient}, \eqref{cat:fold:eq:axial-and-transverse-support-small}, and \eqref{cat:fold:eq:relative-bad-set-little},
\begin{equation}\label{cat:fold:eq:scalar-gap-inputs}
\|\nabla b_j\|_{L^2(D)}+\|\partial_zg_j\|_{L^2(D)}=o(\beta_j).
\end{equation}
Set \(\zeta_j=|\partial_xg_j|-k_j\). Since \(|\zeta_j|\le\min_\sigma|\nabla u_{j,i}-\sigma k_je_x|+|\nabla b_j|\), \eqref{cat:fold:eq:relative-reconstruction} implies \(\int_D\zeta_j^2=O(\beta_j^2)\). By \eqref{cat:fold:eq:relative-graph} and \eqref{cat:fold:eq:relative-bad-set-little}, for \(0\le\chi\in C_c^\infty(D_1)\),
\begin{equation}\label{cat:fold:eq:remaining-gap-gradient-energy}
\int(\chi\circ\pi)\mathcal T_{k_j}\mathrm d\mu_j \le C\int\chi \zeta_j^2+o(\beta_j^2).
\end{equation}
Put
\begin{equation}\label{cat:fold:eq:moving-tail-scalar-bounds}
w_j:=\frac{(\partial_xg_j)^2-k_j^2}{2} =k_j\zeta_j+\frac12\zeta_j^2\ge-\frac12k_j^2.
\end{equation}
By \eqref{cat:fold:eq:scalar-gap-inputs} and Lemma~\ref{cat:fold:prop:fixed-tilt-tail-small},
\begin{align}
&\int_Dw_j^2 \label{cat:fold:eq:area-factor-smallness}\\
\le{}& C_L\left(k_j^2\int_D\zeta_j^2 +\int_{\pi^{-1}(D)}\mathcal T_{k_j}^2\mathrm d\mu_j +|D\setminus K_j|+\|\nabla b_j\|_{L^2(D)}^2\right) =o(\beta_j^2).\notag
\end{align}
Recall that
\[
B(p)=\sqrt{1+|p|^2}I-\frac{p\otimes p}{\sqrt{1+|p|^2}}.
\]
Since \(B(-p)=B(p)\) and \(\nabla u_{j,i}=\nabla b_j+(-1)^i\nabla g_j\),
\[
\frac12\sum_i B(\nabla u_{j,i})=B(\nabla g_j)+O_L(|\nabla b_j|^2).
\]
With \(W_j=\sqrt{1+k_j^2}\), we have \(1+(\partial_xg_j)^2=W_j^2+2w_j\), so
\[
\frac12\sum_i \begin{pmatrix}B_{xx}(\nabla u_{j,i})\\B_{zz}(\nabla u_{j,i})\end{pmatrix} ={} \begin{pmatrix}W_j^{-1}-W_j^{-3}w_j\\W_j+W_j^{-1}w_j\end{pmatrix} +O_L\left(|\nabla b_j|^2+|\partial_zg_j|^2+w_j^2\right).
\]
By \eqref{cat:fold:eq:scalar-gap-inputs} and \eqref{cat:fold:eq:area-factor-smallness}, the remainder is \(o(\beta_j^2)\) in \(L^1(D)\). For \(Y=(-\partial_x\psi,\partial_z\psi)\), \(\psi\in C_c^\infty(D_1)\), \(\operatorname{tr}(B dY)=-B_{xx}\partial_{xx}\psi+B_{zz}\partial_{zz}\psi\). By \eqref{cat:fold:eq:relative-stress-error} and \eqref{cat:fold:eq:relative-bad-set-little}, there is a diagonal matrix \(\mathcal R_j\) of measures such that, for every \(D\Subset D_1\),
\begin{equation}\label{cat:fold:eq:moving-tail-elliptic-equation}
L_jw_j=\operatorname{div}\operatorname{div}\mathcal R_j, \qquad L_j=\partial_x^2+W_j^2\partial_z^2, \qquad |\mathcal R_j|(D)=o(\beta_j^2).
\end{equation}

By Lemma~\ref{cat:fold:prop:relative-offfold-compactness} and \eqref{cat:fold:eq:relative-bad-set-little}, \(\zeta_j/\beta_j\to0\) in \(L^2_{\mathrm{loc}}(D_1\setminus\{x=0\})\). Since \(\|\zeta_j\|_{L^2(D)}=O(\beta_j)\), we have \(\zeta_j/\beta_j\rightharpoonup0\) in \(L^2_{\mathrm{loc}}(D_1)\).

Since \(\zeta_j/k_j\ge-1\), Lemma~\ref{lem:quadratic-elliptic-convergence} below applies to \(\zeta_j/k_j\) with \(\varepsilon_j=\beta_j/k_j\to0\). Thus \(\int_D\zeta_j^2=o(\beta_j^2)\), and \eqref{cat:fold:eq:remaining-gap-gradient-energy} proves the proposition.
\end{proof}

\begin{lemma}\label{lem:quadratic-elliptic-convergence}
Let \(0<\varepsilon_j\to0\), \(c_j>0\) with \(c_j\to1\), and \(\zeta_j\in L^2_{\mathrm{loc}}(D_1)\) satisfy \(\zeta_j\ge-1\) and \(\zeta_j/\varepsilon_j\rightharpoonup0\) in \(L^2_{\mathrm{loc}}(D_1)\). Put
\[
w_j=\zeta_j+\tfrac12\zeta_j^2, \qquad L_j=\partial_x^2+c_j\partial_z^2.
\]
Suppose that \(\mathcal R_j\) are diagonal matrices of Radon measures on \(D_1\) with
\[
L_jw_j=\operatorname{div}\operatorname{div}\mathcal R_j, \qquad |\mathcal R_j|(D)=o(\varepsilon_j^2)
\]
for every \(D\Subset D_1\). Then \(\zeta_j/\varepsilon_j\to0\) in \(L^2_{\mathrm{loc}}(D_1)\).
\end{lemma}
\begin{proof}
Choose \(D\Subset\Omega\Subset D_1\). For a fundamental solution \(\Phi_j\) of \(L_j\), put
\[
P_j=\sum_{a\in\{x,z\}}\partial_{aa}\left(\Phi_j*(\mathcal R_{j,aa}\llcorner\Omega)\right), \qquad h_j=w_j-P_j.
\]
Then \(L_jh_j=0\) on \(\Omega\). Since \(\zeta_j/\varepsilon_j\rightharpoonup0\) in \(L^2(\Omega)\),
\[
\left\|\frac{w_j-\zeta_j}{\varepsilon_j}\right\|_{L^1(\Omega)}=\frac{\|\zeta_j\|_{L^2(\Omega)}^2}{2\varepsilon_j}=O(\varepsilon_j).
\]
Thus \(w_j/\varepsilon_j\to0\) in \(\mathcal D'(\Omega)\). Since \(|\mathcal R_j|(\Omega)=o(\varepsilon_j^2)\), also \(h_j/\varepsilon_j\to0\) in \(\mathcal D'(\Omega)\), and interior estimates give \(\|h_j\|_{L^\infty(D)}=o(\varepsilon_j)\).

With \(z=\sqrt{c_j}\widehat z\), \(L_j=\partial_x^2+\partial_{\widehat z}^2\). Since \(P_j=w_j-h_j\in L^1_{\mathrm{loc}}(\Omega)\), the weak \((1,1)\) estimate for Calder\'on--Zygmund operators, applied to smooth approximations of \(\mathcal R_j\llcorner\Omega\), gives
\[
\int_D\min\{P_j^2,1\}=2\int_0^1t|\{y\in D:|P_j(y)|>t\}|\mathrm d t\le C|\mathcal R_j|(\Omega)=o(\varepsilon_j^2),
\]
with \(C\) independent of \(j\).
Since \(\zeta_j\ge-1\),
\[
w_j\ge-\tfrac12, \qquad \zeta_j^2\le4\min\{w_j^2,1\}+2(w_j-1)_+.
\]
Together with \(w_j=P_j+h_j\), this gives
\begin{equation}\label{eq:quadratic-elliptic-tail}
\int_D\zeta_j^2\le2\int_D(w_j-1)_++o(\varepsilon_j^2).
\end{equation}

Let \(F_j\) be the set of Lebesgue points of \(\zeta_j^2\) in \(\{w_j>1\}\). On \(F_j\), \(\zeta_j^2>1/2\) and \(w_j\le2\zeta_j^2\). Put \(\nu_j=\zeta_j^2\mathcal L^2\llcorner\Omega\), and choose \(0<\eta_j\to0\) with \(|\mathcal R_j|(\Omega)=o(\eta_j\varepsilon_j^2)\). By \eqref{cat:fold:eq:relative-axial-maximal}, there are sets \(G_j\) with \(\nu_j(\Omega\setminus G_j)=o(\varepsilon_j^2)\) such that, for \(y\in G_j\) and \(r>0\),
\[
|\mathcal R_j|(D_r(y)\cap\Omega)\le\eta_j\nu_j(D_r(y)).
\]
Fix a sufficiently small \(0<\delta<\pi/2\). For every \(y\in D\cap F_j\cap G_j\), let \(r=r_j(y)>0\) be the largest radius with \(\nu_j(D_r(y))\ge\delta r^2\). Since \(\nu_j(\Omega)=O(\varepsilon_j^2)\),
\[
r\le C_\delta\varepsilon_j, \qquad \nu_j(D_r(y))=\delta r^2, \qquad \nu_j(D_{2r}(y))<4\delta r^2.
\]
Thus \(D_{2r}(y)\Subset\Omega\) for large \(j\).

For \(y_j\in D\cap F_j\cap G_j\), put \(r_j=r_j(y_j)\) and \(v_j(\xi)=w_j(y_j+r_j\xi)\). On \(D_2\),
\[
\begin{gathered}
v_j\ge-\tfrac12, \qquad \|v_j\|_{L^1(D_2)}\le C(\sqrt\delta+\delta), \\
L_jv_j=\operatorname{div}\operatorname{div}\widetilde{\mathcal R}_j, \qquad |\widetilde{\mathcal R}_j|(D_2)=\frac{|\mathcal R_j|(D_{2r_j}(y_j))}{r_j^2}\le4\eta_j\delta\longrightarrow0.
\end{gathered}
\]
After passing to a subsequence, \(v_j\to v_\infty\) in distributions on \(D_2\), where \(\Delta v_\infty=0\). Put
\[
\widetilde P_j=\sum_{a\in\{x,z\}}\partial_{aa}\left(\Phi_j*(\widetilde{\mathcal R}_{j,aa}\llcorner D_2)\right), \qquad \widetilde h_j=v_j-\widetilde P_j.
\]
Then \(L_j\widetilde h_j=0\) on \(D_2\), \(\widetilde P_j\to0\) locally in measure, and \(\widetilde h_j\to v_\infty\) locally uniformly. For \(0\le\psi\in C_c^\infty(D_2)\),
\[
\int\psi|v_j-v_\infty|=\int\psi(v_j-v_\infty)+2\int\psi(v_\infty-v_j)_+\longrightarrow0,
\]
since \((v_\infty-v_j)_+\le(v_\infty+1/2)_+\). By the mean-value inequality,
\[
\|v_\infty\|_{L^\infty(D_1)}
\le C\|v_\infty\|_{L^1(D_2)}
\le C(\sqrt\delta+\delta)<1.
\]
Hence \(\int_{D_1}(v_j-1)_+\to0\), so
\[
\postdisplaypenalty=10000
\int_{D_{r_j(y)}(y)}(w_j-1)_+=o(r_j(y)^2)=o(1)\nu_j(D_{r_j(y)}(y))
\]
uniformly for \(y\in D\cap F_j\cap G_j\).

By the Besicovitch covering theorem, there are disks \(D_{r_j(y_i)}(y_i)\) with bounded overlap covering almost every point of \(D\cap F_j\cap G_j\). Hence
\[
\begin{aligned}
\int_D(w_j-1)_+\le&\,2\nu_j(\Omega\setminus G_j)+o(1)\sum_i\nu_j(D_{r_j(y_i)}(y_i))\\
\le&\,o(\varepsilon_j^2)+o(1)\nu_j(\Omega)=o(\varepsilon_j^2).
\end{aligned}
\]
By \eqref{eq:quadratic-elliptic-tail}, \(\int_D\zeta_j^2=o(\varepsilon_j^2)\).
\end{proof}

\subsection{Comparison between scales}\label{cat:sec:scale-comparison}

Recall \(F(R)=e(R)+E_h(R)+E_C(R)\) from \eqref{cat:fold:eq:geometric-buffer-error} and \(q(R)=R^{-4}\int_{M\cap B_R}(X\cdot N)^2\mathrm d\mu\).

\begin{proof}[Proof of Proposition~\ref{cat:fold:prop:geometric-buffer-improvement}]
We argue by contradiction. Suppose that there are \(\varepsilon>0\) and \(R_j\to\infty\) such that
\[
F(R_j)>\varepsilon F(\Lambda R_j)+j q(\Lambda R_j).
\]
Put \(s_j=\Lambda R_j\) and \(\beta_j=F(R_j)^{1/2}>0\). Then \(F(s_j)\le\varepsilon^{-1}\beta_j^2\), \(q(s_j)=o(\beta_j^2)\) and, by Lemma~\ref{cat:fold:lem:first-scale-relative-small} and \(a_{R_j}/a_{s_j}\to1\), \(\beta_j=o(a_{s_j})\). Thus \eqref{cat:fold:eq:secondary-sequence-assumptions} holds. After passing to a subsequence, take \(Q_j,\Sigma_j,u_{j,i},f_{j,i}\) from Section~\ref{cat:sec:affine-limits}. By Proposition~\ref{cat:fold:prop:strong-affine-approximation},
\begin{equation}\label{cat:fold:eq:buffered-vanishing-local-outputs}
\sum_i\int_{D_{1/16}}|u_{j,i}-f_{j,i}|^2\mathrm d A +\int_{\Sigma_j\cap B_{1/16}}\mathcal T_{k_j}\mathrm d\mu_j =o(\beta_j^2).
\end{equation}

Write \(s=s_j\), \(R=R_j\). Lemma~\ref{cat:lem:initial-comparison}, applied at \(R,2R,\ldots,s/2\), gives
\begin{equation}\label{cat:fold:eq:buffered-base-angle}
|\Pi_R-\Pi_s|\le C a_s.
\end{equation}
Over \(P_R\), write
\[
R^{-1}C_s=\Lambda Q_j^{-1}\left(\bigcup_i\operatorname{graph}f_{j,i}\right)
=\bigcup_i\operatorname{graph}\widetilde f_i, \qquad \widetilde f_1\le\widetilde f_2.
\]
By \eqref{cat:fold:eq:buffered-base-angle} and \(\|Q_j-I\|=o(a_s)\), \(\Lambda Q_j^{-1}(\operatorname{graph}u_{j,i})\) are graphs over \(D_2\subset P_R\), contained in \(\{X:\operatorname{dist}(X,P_R)\le Ca_s\}\). Apply \eqref{cat:fold:eq:relative-bad-set} to \(\Sigma_j\) with \(C=\{h=\pm k_jx\}\) and to \(R^{-1}M\) with \(C=R^{-1}C_s\). Since \(E_C(s)\le\varepsilon^{-1}\beta_j^2\) and \(D_{2/\Lambda}\Subset D_{1/16}\), \eqref{cat:fold:eq:reconstruction-height-comparison}, \eqref{cat:fold:eq:graph-height-comparison}, and \eqref{cat:fold:eq:buffered-vanishing-local-outputs} give
\[
\begin{aligned}
\|U_R-\widetilde f\|_{L^2_\omega(D_1)}^2 &\le\frac{Ca_s^2}{R^2}\int_{M\cap B_{8s}}\mathcal T_{C_s}^2\mathrm d\mu +C\sum_i\int_{D_{2/\Lambda}}|u_{j,i}-f_{j,i}|^2\mathrm d A \\
&\le Ca_s^2\beta_j^2+o(\beta_j^2)=o(\beta_j^2).
\end{aligned}
\]
For an affine plane \(P\), let \(\Pi_P\) be the orthogonal projection onto the parallel plane through \(0\). By \eqref{cat:fold:eq:buffered-height-fit-definition} and \eqref{cat:fold:eq:fit-stability}, after relabeling the planes \(P_{R,i}\) in \(C_R\) and \(P_{s,i}\) in \(C_s\),
\[
E_h(R)+\max_i|\Pi_{P_{R,i}}-\Pi_{P_{s,i}}|^2 \le C\|U_R-\widetilde f\|_{L^2_\omega(D_1)}^2=o(\beta_j^2).
\]
Since \(8/\Lambda<1/16\), \eqref{eq:area-bound} and \eqref{cat:fold:eq:buffered-vanishing-local-outputs} give
\[
E_C(R)\le C\max_i|\Pi_{P_{R,i}}-\Pi_{P_{s,i}}|^2 +\frac2{R^2}\int_{M\cap B_{8R}}\mathcal T_{C_s}\mathrm d\mu =o(\beta_j^2).
\]

Choose \(\varphi\) as in Lemma~\ref{cat:fold:prop:axial-defect-exclusion}, supported in \(B_{1/\Lambda}\), with \(\int_{\{h=0\}}\varphi\mathrm d A=1\). By the monotonicity formula and \eqref{cat:fold:eq:strong-radial-mass},
\[
\begin{aligned}[b]
\beta_j^2=F(R_j)&=e(R_j)+o(\beta_j^2) \\
&\le2-\int\varphi\mathrm d\mu_j+o(\beta_j^2)=o(\beta_j^2).
\end{aligned}\qedhere
\]
\end{proof}

\begin{proof}[Proof of Theorem~\ref{cat:thm:folded-exclusion}]
Fix a sufficiently large \(R_*\), put \(R_n=\Lambda^nR_*\), and write \(e_n=e(R_n)\), \(q_n=q(R_n)\) and \(F_n=F(R_n)\). By \eqref{eq:monotonicity} and \eqref{cat:fold:eq:geometric-buffer-improvement} with \(\varepsilon=1/2\),
\[
F_n \le\tfrac12F_{n+1}+Cq_{n+1},\quad q_{n+1} \le\Lambda^{-4}q_n+\pi(e_n-e_{n+1}).
\]
Since \(0\le e_n\le F_n\to0\) by Lemma~\ref{cat:fold:lem:first-scale-relative-small}, summing these inequalities from \(k=n\) to \(N\) and letting \(N\to\infty\) gives, for some \(K>1\),
\[
T_n:=\sum_{k\ge n}(F_k+q_k) \le C\sum_{k\ge n}q_k \le K(F_n+q_n)=K(T_n-T_{n+1})<\infty.
\]
Thus \(T_{n+1}\le(1-K^{-1})T_n\). For \(R_n\le R\le R_{n+1}\),
\[
q(R)\le\Lambda^4q_{n+1}\le\Lambda^4T_{n+1}=O(R^{-\gamma})
\]
for some \(\gamma>0\).

The smallest eigenvalue of \(C(R)\) is
\[
\lambda_R=\mathbf{n}_R\cdot C(R)\mathbf{n}_R
=R^{-4}\int_M(X\cdot\mathbf{n}_R)^2\omega(|X|/R)\mathrm d\mu.
\]
By \eqref{eq:area-bound} and \eqref{cat:eq:covariance-derivative},
\[
R|\lambda_R'|=R|\mathbf{n}_R\cdot C'(R)\mathbf{n}_R|\le C\left(\sqrt{q(R)}+q(R)\right).
\]
Since \(0\le\lambda_R\le Ca_R^2\to0\), integration over \([R,\infty)\) gives
\[
\lambda_R\le C\int_R^\infty\left(\sqrt{q(t)}+q(t)\right)\frac{\mathrm d t}{t} =O(R^{-\gamma/2}).
\]

For \(t=R/32\), \(C(t)\ge c\Pi_t\) by Lemma~\ref{cat:lem:initial-comparison}. Since \(\omega\) is positive on \(B_{1/2}\),
\[
\|\Pi_t-\Pi_R\|^2\le C \mathbf{n}_R\cdot C(t)\mathbf{n}_R\le C\lambda_R.
\]
By \eqref{eq:area-bound},
\[
a_t^2\le 2t^{-4}\int_{M\cap B_{8t}}(X\cdot \mathbf{n}_R)^2\mathrm d\mu+C\|\Pi_t-\Pi_R\|^2 \le C\lambda_R=O(R^{-\gamma/2}).
\]
This contradicts \(a_t\ge c_\varepsilon t^{-\varepsilon}\) from Theorem~\ref{cat:thm:regimes}, with \(0<\varepsilon<\gamma/4\).
\end{proof}

\section{Reflection symmetry and catenoid rigidity}\label{sec:catenoid-rigidity}

Theorems~\ref{cat:thm:regimes} and~\ref{cat:thm:folded-exclusion} give
\begin{equation}\label{cat:eq:constant-doubling-input}
a_{2R}/a_R\longrightarrow1/2.
\end{equation}
Thus every limit \(v\) in Theorem~\ref{cat:thm:compat-arbitrary-limits} equals \(b_*=(2|D_8|)^{-1/2}\), and \(\mathsf H_R\le C_\varepsilon R^\varepsilon\) for every \(\varepsilon>0\). Also \(\mathsf H_R=Ra_R\to\infty\) by Lemma~\ref{cat:lem:physical-height}.

\subsection{Projection onto the limiting plane}\label{cat:sec:fixed-plane-capacity}

Fix \(0\le\eta\in C_c^\infty(D_1)\) radial with \(\int\eta=1\), and put \(\eta_r(x)=r^{-2}\eta(x/r)\).
\begin{proposition}\label{cat:prop:constant-ratio-geometry}
The planes \(P_R\) converge to a plane \(P_\infty\), and orthogonal projection \(\pi:\mathbb R^3\to P_\infty\) restricts to a proper map on \(M\). In coordinates \(X=(x,z)\), with \(x=(x_1,x_2)\) and \(P_\infty=\{z=0\}\), write \(N=(\pi N,N_z)\). There is \(z_\infty\in\mathbb R\) such that \(u=z-z_\infty\) satisfies
\begin{equation}\label{cat:eq:localized-first-moment-decay}
\sup_{|x|\le LR}|(\eta_{\rho R}*\pi_\#(u\mu))(x)| \le C_{L,\rho}\frac{\mathsf H_R^3}{R^2}
\end{equation}
for every \(L,\rho>0\) and all sufficiently large \(R\).
\end{proposition}
\begin{proof}
By \eqref{cat:eq:drift} and Lemma~\ref{cat:lem:initial-comparison}, \(P_R\to P_\infty\) with \(\|\Pi_R-\Pi_\infty\|\le Ca_R^3\). By Lemma~\ref{cat:lem:height-tilt},
\[
\sup_{M\cap B_R}|z|\le C\mathsf H_R=o(R).
\]
Thus \(M\cap\pi^{-1}(D_{LR})\subset B_{2LR}\) for large \(R\). Since \(M\) is proper, so is \(\pi|_M\). By \eqref{eq:area-bound},
\begin{equation}\label{cat:eq:constant-ratio-cylinder-bounds}
\sup_{M\cap\pi^{-1}(D_{LR})}|z|\le C_L\mathsf H_R, \qquad \mu(\pi^{-1}(D_{LR}))\le C_LR^2.
\end{equation}

Put
\[
f_R=\eta_R*\pi_\#\mu, \qquad m_R=\frac{\eta_R*\pi_\#(z\mu)}{f_R}.
\]
Since \(f_R(R \cdot)\to2\) locally uniformly, Proposition~\ref{cat:prop:center} gives
\[
\sup_{R\le t\le2R}\sup_{D_{LR}}|m_t-m_R| \le C_L\mathsf H_R^3/R^2,\quad \sup_{D_{LR}}|\nabla m_R| \le C_L\mathsf H_R/R.
\]
Since \(\mathsf H_{2R}/\mathsf H_R\to1\),
\[
\sup_{t\ge R}\sup_{D_{LR}}|m_t-m_R|
\le C_L\sum_{k=0}^\infty\frac{\mathsf H_{2^kR}^3}{(2^kR)^2}
\le C_L\frac{\mathsf H_R^3}{R^2}.
\]
Hence \(m_R\) converges locally uniformly. Since \(\nabla m_R\to0\), the limit is a constant \(z_\infty\in\mathbb R\), with
\[
\sup_{D_{LR}}|m_R-z_\infty| \le C_L\mathsf H_R^3/R^2.
\]
Since \(\mathsf H_{\rho R}\le C_\rho\mathsf H_R\) for \(\rho>0\), we have on \(D_{LR}\)
\[
|\eta_{\rho R}*\pi_\#(u\mu)|=f_{\rho R}|m_{\rho R}-z_\infty|\le C_{L,\rho}\frac{\mathsf H_R^3}{R^2}.
\]
\end{proof}

Let \(\Omega_0\subset P_\infty\) be the open set of \(x\) for which \(M\cap\pi^{-1}(x)\) has exactly two points, both with \(|N_z|>1/2\), and put \(E=P_\infty\setminus\Omega_0\).

\begin{corollary}
\label{cat:cor:constant-ratio-outer-fibers}
For every \(L>0\),
\begin{equation}\label{cat:eq:constant-ratio-bad-fiber-collapse}
\sup_{M\cap\pi^{-1}(E\cap D_{LR})} \frac{|u|}{\mathsf H_R}\longrightarrow0\qquad\text{as }R\to\infty.
\end{equation}
\end{corollary}
\begin{proof}
Using \(|\pi N|\le|\Pi_RN|+Ca_R^3\), \eqref{cat:eq:constant-ratio-cylinder-bounds}, and Proposition~\ref{cat:prop:strong-height-limits} with \(v=b_*\),
\begin{equation}\label{cat:eq:constant-ratio-tilt}
\int_{\pi^{-1}(D_{LR})}|\pi N|^2 \mathrm d\mu \le2\int_{M\cap B_{2LR}}|\Pi_RN|^2 \mathrm d\mu +C_La_R^4\mathsf H_R^2 =o(\mathsf H_R^2).
\end{equation}
For fixed \(0<\rho<1\), \eqref{cat:eq:localized-first-moment-decay} and \(f_{\rho R}\ge1\) on \(D_{LR}\) give
\[
\sup_{D_{LR}}|m_{\rho R}-z_\infty|\le C_{L,\rho}\frac{\mathsf H_R^3}{R^2}=o(\mathsf H_R).
\]
By \eqref{cat:eq:constant-ratio-cylinder-bounds}, \eqref{cat:eq:constant-ratio-tilt}, and Proposition~\ref{cat:prop:compat-high-graphs} applied to \(R^{-1}(M-z_\infty e_3)\), \eqref{cat:eq:constant-ratio-bad-fiber-collapse} follows.
\end{proof}

\begin{proposition}
\label{cat:prop:constant-ratio-capacity}
There are \(\xi_R\in C_c^\infty(D_5)\) such that
\[
0\le\xi_R\le1, \qquad \xi_R=1\text{ near }(R^{-1}E)\cap\overline D_4,\quad \|\xi_R\|_{W^{1,2}(D_5)}\longrightarrow0.
\]
For all large \(R\), some \(t_R\in(1,2)\) satisfies \(\partial D_{t_RR}\subset\Omega_0\).
\end{proposition}
\begin{proof}
For \(x\in\Omega_0\), put \(u_+(x)=\max\{u(X):X\in M,\ \pi(X)=x\}\). Set \(\chi_R=1\) on \(D_5\cap R^{-1}E\). On \(D_5\cap R^{-1}\Omega_0\), put
\[
\chi_R(y)=\min\left\{1,\left(2-\frac{2u_+(Ry)}{b_*\mathsf H_R}\right)_+\right\}.
\]
For large \(R\), \(\chi_R=1\) near \(R^{-1}E\) in \(D_5\), by \eqref{cat:eq:constant-ratio-bad-fiber-collapse} and properness of \(\pi|_M\). Since \(u_+\in C^\infty(\Omega_0)\), \(\chi_R\) is locally Lipschitz.

Proposition~\ref{cat:prop:lipschitz-approximation} with \(q=2\) and \eqref{cat:eq:constant-ratio-tilt} give \(|D_5\cap R^{-1}E|\to0\). Since \(R\|\Pi_R-\Pi_\infty\|+|z_\infty|=o(\mathsf H_R)\), Lemma~\ref{cat:lem:compat-weak} and Theorem~\ref{cat:thm:regimes} imply
\[
\int_{D_5\cap R^{-1}\Omega_0}\left|\frac{u_+(Ry)}{\mathsf H_R}-b_*\right|^2\mathrm d A(y)\longrightarrow0.
\]
Thus \(\chi_R\to0\) in \(L^2(D_5)\). The area formula and \eqref{cat:eq:constant-ratio-tilt} give
\[
\int_{D_5}|\nabla\chi_R|^2\mathrm d A \le\frac{C}{\mathsf H_R^2} \int_{\pi^{-1}(D_{5R})}|\pi N|^2 \mathrm d\mu \longrightarrow0.
\]

Choose \(0\le\xi_R\le1\) in \(C_c^\infty(D_5)\) with \(\xi_R=1\) near \((R^{-1}E)\cap\overline D_4\) and \(\|\xi_R\|_{W^{1,2}(D_5)}\le C\|\chi_R\|_{W^{1,2}(D_5)}\).
For some \(t_R\in(1,2)\),
\[
\sup_{\partial D_{t_R}}|\xi_R| \le C\|\xi_R\|_{W^{1,2}(\partial D_{t_R})} \le C\|\xi_R\|_{W^{1,2}(D_5)}\longrightarrow0.
\]
Since \(\xi_R=1\) on \((R^{-1}E)\cap D_4\), \(\partial D_{t_RR}\subset\Omega_0\) for large \(R\).
\end{proof}

\subsection{Reflection symmetry}\label{cat:sec:reflection}

\begin{proposition}
\label{cat:prop:compact-midpoint-principle}
Let \(\Sigma\subset\overline D_r\times\mathbb R\) be a smooth compact embedded minimal surface, possibly disconnected, with \(\Sigma\setminus\partial\Sigma\subset D_r\times\mathbb R\). Suppose
\[
\partial\Sigma=\operatorname{graph}w_-\cup\operatorname{graph}w_+,
\]
where \(w_\pm:\partial D_r\to\mathbb R\) are smooth and \(w_-<w_+\). For any two distinct points \((x,s),(x,t)\in\Sigma\),
\begin{equation}\label{cat:eq:compact-midpoint-interval}
\min_{\partial D_r}\frac{w_++w_-}{2} \le\frac{s+t}{2} \le\max_{\partial D_r}\frac{w_++w_-}{2}.
\end{equation}
If \(w_++w_-=2c\), then \(\Sigma\) is invariant under reflection in \(z=c\), and \(\Sigma\cap\{z>c\}\) and \(\Sigma\cap\{z<c\}\) are smooth graphs.
\end{proposition}
\begin{proof}
Put \(\beta=\max_{\partial D_r}(w_++w_-)/2\). For \(c>\beta\), we have \(w_-<c\) and \(2c-w_+>w_-\). By \cite[Theorem 1, p.794]{SchoenSymmetry}, applied to \(\Sigma-ce_3\), \(2c-s\ge t\) whenever \((x,s),(x,t)\in\Sigma\) and \(s>c>t\).

If \(s>t\) and \((s+t)/2>\beta\), taking \(\max\{\beta,t\}<c<(s+t)/2\) gives \(s+t\le2c<s+t\). Replacing \(z\) by \(-z\) gives the lower bound. If \(w_++w_-=2c\), apply \cite[Theorem 2, p.795]{SchoenSymmetry} to \(\Sigma-ce_3\).
\end{proof}

\begin{theorem}
\label{cat:thm:constant-ratio-reflection}
\label{cat:thm:constant-ratio-bigraph}
There is a nonempty locally finite family of pairwise disjoint compact strictly convex domains \(K_i\) with smooth boundary such that \(\Omega=\mathbb R^2\setminus\bigcup_i K_i\) is connected and unbounded. There is \(q\in C^\infty(\Omega)\cap C^0(\overline\Omega)\) such that \(q>0\) in \(\Omega\), \(q=0\) on \(\partial\Omega\), and
\begin{equation}\label{cat:eq:constant-ratio-bigraph}
M=\{(x,z_\infty\pm q(x)):x\in\overline\Omega\}.
\end{equation}
\(N_z=0\) on \(M\cap\{u=0\}\).
\end{theorem}
\begin{proof}
For \(R>0\), define
\[
L(R)=\sup_{\substack{X,Y\in M, X\ne Y\\ \pi X=\pi Y\in\overline D_R}} \frac{|u(X)+u(Y)|}{2}.
\]
\(L\) is nondecreasing, and \eqref{cat:eq:constant-ratio-cylinder-bounds} gives \(L(R)\le C\mathsf H_R\) for large \(R\). Suppose \(L\not\equiv0\), so \(L(R)\ge c>0\) for large \(R\). Let \(u_+>u_-\) be the two values of \(u\) over \(\Omega_0\), and choose \(\xi_R\) as in Proposition~\ref{cat:prop:constant-ratio-capacity}. On \(R^{-1}\Omega_0\cap D_4\), put
\[
v_R(y)=(1-\xi_R(y))\frac{u_+(Ry)+u_-(Ry)}{2L(8R)},
\]
and set \(v_R=0\) on \(D_4\setminus R^{-1}\Omega_0\). Since \(|u_+(Ry)+u_-(Ry)|\le2L(8R)\), we have \(|v_R|\le1\).

Apply Proposition~\ref{cat:prop:lipschitz-approximation} to \((4R)^{-1}M\) with fixed small \(\lambda>0\). Then \(E\cap D_{4R}\subset D_{4R}\setminus(4R K_\lambda)\). By \eqref{cat:eq:lipschitz-exceptional-set} with \(q=2\) and \eqref{cat:eq:constant-ratio-tilt},
\[
\mu\left(\pi^{-1}(E\cap D_{4R})\right)\le C \mathsf H_R^2.
\]
By the area formula, \eqref{cat:eq:constant-ratio-cylinder-bounds}, and \eqref{cat:eq:constant-ratio-tilt}, we have for \(\varphi\in C_c^\infty(D_{4R})\)
\[
\left|\int_{\Omega_0}(u_++u_-)\varphi\mathrm d A -\int_Mu(\varphi\circ\pi)\mathrm d\mu\right| \le C\mathsf H_R^3\|\varphi\|_{L^\infty(D_{4R})}.
\]
Together with \eqref{cat:eq:localized-first-moment-decay} and \(|u_+(Ry)+u_-(Ry)|\le2L(8R)\), this gives for \(0<\rho<1\)
\begin{equation}\label{cat:eq:constant-ratio-midpoint-moment-limit}
\|\eta_\rho*v_R\|_{L^\infty(D_{4-\rho})} \le C_\rho\left(\frac{\mathsf H_R^3}{R^2L(8R)}+\|\xi_R\|_{L^2(D_5)}\right)\longrightarrow0.
\end{equation}

On \(\Omega_0\), \(|\nabla u_\pm|<\sqrt3\). The function \(G(p)=p/\sqrt{1+|p|^2}\) is odd and satisfies \(cI\le dG(p)\le I\) for \(|p|\le\sqrt3\). Hence, for \(|p|,|q|\le\sqrt3\),
\[
(G(p)+G(q))\cdot(p+q) \ge c|p+q|^2,\quad |G(p)+G(q)| \le|p+q|.
\]
For \(\psi\in C_c^\infty(D_4)\), multiply the sum of the minimal surface equations for \(u_\pm(R \cdot)/R\) by \(\psi^2(1-\xi_R)v_R\) and integrate by parts. Using the bound \(|u_+(Ry)+u_-(Ry)|\le2L(8R)\) in the terms containing \(\nabla\xi_R\), we obtain
\begin{equation}\label{cat:eq:constant-ratio-midpoint-caccioppoli}
\int\psi^2|\nabla v_R|^2 \le C\int |v_R|^2|\nabla \psi|^2 +C\int\psi^2|\nabla \xi_R|^2.
\end{equation}
By \eqref{cat:eq:constant-ratio-midpoint-caccioppoli}, \((v_R)\) is bounded in \(W^{1,2}(D_3)\) and hence precompact in \(L^2(D_3)\). Since every limit is zero by \eqref{cat:eq:constant-ratio-midpoint-moment-limit}, \(v_R\to0\) in \(L^2(D_3)\). For some \(\tau_R\in(1,2)\),
\[
\begin{aligned}
\sup_{\partial D_{\tau_R}}(|v_R|+\xi_R) &\le C\|(v_R,\xi_R)\|_{W^{1,2}(D_2)} \\
&\le C\left(\|v_R\|_{L^2(D_3)} +\|\xi_R\|_{W^{1,2}(D_3)}\right)\longrightarrow0.
\end{aligned}
\]
Since \(\xi_R=1\) on \(R^{-1}E\cap D_4\), we have \(\partial D_{\tau_RR}\subset\Omega_0\) for large \(R\). Proposition~\ref{cat:prop:compact-midpoint-principle}, applied to \(M\cap\pi^{-1}(\overline D_{\tau_RR})\), gives
\[
\frac{L(R)}{L(8R)} \le\frac{\sup_{\partial D_{\tau_R}}|v_R|} {1-\sup_{\partial D_{\tau_R}}\xi_R}\longrightarrow0.
\]
For large \(R\),
\[
L(R)\le8^{-n}L(8^nR)\le C8^{-n}\mathsf H_{8^nR}\longrightarrow0\qquad\text{as }n\to\infty,
\]
since \(\mathsf H_s=o(s)\). Hence \(L\equiv0\).

By Proposition~\ref{cat:prop:constant-ratio-capacity}, choose \(t_R\in(1,2)\) such that \(\partial D_{t_RR}\subset\Omega_0\). Since \(L\equiv0\), \(u_++u_-=0\) on \(\partial D_{t_RR}\). Applying Proposition~\ref{cat:prop:compact-midpoint-principle} with \(c=z_\infty\) to \(M\cap\pi^{-1}(\overline D_{t_RR})\) and letting \(R\to\infty\), we find that \(M\) is symmetric about \(u=0\) and that \(M\cap\{\pm u>0\}\) are smooth graphs.

Put \(\Gamma=M\cap\{u=0\}\). Since \(M\) is connected and symmetric about \(u=0\), \(\Gamma\ne\varnothing\). If \(N_z\ne0\) at a point of \(\Gamma\), write \(u=f(x)\) near that point. By symmetry, \(f=-f=0\). Since \(\Delta_Mu=0\), we have \(u\equiv0\) on \(M\) by unique continuation, contrary to the assumption that \(M\) has density two at infinity. Thus \(N_z=0\) on \(\Gamma\), so \(M\cap\{u\ge0\}\) is a smooth surface with boundary \(\Gamma\).

The map \((x,z_\infty+u)\mapsto(x,z_\infty+|u|)\) sends \(M\) onto \(M\cap\{u\ge0\}\), which is therefore connected. Hence \(M^+=M\cap\{u>0\}\) is a graph \(u=q(x)\) over a connected domain \(\Omega\), with \(q>0\). Since \(L\equiv0\), distinct \(X,Y\in M\) with \(\pi X=\pi Y\) satisfy \(u(X)=-u(Y)\ne0\). Since \(\pi|_M\) is proper by Proposition~\ref{cat:prop:constant-ratio-geometry}, \(\pi:M\cap\{u\ge0\}\to\overline\Omega\) is a homeomorphism. Thus \(q\in C^0(\overline\Omega)\) and \(q=0\) on \(\partial\Omega=\pi(\Gamma)\).

Since \(\partial D_{t_RR}\subset\Omega\), \(\partial\Omega\cap D_{t_RR}\) is a compact smooth \(1\)-manifold without boundary, hence a finite disjoint union of Jordan curves. These disks exhaust \(\mathbb R^2\), so \(\partial\Omega\) is a locally finite union of Jordan curves.

Choose the unit normal \(N\) on \(M^+\) so that \(N_z>0\). Along \(\Gamma\), \(N\) points out of \(\Omega\), and \(N_z=0\). The inward unit conormal of \(M^+\) is \(e_3\). For an arclength parametrization \(\gamma\) of \(\partial\Omega\), \(\operatorname{tr}(dN)=0\) and Hopf's lemma for \(\Delta_MN_z=-|dN|^2N_z\) give
\begin{equation}\label{cat:eq:constant-ratio-waist-curvature}
\gamma''\cdot N =-\langle dN(\gamma'),\gamma'\rangle =\langle dN(e_3),e_3\rangle =\partial_{e_3}N_z>0.
\end{equation}
Each boundary component therefore bounds a compact strictly convex domain \(K_i\). Since \(\Omega\) is connected and unbounded, the \(K_i\) are pairwise disjoint and \(\Omega=\mathbb R^2\setminus\bigcup_i K_i\).
\end{proof}

\subsection{The overdetermined problem}\label{cat:sec:planar-rigidity}
\leavevmode\par\nobreak \makeatletter\@afterheading\makeatother

Following Traizet~\cite[Theorem~10]{TraizetRoof}, we construct \(F\) so that \(q\circ F^{-1}\) is harmonic.

\begin{proposition}
\label{cat:prop:constant-ratio-roof-map}
There is a homeomorphism \(F:\mathbb R^2\to\mathbb R^2\) with \(F-\mathrm{id}\) constant on each \(K_i\). Put \(\widetilde K_i=F(K_i)\) and \(\Omega'=F(\Omega)\), and let \(\mathbf{n}\) be the inward unit normal to \(\Omega'\). Then \(F:\Omega\to\Omega'\) is a diffeomorphism, and \(U=q\circ F^{-1}\) satisfies \(U(y)\le C_\varepsilon(1+|y|)^\varepsilon\) for every \(\varepsilon>0\) and
\begin{equation}\label{cat:eq:constant-ratio-roof-problem}
U>0, \quad \Delta U=0, \quad |\nabla U|<1\quad \text{ in }\Omega',\quad U=0, \quad \nabla U=\mathbf{n}\quad \text{ on }\partial\Omega'.
\end{equation}
\end{proposition}
\begin{proof}
Put \(p=\nabla q\), and recall that \(W(p)=\sqrt{1+|p|^2}\) and \(B(p)=W(p)I-p\otimes p/W(p)\). For \(Y\in C_c^\infty(\mathbb R^2;\mathbb R^2)\),
\[
0=\frac12\int_M\operatorname{div}_M(Y(x),0) \mathrm d\mu =\int_\Omega\operatorname{tr}\left(B(p) dY\right)\mathrm d A.
\]
Thus \(\operatorname{div}(\mathbf1_\Omega B(p))=0\) in \(\mathcal D'(\mathbb R^2)\). By the Poincar\'e lemma, there is a convex function \(\Phi\) satisfying
\[
\nabla^2 \Phi=\mathbf1_\Omega\operatorname{cof}B(p) =\mathbf1_\Omega\frac{I+p\otimes p}{W(p)}\ge0.
\]

Near \(\partial K_i\times\{z_\infty\}\), \(M\) has a parametrization
\[
X(s,t)=\left(\gamma(s)+g(s,t)N(s),z_\infty+t\right), \qquad g(s,0)=g_t(s,0)=0,
\]
where \(\gamma\) is an arclength parametrization of \(\partial K_i\) and \(N\) points into \(K_i\). By \eqref{cat:eq:constant-ratio-waist-curvature}, \(g_{tt}(s,0)=-\gamma''(s)\cdot N(s)<0\). Set \(d=-g(s,t)\), so \(x=\gamma(s)-dN(s)\) and \(d=\operatorname{dist}(x,\partial\Omega)\). Since \(g\) is even in \(t\),
\begin{equation}\label{cat:eq:constant-ratio-fold}
q(s,d)=\sqrt d c(s,d), \qquad c(s,0)>0,
\end{equation}
where \(c\) is smooth up to \(d=0\).
Since \(\nabla^2\Phi=0\) in \(K_i^\circ\) and \(\nabla^2\Phi=O(d^{-1/2})\) near \(\partial K_i\), we have \(\Phi\in C^1(\mathbb R^2)\).

We may assume that \(\nabla\Phi(0)=0\). By convexity, the map \(F(x)=x+\nabla\Phi(x)\) satisfies
\begin{equation}\label{cat:eq:constant-ratio-map-control}
(F(x)-F(x'))\cdot(x-x')\ge|x-x'|^2.
\end{equation}
For each \(y\in\mathbb R^2\), the function \(|x|^2/2+\Phi(x)-y\cdot x\) attains its minimum at a point \(x\) satisfying \(F(x)=y\). By \eqref{cat:eq:constant-ratio-map-control}, \(F\) is therefore a homeomorphism with a \(1\)-Lipschitz inverse. Since \(\nabla\Phi\) is constant on each \(K_i\), \(F\) translates \(K_i\) onto \(\widetilde K_i\).

On \(\Omega\),
\begin{equation}\label{cat:eq:constant-ratio-roof-differential}
\begin{aligned}
dF&=I+\frac{I+p\otimes p}{W(p)}\ge I, &dF p&=(1+W(p))p, \\
\det dF&=\frac{(1+W(p))^2}{W(p)}, &\nabla U(F(x))&=(dF)^{-1}p=\frac{p}{1+W(p)}.
\end{aligned}
\end{equation}
For \(\varphi\in C_c^\infty(\Omega')\),
\[
\int_{\Omega'}\nabla U\cdot\nabla\varphi\mathrm d A(y) =\int_\Omega\frac{p}{W(p)}\cdot\nabla(\varphi\circ F)\mathrm d A(x)=0.
\]
Since \(F\) translates each \(\partial K_i\), \eqref{cat:eq:constant-ratio-fold} gives, as \(\Omega\ni x\to x_0\in\partial\Omega\),
\[
\nabla U(F(x))\to \mathbf{n}(F(x_0)).
\]

Since \(F^{-1}\) is \(1\)-Lipschitz and \(F(0)=0\), we have \(|F^{-1}(y)|\le|y|\). By \eqref{cat:eq:constant-ratio-cylinder-bounds},
\[
U(y)=q(F^{-1}(y))\le C_\varepsilon(1+|y|)^\varepsilon.
\]
\end{proof}

Put \(U=0\) on \((\Omega')^c\) and \(\lambda=\mathcal H^1\llcorner\partial\Omega'\).

\begin{proposition}
\label{cat:prop:constant-ratio-roof-field}
\(\Delta U=\lambda\) in \(\mathcal D'(\mathbb R^2)\) and
\begin{equation}\label{cat:eq:constant-ratio-charge-bounds}
\lambda(D_R)+R\int_{|y|>R}\frac{\mathrm d\lambda(y)}{|y|} \le C_\varepsilon R^\varepsilon
\end{equation}
for \(R\ge1\) and \(0<\varepsilon<1\). For \(x\notin\partial\Omega'\),
\begin{equation}\label{cat:eq:constant-ratio-roof-field}
\nabla U(x)=\frac1{2\pi}\int \frac{x-y}{|x-y|^2} \mathrm{d}\lambda(y).
\end{equation}
Locally uniformly on \(\partial\Omega'\),
\begin{equation}\label{cat:eq:constant-ratio-field-trace}
\lim_{\delta\downarrow0}\int_{|x-y|>\delta} \frac{x-y}{|x-y|^2} \mathrm{d}\lambda(y) =\pi \mathbf{n}(x).
\end{equation}
\end{proposition}
\begin{proof}
Integration by parts, with \(\nabla U=\mathbf{n}\) on \(\partial\Omega'\), gives \(\Delta U=\lambda\ge0\). Hence
\[
\lambda(D_R) \le\frac{C}{R^2}\int_{D_{2R}}U \le C_\varepsilon R^\varepsilon,\quad R\int_{|y|>R}\frac{\mathrm d\lambda(y)}{|y|} \le R\int_R^\infty\frac{\lambda(D_t)}{t^2} \mathrm dt \le C_\varepsilon R^\varepsilon.
\]

By \eqref{cat:eq:constant-ratio-charge-bounds} and the growth bound for \(U\),
\[
R^{-1}\left(\sup_{D_{2R}}U+\lambda(D_{2R})\right)\longrightarrow0.
\]
For \(x\in\mathbb R^2\), applying Green's formula to \(U\) and \((x-y)(|x-y|^{-2}-R^{-2})\) on \(\{\delta<|y-x|<R\}\), and letting \(R\to\infty\), we obtain
\[
\int_{|x-y|>\delta}\frac{x-y}{|x-y|^2}\mathrm d\lambda(y)=\int_{S^1}e\left(e\cdot\nabla U(x+\delta e)+\frac{U(x+\delta e)-U(x)}{\delta}\right)\mathrm d e.
\]
For \(x\in\partial\Omega'\) and \(e\cdot \mathbf{n}(x)\ne0\), \eqref{cat:eq:constant-ratio-roof-problem} and \(U=0\) on \((\Omega')^c\) give
\[
\frac{U(x+\delta e)}{\delta}\longrightarrow\left(e\cdot \mathbf{n}(x)\right)_+,\qquad \nabla U(x+\delta e)\longrightarrow \mathbf{n}(x)\mathbf1_{\{e\cdot \mathbf{n}(x)>0\}}.
\]
Thus
\[
\lim_{\delta\downarrow0}\int_{|x-y|>\delta}\frac{x-y}{|x-y|^2}\mathrm d\lambda(y)=\begin{cases}2\pi\nabla U(x),&x\notin\partial\Omega',\\ \pi \mathbf{n}(x),&x\in\partial\Omega'.\end{cases}
\]
The convergence on \(\partial\Omega'\) is locally uniform by \eqref{cat:eq:constant-ratio-roof-problem}.
\end{proof}

\begin{lemma}
\label{cat:lem:constant-ratio-sparse-line}
If \((\Omega')^c\) is unbounded, there is an affine line in \(\Omega'\) separating two points of \((\Omega')^c\).
\end{lemma}
\begin{proof}
Suppose otherwise. We may assume that \(0\in(\Omega')^c\). Since \((\Omega')^c\) is unbounded, for each \(R\ge1\) there is \(y_R\in(\Omega')^c\) with \(|y_R|>2R\). If \(e\in S^1\) and \(e\cdot y_R>|y_R|/2\), every line \(e\cdot y=t\), \(0<t<R\), separates \(0\) and \(y_R\) and must meet \(\partial\Omega'\). By the coarea formula and \eqref{cat:eq:constant-ratio-charge-bounds},
\[
\begin{aligned}[b]
R &\le C\int_{S^1}\int_0^R \#\left(\partial\Omega'\cap\{e\cdot y=t\}\right)\mathrm d t\mathrm d e \\
&\le C\int \left|\{e\in S^1:0<e\cdot y<R\}\right|\mathrm d\lambda(y) \\
&\le C\lambda(D_{2R}) +CR\int_{|y|>2R}\frac{\mathrm d\lambda(y)}{|y|} =o(R).
\end{aligned}\qedhere
\]
\end{proof}

\begin{proposition}
\label{cat:prop:planar-rigidity}
After a translation,
\[
\Omega'=\{|y|>b\}, \qquad U(y)=b\log(|y|/b)
\]
for some \(b>0\).
\end{proposition}
\begin{proof}
Suppose that \((\Omega')^c\) is unbounded. By Lemma~\ref{cat:lem:constant-ratio-sparse-line}, we may choose coordinates such that \(\{y_1=0\}\subset\Omega'\) and
\[
\Gamma_\pm:=\partial\Omega'\cap\{\pm y_1>0\}\ne\varnothing.
\]
Let \(\eta_R(y)=\eta(|y|/R)\), where \(\eta\in C^\infty([0,\infty))\), \(\eta'\le0\), \(\eta=1\) on \([0,1]\), and \(\eta=0\) on \([2,\infty)\). By the divergence theorem,
\[
\int_{\Gamma_+}\eta_R\mathbf{n}\cdot e_1 \mathrm d\lambda =\sum_{\widetilde K_i\subset\{y_1>0\}} \int_{\widetilde K_i}\partial_1\eta_R\mathrm d A\le0.
\]
For \(x\in D_{2R}\) and \(y\ne x\),
\[
\frac{|\eta_R(x)-\eta_R(y)|}{|x-y|}\le\frac{C}{R+|y|}.
\]
Fix \(0<\varepsilon<1/2\). By \eqref{cat:eq:constant-ratio-field-trace} and \eqref{cat:eq:constant-ratio-charge-bounds},
\begin{equation}\label{cat:eq:constant-ratio-force-tail-error}
\begin{aligned}
&\iint_{\Gamma_+\times\Gamma_-} \eta_R(x)\frac{x_1-y_1}{|x-y|^2} \mathrm d\lambda(x) \mathrm d\lambda(y) \\
\le{}& -\frac12\iint_{\Gamma_+\times\Gamma_+} (\eta_R(x)-\eta_R(y))\frac{x_1-y_1}{|x-y|^2} \mathrm d\lambda(x) \mathrm d\lambda(y) \\
\le{}& C\lambda(D_{2R})\int\frac{\mathrm d\lambda(y)}{R+|y|} \le C_\varepsilon R^{2\varepsilon-1}\longrightarrow0.
\end{aligned}
\end{equation}

The integral over \(\Gamma_+\times\Gamma_-\) is positive for large \(R\) and nondecreasing in \(R\), a contradiction. Thus \(K=(\Omega')^c\) is bounded and is a finite union of the \(\widetilde K_i\), by local finiteness.

Since \(x\cdot(x-y)+y\cdot(y-x)=|x-y|^2\), \eqref{cat:eq:constant-ratio-field-trace} and the divergence theorem give
\begin{equation}\label{cat:eq:constant-ratio-finite-virial}
2\pi|K| =\lim_{\delta\downarrow0}\iint_{|x-y|>\delta} \frac{x\cdot(x-y)}{|x-y|^2} \mathrm d\lambda(x) \mathrm d\lambda(y) =\tfrac12\lambda(\partial K)^2\ge2\pi|K|.
\end{equation}
Equality in the isoperimetric inequality implies that \(K\) is a disk. After a translation, \(\Omega'=\{|y|>b\}\) for some \(b>0\). By \eqref{cat:eq:constant-ratio-roof-field} and \(U|_{\partial\Omega'}=0\),
\[
\nabla U(y)=\frac{by}{|y|^2}, \qquad U(y)=b\log(|y|/b).\qedhere
\]
\end{proof}

\begin{proof}[Proof of Theorem~\ref{thm:catenoid-classification}]
By Proposition~\ref{cat:prop:planar-rigidity}, we may assume that \(\Omega'=\{|y|>b\}\). By \eqref{cat:eq:constant-ratio-roof-differential}, with \(r=|y|>b\),
\[
d(F^{-1})=\tfrac12(1+|\nabla U|^2)I-\nabla U\otimes \nabla U =d_y\left[\tfrac12(1+b^2/r^2)y\right].
\]
Hence, for some \(x_0\in\mathbb R^2\),
\[
|F^{-1}(y)-x_0|=\tfrac12(r+b^2/r)=b\cosh(U(y)/b).\qedhere
\]
\end{proof}

\begin{proof}[Proof of Theorem~\ref{thm:classification}]
This follows from Lemma~\ref{lem:cones} and Theorems~\ref{thm:transverse-scherk-classification} and~\ref{thm:catenoid-classification}.
\end{proof}


\begin{thebibliography}{99}

\bibitem{Allard}
W.~K.~Allard, \href{https://doi.org/10.2307/1970868}{\emph{On the first variation of a varifold}}, Ann. of Math. (2) \textbf{95} (1972), 417--491.

\bibitem{AllardAlmgren}
W.~K.~Allard and F.~J.~Almgren, Jr., \href{https://doi.org/10.2307/2006984}{\emph{On the radial behavior of minimal surfaces and the uniqueness of their tangent cones}}, Ann. of Math. (2) \textbf{113} (1981), no.~2, 215--265.

\bibitem{ABR}
S.~Axler, P.~Bourdon, and W.~Ramey, \href{https://www.axler.net/HFT.pdf}{\emph{Harmonic Function Theory}}, second edition, Graduate Texts in Mathematics, vol.~137, Springer, New York, 2001.

\bibitem{BeckerKahn}
S.~T.~Becker-Kahn, \href{https://doi.org/10.4310/jdg/1506650421}{\emph{Transverse singularities of minimal two-valued graphs in arbitrary codimension}}, J. Differential Geom. \textbf{107} (2017), no.~2, 241--325.

\bibitem{BDD}
C.~Brena, S.~Decio, and C.~De Lellis, \href{https://doi.org/10.4171/RLM/1077}{\emph{Remarks and conjectures on stationary varifolds dedicated to Enrico Bombieri}}, Rend. Lincei Mat. Appl. \textbf{36} (2025), no.~3, 409--443.

\bibitem{Collin}
P.~Collin, \href{https://doi.org/10.2307/2951822}{\emph{Topologie et courbure des surfaces minimales proprement plong\'ees de \(\mathbb R^3\)}}, Ann. of Math. (2) \textbf{145} (1997), no.~1, 1--31.

\bibitem{Delort}
J.-M.~Delort, \href{https://www.numdam.org/item/SEDP_1990-1991____A2_0.pdf}{\emph{Existence de nappes de tourbillon pour l'\'equation d'Euler sur le plan}}, S\'eminaire \'Equations aux d\'eriv\'ees partielles, 1990--1991, expos\'e II, 1--12.

\bibitem{DiBenedetto}
E.~DiBenedetto, \href{https://doi.org/10.1007/978-1-4939-4005-9}{\emph{Real Analysis}}, second edition, Birkh\"auser Advanced Texts: Basler Lehrb\"ucher, Birkh\"auser/Springer, New York, 2016.

\bibitem{DS}
S.~Donaldson and S.~Sun, \href{https://doi.org/10.4310/jdg/1506650422}{\emph{Gromov--Hausdorff limits of K\"ahler manifolds and algebraic geometry, II}}, J. Differential Geom. \textbf{107} (2017), no.~2, 327--371.

\bibitem{HHW}
R.~Haslhofer, O.~Hershkovits, and B.~White, \href{https://doi.org/10.1353/ajm.2023.a902954}{\emph{Moving plane method for varifolds and applications}}, Amer. J. Math. \textbf{145} (2023), no.~4, 1051--1076.

\bibitem{LazardHollyMeeks}
H.~Lazard-Holly and W.~H.~Meeks III, \href{https://doi.org/10.1007/PL00005796}{\emph{Classification of doubly-periodic minimal surfaces of genus zero}}, Invent. Math. \textbf{143} (2001), no.~1, 1--27.

\bibitem{LSS}
J.~D.~Lotay, F.~Schulze, and G.~Sz\'ekelyhidi, \href{https://doi.org/10.1007/s10240-023-00143-5}{\emph{Ancient solutions and translators of Lagrangian mean curvature flow}}, Publ. Math. Inst. Hautes \'Etudes Sci. \textbf{140} (2024), 1--35.

\bibitem{Luo}
H.-P.~Luo, \emph{Desingularizing the intersection between a catenoid and a plane}, Ph.D. thesis, University of Massachusetts Amherst, 1997.

\bibitem{Meeks2003}
W.~H.~Meeks III, \href{https://doi.org/10.4310/SDG.2003.v8.n1.a10}{\emph{Geometric results in classical minimal surface theory}}, Surv. Differ. Geom. \textbf{8} (2003), 269--306.

\bibitem{MeeksRosenberg}
W.~H.~Meeks III and H.~Rosenberg, \href{https://doi.org/10.1007/BF01389046}{\emph{The global theory of doubly periodic minimal surfaces}}, Invent. Math. \textbf{97} (1989), no.~2, 351--379.

\bibitem{MW}
W.~H.~Meeks III and M.~Wolf, \href{https://doi.org/10.1090/S0894-0347-06-00537-6}{\emph{Minimal surfaces with the area growth of two planes: the case of infinite symmetry}}, J. Amer. Math. Soc. \textbf{20} (2007), no.~2, 441--465.

\bibitem{MichaelSimon}
J.~H.~Michael and L.~M.~Simon, \href{https://doi.org/10.1002/cpa.3160260305}{\emph{Sobolev and mean-value inequalities on generalized submanifolds of \(\mathbb R^n\)}}, Comm. Pure Appl. Math. \textbf{26} (1973), no.~3, 361--379.

\bibitem{Necas}
J.~Ne\v{c}as, \href{https://doi.org/10.1007/978-3-642-10455-8}{\emph{Direct Methods in the Theory of Elliptic Equations}}, Springer Monographs in Mathematics, Springer, Heidelberg, 2012.

\bibitem{PerezICM}
J.~P\'erez, \href{https://doi.org/10.1137/25M1803188}{\emph{Minimal surfaces of finite genus: classification, dynamics and laminations}}, in \emph{Proceedings of the International Congress of Mathematicians 2026}, vol.~4, 2026, 215--232.

\bibitem{PerezTraizet}
J.~P\'erez and M.~Traizet, \href{https://doi.org/10.1090/S0002-9947-06-04094-3}{\emph{The classification of singly periodic minimal surfaces with genus zero and Scherk-type ends}}, Trans. Amer. Math. Soc. \textbf{359} (2007), no.~3, 965--990.

\bibitem{SchoenSymmetry}
R.~M.~Schoen, \href{https://doi.org/10.4310/jdg/1214438183}{\emph{Uniqueness, symmetry, and embeddedness of minimal surfaces}}, J. Differential Geom. \textbf{18} (1983), 791--809.

\bibitem{SimonAsymptotics}
L.~Simon, \href{https://doi.org/10.2307/2006981}{\emph{Asymptotics for a class of non-linear evolution equations, with applications to geometric problems}}, Ann. of Math. (2) \textbf{118} (1983), no.~3, 525--571.

\bibitem{SimonCylindrical}
L.~Simon, \href{https://doi.org/10.4310/jdg/1214454484}{\emph{Cylindrical tangent cones and the singular set of minimal submanifolds}}, J. Differential Geom. \textbf{38} (1993), no.~3, 585--652.

\bibitem{SimonGMT}
L.~Simon, \href{https://web.stanford.edu/class/math285/ts-gmt.pdf}{\emph{Introduction to Geometric Measure Theory}}, Tsinghua lecture notes, 2014.

\bibitem{TraizetRoof}
M.~Traizet, \href{https://doi.org/10.1007/s00039-014-0268-5}{\emph{Classification of the solutions to an overdetermined elliptic problem in the plane}}, Geom. Funct. Anal. \textbf{24} (2014), no.~2, 690--720.

\bibitem{White}
B.~White, \href{https://doi.org/10.1215/00127094-2009-019}{\emph{Currents and flat chains associated to varifolds, with an application to mean curvature flow}}, Duke Math. J. \textbf{148} (2009), no.~1, 41--62.

\bibitem{Zhou}
J.~Zhou, \href{https://doi.org/10.1093/imrn/rnaa340}{\emph{Topology of surfaces with finite Willmore energy}}, Int. Math. Res. Not. IMRN (2022), no.~9, 7100--7151.

\end{thebibliography}
\end{document}